\documentclass[mnsc,nonblindrev]{informs3}
\usepackage{comment}
\usepackage{algorithm, algorithmic}
\usepackage[font=small,skip=0pt]{caption}
\usepackage{subcaption}
\usepackage{multirow}
\usepackage{enumitem}
\usepackage{graphicx}
\usepackage{hyperref}
\usepackage[toc,page]{appendix}
\usepackage{tikz}
\usepackage{pgfplots}
\pgfplotsset{compat=1.17}
\OneAndAHalfSpacedXI 

\usepackage{natbib}
 \bibpunct[, ]{(}{)}{,}{a}{}{,}%
 \def\bibfont{\small}%
 \def\newblock{\ }%
\TheoremsNumberedThrough     
\ECRepeatTheorems

\EquationsNumberedThrough    

\begin{document}



\RUNTITLE{Risk-Averse Stochastic Facility Location with or without Prioritization}

\TITLE{On the Value of Multistage Risk-Averse Stochastic Facility Location with or without Prioritization}

\ARTICLEAUTHORS{%
\AUTHOR{Xian Yu}
\AFF{Department of Industrial and Operations Engineering, University of Michigan at Ann Arbor, USA; \EMAIL{yuxian@umich.edu}} 
\AUTHOR{Siqian Shen}
\AFF{Corresponding author; Department of Industrial and Operations Engineering, University of Michigan at Ann Arbor, USA; \EMAIL{siqian@umich.edu}
}
} 

\ABSTRACT{%
We consider a multiperiod stochastic capacitated facility location problem under uncertain demand and budget in each period. Using a scenario tree representation of the uncertainties, we formulate a multistage stochastic integer program to dynamically locate facilities in each period and compare it with a two-stage approach that determines the facility locations up front. In the multistage model, in each stage, a decision maker optimizes facility locations and recourse flows from open facilities to demand sites, to minimize certain risk measures of the cost associated with current facility location and shipment decisions. When the budget is also uncertain, a popular modeling framework is to prioritize the candidate sites \citep{kocc2015prioritization}. In the two-stage model, the priority list is decided in advance and fixed through all periods, while in the multistage model, the priority list can change adaptively. In each period, the decision maker follows the priority list to open facilities according to the realized budget, and optimizes recourse flows given the realized demand. Using expected conditional risk measures (ECRMs), we derive tight lower bounds for the gaps between the optimal objective values of risk-averse multistage models and their two-stage counterparts in both settings with and without prioritization. Moreover, we propose two approximation algorithms to efficiently solve risk-averse two-stage and multistage models without prioritization, which are asymptotically optimal under an expanding market assumption. We also design a set of super-valid inequalities for risk-averse two-stage and multistage stochastic programs with prioritization to reduce the computational time. We conduct numerical studies using both randomly generated and real-world instances with diverse sizes, to demonstrate the tightness of the analytical bounds and efficacy of the approximation algorithms and prioritization cuts. We find that the gaps between risk-averse multistage and two-stage models increase as the variations of the uncertain parameters increase, and stagewise dependent scenario trees attain much higher gaps than the stagewise independent ones. 
}%



\maketitle

\section{Introduction}
\label{sec:intro}
The capacitated facility location problem, aiming to build facilities in potential locations to meet customers' demand, is one of the most classical optimization problems solved in a broad spectrum of applications, including locating warehouses in supply chains \citep{perl1985warehouse, aghezzaf2005capacity}, shelters in disaster relief networks \citep{rawls2010pre, balcik2008facility}, car rental facilities \citep{lin2011strategic, garcia2012optimizing}, 
and so on. In these applications, demand fluctuates spatially and temporally, resulting in uncertain operational cost of using the facilities over time. To adapt to new demand, service providers need to relocate facilities or expand existing capacities, which could be costly and inefficient. Moreover, the budget for opening facilities may be uncertain and paid by installment per period. Therefore, estimating and utilizing uncertain demand and budget information in the decision processes of locating facilities is crucial for the purposes of cost reduction and quality-of-service improvement.

In this paper, we focus on a finite time horizon (e.g., months or years) {of planning the locations of facilities}, where customers' demand {and budget for opening facilities} in each period are modeled by a {joint} random vector. {When only the demand is uncertain,} a decision maker optimizes when and where to open facilities and how to supply products from open facilities to customers, to minimize a certain risk measure of the cost of locating facilities and shipping products over multiple periods. We compare two modeling frameworks: a two-stage stochastic optimization model and a multistage stochastic dynamic program. In the two-stage model, we determine facility locations for all periods at the beginning of the time horizon as ``here-and-now'' decisions and decide optimal shipping assignments as ``wait-and-see'' recourse for each sample of the multiperiod demand. In the multistage model, the uncertain demand values are revealed gradually over time and the facility-location decisions are adapted to this process. Specifically, we determine the current period's facility opening and product shipping plans, given demand information revealed up to each current period. 

\paragraph{Prioritization for Facility Location.} In practice, when the budget for opening facilities is also uncertain,  practitioners may prefer a rank-ordered list of solutions (e.g., candidate facility sites), termed as \textit{priority list}, and choose those that have higher priority until using up all the current budget \citep{kocc2015prioritization}. The prioritization processes include (i) placing candidate sites into a priority list before the uncertainty is revealed and, (ii) after realizing the values of uncertain parameters, the site-selection decisions must be consistent with the priority list. In other words, prioritization requires the selection of candidate sites to be \textit{nested} with respect to different budget values. There are several appealing properties of using prioritization. First, because the optimal facility locations are nested, we can incrementally locate facilities when the budget is not pre-specified. Note that when the budget changes, we may need to relocate some existing facilities if the optimal facility locations are not nested. Second, the priority list constructs an easy-to-implement policy, where it always remains feasible under different realizations of the uncertain demand and budget as long as the budget is set large enough to cover all the demand in each scenario. 
In a two-stage setting, the priority list is decided in advance and fixed through all periods, while in a multistage setting, the priority list can change adaptively. To the best of our knowledge, this paper is the first to propose a multistage stochastic dynamic programming framework to model facility location with prioritization under uncertain budget and demand, and derive analytical bounds by comparing it with a two-stage counterpart.

\paragraph{Applications of Two-Stage and Multistage Facility Location.} Both two-stage and multistage decision frameworks are commonly used in the stochastic facility location literature and its wide applications. For example, power system operators need to build new transmission lines over {multiple periods of years to satisfy} growing demand {for electricity} \citep{bruno2016risk}. A two-stage model can be solved to determine locations of transmission lines for each year up front, given forecasted demand. Alternatively, at the end of each year, the operators can decide new lines to build or expand, as well as power generation. {We can also rank the candidate transmission lines to form a priority list, which can be followed and executed more easily by power system operators when they are given specific budget for each year.} Another example  is in shared-mobility market penetration, where future carsharing demand is obscure and can fluctuate dependent on technology maturity and public acceptance \citep{lu2018optimizing,zhang2021values}. In such a case, there exists little market information that can be used for accurate demand forecast, and decision makers may choose to take multiple stages to locate car rental facilities, as well as charging stations for electric vehicles, rather than commit all the resources up front without accurate demand information. Note that in the context of facility location with prioritization, a two-stage model will first decide the priority list and then adapt resource allocation for specific demand and budget realizations in the second stage; a multistage model may have more flexibility to update the priority list and locate facilities after observing the randomness in each period so as to avoid relocation and demolition cost.

A natural question is then about the performance of the above two decision-making frameworks. Specifically, in this paper, we are interested in comparing the objective values and computational effort between solving two-stage and multistage models {for stochastic facility location with or without prioritization}. Noting that the multistage models have {larger feasible regions} and thus will always have better cost-wise objective values, we aim to bound the gap between optimal objective values of the two-stage and multistage models given specific risk measures of the cost and characteristics of the uncertainty. \cite{huang2009value} are the first to show analytical bounds for the value of multistage stochastic programming (VMS) compared to the two-stage approach for capacity planning problems with an \textit{expectation-based} objective function. \citet{maggioni2012analyzing, birge2011introduction} propose the concept of the  Value of Stochastic Solution (VSS) and present bounds on the potential benefit from solving a two-stage stochastic program over a deterministic counterpart based on the mean values of uncertain parameters. \citet{maggioni2014bounds, escudero2007value, nickel2012multi} extend the measure of uncertain information from two-stage to multistage stochastic programs, compared with their deterministic counterparts.  To our best knowledge, it remains an open question to bound the gap between two-stage and multistage facility location models, if using \textit{risk-averse} objective functions, {or under a prioritization setting}. In this paper, we consider a class of coherent risk measures (i.e., ECRMs) and provide \textit{tight} lower bounds of the gap between the optimal objective values of these two models {with or without prioritization}.

As the two-stage and multistage stochastic mixed-integer programs are known to be computationally intractable, we also develop approximation algorithms for solving the two models without prioritization and a set of cutting planes for solving the two models with prioritization. It turns out that the approximation schemes are asymptotically optimal under increasing demand (e.g., when the market is expanding or launching new businesses). When the budget is also uncertain, we propose a set of super-valid inequalities, which may rule out some feasible solutions but ensure that at least one optimal solution remains.

The main contributions of this work are summarized in Table \ref{tab:contribution} below.
\begin{table}[ht!]
  \centering
  \caption{Main contributions of this work}
  \resizebox{\textwidth}{!}{
    \begin{tabular}{lll}
    \hline
          & Without Prioritization & With Prioritization \\
          \hline
    Uncertainty & Demand & Demand and budget\\
    Implemented decisions & Facility locations & Priority lists\\
    VMS   &   ${\rm VMS_R}=z_R^{TS}-z_R^{MS}$    & ${\rm VMS_P}=z_P^{TS}-z_P^{MS}$ \\
    \multirow{3}[0]{*}{${\rm VMS^{LB}}$} &   \textit{Tight} lower bound ${\rm VMS_R^{LB}}$ (Theorem \ref{thm:risk-averse-bound})    &  \textit{Tight} lower bound  ${\rm VMS_P^{LB}}$ (Theorem \ref{thm:risk-averse-bound-prior})\\
          &   ${\rm VMS_R^{LB1}}$ relying on LP relaxations (Corollary \ref{cor:LP})    & ${\rm VMS_P^{LB1}}$ relying on parameters (Remark \ref{remark:parameter}) \\
          &  ${\rm VMS_R^{LB2}}$ relying on parameters (Corollary \ref{cor:parameter})     &  \\
    Computation & Approximation algorithms (Section \ref{sec:approx}) & Prioritization cuts (Section \ref{sec:prioritization-cuts}) \\
    Out-of-sample test & Not applicable & Rolling horizon approach (Section \ref{sec:risk-profile}) \\
    \hline
    \end{tabular}%
    }
  \label{tab:contribution}%
\end{table}%

    The remainder of the paper is organized as follows. 
    In Section \ref{sec:value-proof}, tight lower bounds are derived for the gaps between the optimal ECRM-based objective values of the two models {without prioritization}. We also propose approximation algorithms with performance guarantee.
    In Section \ref{sec:VMS-prior}, we provide tight lower bounds for the gaps between the two models with prioritization, and derive a set of cutting planes to speed up the computation.
    In Section \ref{sec:compu}, numerical studies are conducted on instances with diverse uncertainty patterns, to show the tightness of our derived bounds, as well as performance of the approximation algorithms and prioritization cuts. 
Section \ref{sec:conclu} concludes the paper and states future research directions. We review the most relevant papers and clarify specific contributions of our paper compared to the existing literature in Appendix \ref{sec:literature}.

\paragraph{Notation.} Throughout the paper, we use bold symbols to denote vectors/matrices and use $[n]$ to denote the set $\{1,2,\ldots,n\}$.
	
	\section{Value of {Risk-Averse} Multistage {Facility Location}}
	\label{sec:value-proof}
		
    {We first consider the case where only the multiperiod demand is uncertain.
    We describe the problem formulations in Section \ref{sec:models}}, and then examine a substructure of our problem and provide its analytical optimal solutions in Section \ref{sec:substructure}. We derive lower bounds for the gaps between the objective values of the risk-averse two-stage and multistage models with ECRM-based objective in Section \ref{sec:VMS}. {We further design approximation algorithms with performance guarantee for the risk-averse two-stage and multistage facility location problems in Section \ref{sec:approx}.}
    
	\subsection{Problem Formulations}
	\label{sec:models}
	
Consider $1,\ldots,T$ stages, $1,\ldots,M$ potential facility locations, and $1,\ldots,N$ customer sites. Let $c_{tij}$ be the cost of shipping one unit of product from facility $i$ to customer $j$ at stage $t$, $f_{ti}$ be the fixed cost of renting a facility in location $i$ {at stage $t$}, $h_{ti}$ be the capacity of facility $i$ {at stage $t$}, and $d_{tj}$ be the demand at customer site $j$ at stage $t$ for all $i\in [M]$, $j \in [N]$, and $t \in [T]$. 

Define decision variables $x_{ti} \in \{0,1\}$ such that $x_{ti} =1$ if we start to {open} facility $i$ at the beginning of stage $t$ (we assume that the facility will remain open until the last stage $T$), and $x_{ti} =0$ otherwise. We also define continuous variables $y_{tij}$ as the amount of flow we ship from facility $i$ to customer $j$ at stage $t$. For a deterministic capacitated facility location problem, the goal is to minimize the total cost of renting facilities and flow, subject to satisfying demand at each customer site in each stage. The mixed-integer programming model is given by:
\begin{subequations}
\label{model:deter}
\begin{align}
\min\quad &\sum_{t=1}^T\sum_{i=1}^M f_{ti}\sum_{\tau=1}^tx_{\tau i}+\sum_{t=1}^T\sum_{i=1}^M\sum_{j=1}^N c_{tij}y_{tij} \label{eq:1}\\
\text{s.t.}\quad& \sum_{i=1}^M y_{tij}= d_{tj}, \ \forall j=1,\ldots,N,\ t=1,\ldots,T \label{eq:2}\\
&\sum_{j=1}^N y_{tij}\le h_{ti}\sum_{\tau=1}^t x_{\tau i},\ \forall i=1,\ldots,M,\ t=1,\ldots,T \label{eq:3}\\
&\sum_{\tau=1}^t x_{\tau i}\le 1, \ \forall i=1,\ldots,M,
 \ t=1,\ldots,T\label{eq:4}\\
&x_{ti}\in \mathbb{Z}_{+}, \ \forall i=1,\ldots,M,\ t=1,\ldots,T\label{eq:5}\\
&y_{tij}\in \mathbb{R}_{+}, \ \forall i=1,\ldots,M,\ j=1,\ldots,N,\ t=1,\ldots,T.\label{eq:6}
\end{align}
\end{subequations}
The objective function \eqref{eq:1} minimizes the total rental cost and operational cost, where $\sum_{\tau=1}^tx_{\tau i}$ represents whether a facility $i$ is open at stage $t$ and the cost is summed over all open facilities. Constraints \eqref{eq:2} require all the demand to be satisfied in each stage. Constraints \eqref{eq:3} indicate that in each period, we can only ship from a facility within its capacity when it is open. Constraints \eqref{eq:4} imply that we cannot {open} more than once in the same location. Because of constraints \eqref{eq:4}, we can relax binary variable $x_{ti}$ to be integer-valued, indicated in \eqref{eq:5}. 

\begin{remark}
In practice, one can also allow closing facilities in later periods. To accommodate this flexibility, we can modify the decision variable $x_{ti}\in\{0,1\}$ such that $x_{ti}=1$ if the facility $i$ is open in stage $t$, and Model \eqref{model:deter} can be revised as follows:
\begin{align*}
\min\quad &\sum_{t=1}^T\sum_{i=1}^M f_{ti}x_{ti}+\sum_{t=1}^T\sum_{i=1}^M\sum_{j=1}^N c_{tij}y_{tij}\\
\text{s.t.}\quad& \sum_{i=1}^M y_{tij}= d_{tj}, \ \forall j=1,\ldots,N,\ t=1,\ldots,T \\
&\sum_{j=1}^N y_{tij}\le h_{ti}x_{ti},\ \forall i=1,\ldots,M,\ t=1,\ldots,T \\
&x_{ti}\in \{0,1\}, \ \forall i=1,\ldots,M,\ t=1,\ldots,T\\
&y_{tij}\in \mathbb{R}_{+}, \ \forall i=1,\ldots,M,\ j=1,\ldots,N,\ t=1,\ldots,T.
\end{align*}
Under this setting, we can also derive similar lower bounds as we show later for the gaps between risk-averse two-stage and multistage models.
However, frequently closing facilities may lead to practical inconvenience, and as a result, we focus on Model \eqref{model:deter} in the analysis of this paper.
\end{remark}

Model \eqref{model:deter} can be rewritten in a vector form below: 
\begin{subequations}\label{model:vector}
\begin{align}
\min_{\substack{\boldsymbol{x}_1,\ldots,\boldsymbol{x}_T\\\boldsymbol{y}_1,\ldots,\boldsymbol{y}_T}}\quad &\sum_{t=1}^T \boldsymbol{f}_{t}^{\mathsf T} \sum_{\tau=1}^t \boldsymbol{x}_{\tau}+\sum_{t=1}^T\boldsymbol{c}_t^{\mathsf T} \boldsymbol{y}_t\\
\text{s.t.}\quad& \boldsymbol{A}\boldsymbol{y}_t= \boldsymbol{d}_t, \ \forall t=1,\ldots,T \label{eq:2-2}\\
&\boldsymbol{B}_{t}\boldsymbol{y}_t\le \sum_{\tau=1}^t \boldsymbol{x}_{\tau},\ \forall t=1,\ldots,T \label{eq:2-3}\\
& \sum_{\tau=1}^t \boldsymbol{x}_{\tau}\le \mathbf{1}, \ \forall t=1,\ldots,T\label{eq:2-4}\\
&\boldsymbol{x}_t\in \mathbb{Z}^{M}_{+}, \ \boldsymbol{y}_t\in \mathbb{R}^{M\times N}_{+},\ \forall t=1,\ldots,T,\nonumber
\end{align}
\end{subequations}
where matrices $\boldsymbol{A}\in \mathbb{R}^{M},\ \boldsymbol{B}_{t}\in\mathbb{R}^N$ correspond to the coefficients of constraints \eqref{eq:2} and \eqref{eq:3}, respectively. 
In Model \eqref{model:vector}, the data we acquire at each stage $t$ is the demand $\boldsymbol{d}_t\in\mathbb{R}_+^N$ for all $t=1,\ldots, T$. We consider that the data series $\{\boldsymbol{d}_2,\ldots,\boldsymbol{d}_T\}$ evolve according to a known probability distribution, and $\boldsymbol{d}_1\in \mathbb R^N_+$ is deterministic (see similar assumptions made in, e.g.,  \cite{huang2009value,shapiro2009lectures,zou2019stochastic}). In practice, $\boldsymbol{d}_1$ can be derived and forecasted as the average of historical demand, based on which we make the first-stage decisions. {Note that in Section \ref{sec:VMS-prior}, we will consider the case where $\{\boldsymbol{d}_1,\boldsymbol{d}_2,\ldots,\boldsymbol{d}_T\}$ are uncertain, and in Stage $0$, an initial priority list needs to be decided before locating facilities and planning shipments in Stage 1.}

\begin{figure}[ht!]
    \centering
    \begin{subfigure}{0.45\textwidth}
    \centering
    		\includegraphics[height = 6.5cm]{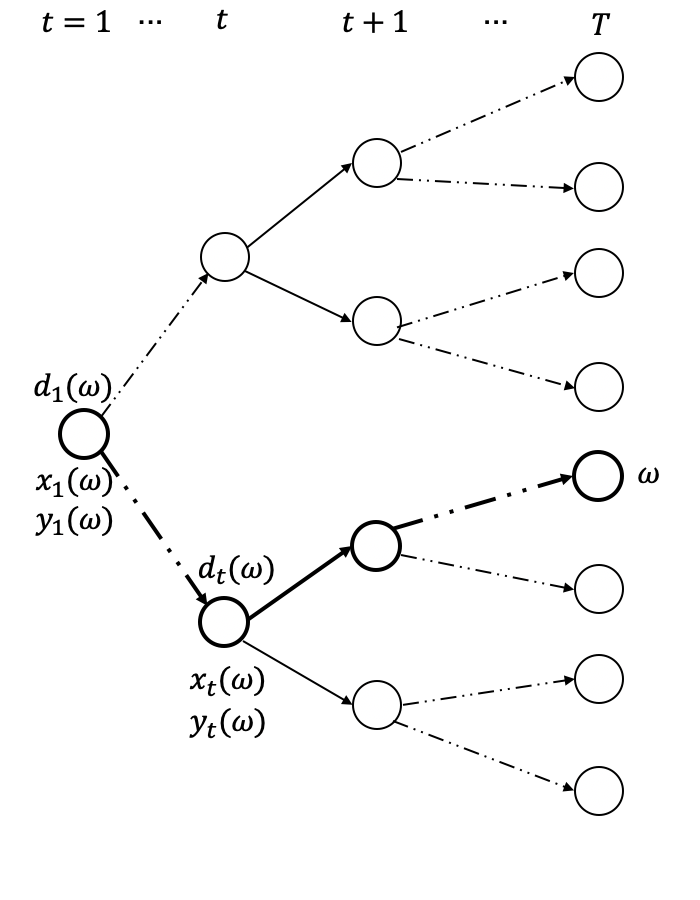}
        \caption{Scenario-path-based notation}
    \end{subfigure}
	\begin{subfigure}{0.45\textwidth}
	\centering
			\includegraphics[height = 6.5cm]{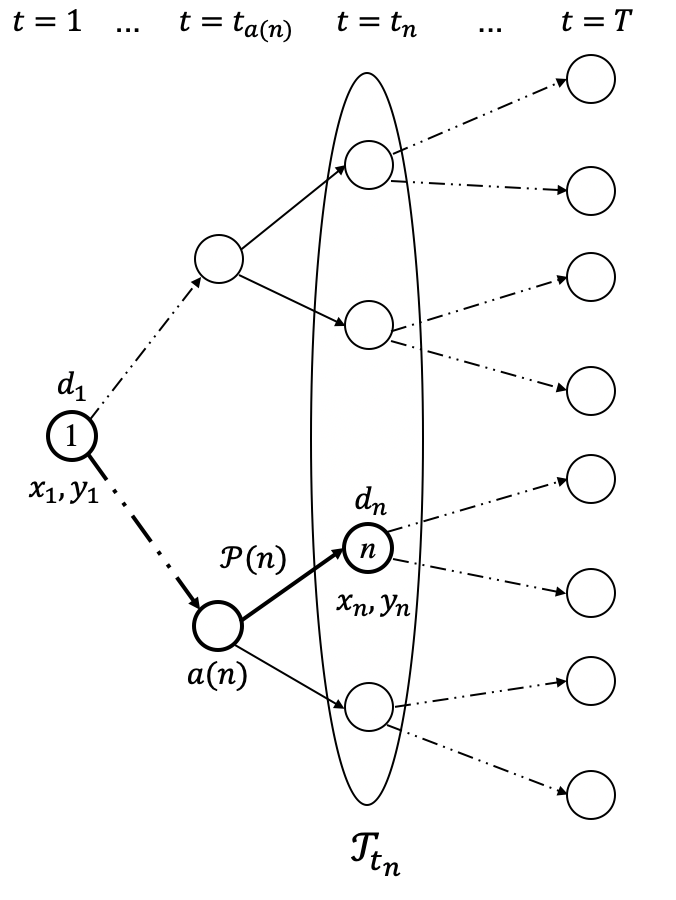}
        \caption{Scenario-node-based notation}
	\end{subfigure}
    \caption{Illustration of a scenario tree and its related notation.}
    \label{fig:scenario-tree}
\end{figure}

To facilitate formulating the stochastic programs, we shall introduce the scenario-path-based notation. We gather one possible path of realizations from the beginning to the end, and denote it as a scenario path $\omega$, i.e., $\boldsymbol{d}(\omega) = (\boldsymbol{d}_1(\omega),\ldots, \boldsymbol{d}_T(\omega))$.  Correspondingly, we use $\boldsymbol{x}_t(\omega),\ \boldsymbol{y}_t(\omega)$ to denote the decision vectors at stage $t$ under scenario $\omega$ (see Figure \ref{fig:scenario-tree}(a)). Consider a discrete distribution and assume that the number of realizations is finite, where such an approximation can be constructed by Monte Carlo sampling if the probability distribution is instead continuous and the resultant problem is called Sample Average Approximation (SAA) \citep{kleywegt2002sample}.  Let $\Omega$ be the support set of $\omega$ with each realization $\omega\in\Omega$ having probability $p(\omega)$ such that $\sum_{\omega\in\Omega}p(\omega)=1$. To ensure the subproblem maintained in each stage is always feasible for all decisions in the constraint set and for every realization of the random data, we refine the scenario set $\Omega$ so that the underlying assumption $\sum_{j=1}^N d_{tj}(\omega)\le\sum_{i=1}^M h_{i}$ holds for every $\omega\in\Omega$.

Instead of using decisions associated with scenario paths, we can equivalently define decision variables for each node in the scenario tree. Let $\mathcal{T}$ be {the set of all nodes in} the scenario tree associated with the underlying stochastic process. Each node $n$ in stage $t>1$ has a unique parent node $a(n)$ in stage $t-1$, and the set of children nodes of a node $n$ is denoted by $\mathcal{C}(n)$. The set $\mathcal{T}_t$ denotes the nodes corresponding to time period $t$, and $t_n$ is the time period corresponding to node $n$. Specially, all nodes in the last stage $\mathcal{T}_T$ are referred as leaf nodes and denoted by $\mathcal{L}$.  The path from the root node to node $n$ is denoted by $\mathcal{P}(n)$. 
Each node $n$ is associated with probability $p_n$. The probabilities of the nodes in each stage sum up to one, i.e., $\sum_{n\in\mathcal{T}_t}p_n = 1,\ \forall t=1,\ldots, T$, and the probabilities of all children nodes sum up to the probability of the parent node, i.e., $\sum_{m\in\mathcal{C}(n)}p_m=p_n,\ \forall n\not\in\mathcal{L}$.
If $n$ is a leaf node, i.e., $n\in\mathcal{L}$, then $\mathcal{P}(n)$ corresponds to a scenario $\omega\in\Omega$ and $p(\omega)=p_n$.
We denote the facility-location decision variable at node $n$ by $\boldsymbol{x}_n$, and the flow decision variable at node $n$ by $\boldsymbol{y}_n$. {Figure \ref{fig:scenario-tree} depicts and compares the scenario-path-based and scenario-node-based notation for a scenario tree representation of $T$-period uncertainty.}

In the following Sections \ref{sec:risktwo} and \ref{sec:riskmulti}, using the above notation, we describe how to formulate the risk-averse two-stage and multistage models, respectively.

\subsubsection{Risk-Averse Two-Stage Formulation}
\label{sec:risktwo}

We first introduce the definition and key properties of coherent risk measures. Consider a probability space $(\Xi,\mathcal{F},P)$. We refer a measurable function $Z: \Xi\to\mathbb{R}$ as a random variable. Let $\mathcal{Z}$ denote a space of $\mathcal{F}$-measurable functions from $\Xi$ to $\mathbb{R}$. With every random variable $Z\in \mathcal{Z}$, we associate a number, denoted as $\rho(Z)$, to indicate our preference between possible realizations of random variables. That is, $\rho(\cdot)$ is a real valued function $\rho:\mathcal{Z}\to\mathbb{R}$, which we call a risk measure.

According to \citet{artzner1999coherent}, a risk measure is a coherent risk measure if it satisfies the following properties:
\begin{enumerate}
\item Monotonicity: If $Z_1,Z_2\in\mathcal{Z}$ and $Z_1\succeq Z_2$, then $\rho(Z_1)\ge\rho(Z_2)$.
\item Convexity: $\rho(\gamma Z_1+(1-\gamma)Z_2)\le \gamma\rho(Z_1)+(1-\gamma)\rho(Z_2)$ for all $Z_1,\ Z_2\in\mathcal{Z}$ and all $\gamma\in[0, 1]$.
\item Translation invariance: If $a\in\mathbb{R}$ and $Z\in\mathcal{Z}$, then $\rho(Z+a)=\rho(Z)+a$.
\item Positive Homogeneity: If $\gamma\ge 0$ and $Z\in\mathcal{Z}$, then $\rho(\gamma Z)=\gamma\rho(Z)$.
\end{enumerate}
Here, $Z_1\succeq Z_2$ if and only if $Z_1(\xi)\ge Z_2(\xi)$ for a.e. $\xi\in\Xi$.

For our problem, we consider a special class of coherent risk measures -- a convex combination of expectation and Conditional Value-at-Risk (CVaR) \citep[see][]{rockafellar2000optimization}, i.e., for $t=2,\ldots,T$,
\begin{equation}
\rho_{t}(Z)=(1-\lambda_t)\mathbb{E}[Z]+\lambda_t \text{CVaR}_{\alpha_t}[Z],
\label{eq:rho}
\end{equation}
where $\lambda_t\in[0,1]$ is a parameter that compromises between optimizing on average and risk control, and $\alpha_t\in(0,1)$ is a number representing the confidence level. Notice that this risk measure is more general than CVaR and it includes CVaR as a special case when $\lambda_t=1$. 

Following the results by \cite{rockafellar2002conditional}, CVaR can be expressed as the following optimization problem:
\begin{equation}
\text{CVaR}_\alpha[Z]:=\inf_{\eta\in\mathbb{R}}\left\lbrace \eta+\frac{1}{1-\alpha}\mathbb{E}[Z-\eta]_{+}\right\rbrace ,\label{eq:cvar}
\end{equation}
where $[a]_{+}:=\max\{a,0\},$ and $\eta$ is an auxiliary variable. The minimum of the right-hand side of the above definition is attained at $\eta^* = \text{VaR}_{\alpha}[Z] := \inf\{v: \mathbb{P}(Z\le v)\ge \alpha\}$. To linearize $[Z-\eta]_{+}$, we replace it by a variable $u$ with two additional constraints: $u\ge 0,\ u\ge Z-\eta$.

In a two-stage stochastic program, we decide the facility locations $\boldsymbol{x}_1,\ldots,\boldsymbol{x}_T$ for all time periods in the first stage, and then evaluate the performance of first-stage decisions over the in-sample scenario set $\Omega$. Specifically, for each realized scenario path $\omega\in\Omega$ of the multiperiod demand $\boldsymbol{d}_1(\omega),\ldots,\boldsymbol{d}_T(\omega)$, we optimize the resource-allocation decision $\boldsymbol{y}_t(\omega)$ and calculate the operational cost $Q_t(\boldsymbol{x},\omega)$ for each period $t=1,\ldots, T$ independently. Note that because demand $\boldsymbol{d}_1$ is deterministic, the operational cost $Q_1(\boldsymbol{x},\omega)=Q_1(\boldsymbol{x})=\boldsymbol{c}^{\mathsf{T}}_{t_1}\boldsymbol{y}_1$ is also deterministic, while $Q_t(\boldsymbol{x},\omega),\ \forall t=2,\ldots,T$ is stochastic with respect to the scenario path $\omega$. The decision-making process is
\begin{eqnarray*}
\underbrace{\text{decision}\ (\boldsymbol{x}_1, \ldots, \boldsymbol{x}_T)}_{\textcolor{black}{\text{Stage 1}}}\to\underbrace{\text{observation}\ (\boldsymbol{d}_1,\ldots, \boldsymbol{d}_T)\to \text{decision}\ (\boldsymbol{y}_1,\ldots,\boldsymbol{y}_T)}_{\textcolor{black}{\text{Stage 2}}}.
\end{eqnarray*}
After fixing the facility locations $\boldsymbol{x}$, the flow decisions $\boldsymbol{y}_t(\omega)$ can be easily computed by solving a single-stage problem \eqref{eq:Q} described later that only depends on the demand realization $\boldsymbol{d}_t(\omega)$ in period $t$, and these decisions are stagewise disconnected, i.e., $\boldsymbol{y}_t(\omega)$ does not affect other stages' decisions. As a result, we consider a multiperiod risk function defined as the summation of the risk in each time period: $\mathbb{F}^{TS}(Q_1,Q_2,\ldots,Q_T)=Q_1+\sum_{t=2}^T\rho_t(Q_t)$, where each $\rho_t,\ t=2,\ldots,T$ is defined in \eqref{eq:rho}.
Then, a scenario-path-based formulation of two-stage risk-averse model with the multiperiod risk measure $\mathbb{F}^{TS}$ can be written as follows:
\begin{align}
z_R^{TS}=\min_{\boldsymbol{x}_1,\cdots,\boldsymbol{x}_T}\quad &\sum_{t=1}^T \boldsymbol{f}_{t}^{\mathsf T} \sum_{\tau=1}^t \boldsymbol{x}_{\tau}+{Q_1(\boldsymbol{x})}+\sum_{t=2}^T \rho_{t}\left(Q_t(\boldsymbol{x},\omega)\right)\label{model:twostagerisk-original}\\
\text{s.t.}\quad& \text{\eqref{eq:2-4}\ (constraints for } \boldsymbol{x})\nonumber\\
&\boldsymbol{x}_t\in \mathbb{Z}^{M}_{+}, \ \forall t=1,\ldots,T,\nonumber
\end{align}
where for each $t=1,\ldots,T,\ \omega\in\Omega$, 
\begin{subequations}\label{eq:Q}
\begin{align}
Q_t(\boldsymbol{x},\omega):=\min_{\boldsymbol{y}_t(\omega)} \quad&\boldsymbol{c}_t^{\mathsf T} \boldsymbol{y}_t(\omega)\nonumber\\
\text{s.t.}\quad& \boldsymbol{A}\boldsymbol{y}_t(\omega)=\boldsymbol{d}_t(\omega)\label{constraint_y1}\\
& \boldsymbol{B}_{t}\boldsymbol{y}_t(\omega)\le \sum_{\tau=1}^t \boldsymbol{x}_{\tau} \label{constraint_y3}\\
& \boldsymbol{y}_t(\omega)\in \mathbb{R}^{M\times N}_{+}. \nonumber
\end{align}
\end{subequations}
Problem \eqref{model:twostagerisk-original} represents the first-stage problem where the objective is to minimize the total cost of locating facilities and the risk measure of the {random operational cost, based on the outcomes of the functions   $Q_t(\boldsymbol{x},\omega)$ in \eqref{eq:Q} for each stage $t$ and scenario $\omega$, given facility-location decision $\boldsymbol{x}$.}

Using \eqref{eq:rho} and \eqref{eq:cvar}, a scenario-node-based formulation of the risk-averse two-stage model is given by:
\small
\begin{subequations}
\label{model:tworiskaverse-scenario-node}
\begin{align}
z^{TS}_R=\min_{\boldsymbol{x}_n,\boldsymbol{y}_n,\eta_n,u_n, n\in \mathcal{T}} \quad& {\sum_{n\in \mathcal{T}} p_n\left(\boldsymbol{f}_{t_n}^{\mathsf T}\sum_{m\in \mathcal{P}(n)}\boldsymbol{x}_m\right) + \boldsymbol{c}_{t_1}^{\mathsf T}\boldsymbol{y}_1+\sum_{n\in\mathcal{T},n\not=1}p_n\left((1-\lambda_{t_n})\boldsymbol{c}_{t_n}^{\mathsf T}\boldsymbol{y}_n + {\lambda_{t_n}} \eta_n + \frac{\lambda_{t_n}}{1-\alpha_{t_n}} u_n\right)}\nonumber\\
\text{s.t.} \quad& \boldsymbol{A}\boldsymbol{y}_n= \boldsymbol{d}_n, \ \forall n\in \mathcal{T}\label{eq:two-stage-constraint_y}\\
& \boldsymbol{B}_{t_n}\boldsymbol{y}_n\le \sum_{m\in \mathcal{P}(n)}\boldsymbol{x}_m,\ \forall n\in \mathcal{T}\\
& \sum_{m\in \mathcal{P}(n)}\boldsymbol{x}_m\le \mathbf{1},\ \forall n\in \mathcal{T}\label{eq:two-stage-constraint_x}\\
& u_n + \eta_{n} \ge \boldsymbol{c}_{t_n}^{\mathsf T} \boldsymbol{y}_n,\ \forall n\in \mathcal{T},\ {n\not=1} \label{eq:two-stage-constraint_eta}\\
& {\boldsymbol{x}_m}=\boldsymbol{x}_n,\ \forall m,\ n\in \mathcal{T}_t, \ t=1,...,T\label{eq:original_twostage_x}\\
& {\boldsymbol{\eta}_m}=\boldsymbol{\eta}_n,\ \forall m,\ n\in \mathcal{T}_t, \ t={2},...,T\label{eq:original_twostage_eta}\\
& \boldsymbol{x}_n\in \mathbb{Z}^{M}_{+},\ \boldsymbol{y}_n\in \mathbb{R}^{M\times N}_{+},\ u_n \ge 0,\ \forall n\in \mathcal{T},\nonumber
\end{align}
\end{subequations}
\normalsize
where constraints \eqref{eq:original_twostage_x}--\eqref{eq:original_twostage_eta} are the ``two-stage'' constraints that enforce the first-stage decisions $\boldsymbol{x}$ and $\eta$ to be identical for all nodes in the same stage. {We denote constraints \eqref{eq:two-stage-constraint_y}--\eqref{eq:two-stage-constraint_x} as $(\boldsymbol{x}_n,\boldsymbol{y}_n)\in X_n,\ \forall n\in\mathcal{T}$, and its scenario-path-based formulation as $(\boldsymbol{x}_t(\omega),\boldsymbol{y}_t(\omega))\in X_t(\omega),\ \forall t\in [T],\ \omega\in\Omega$.}

\subsubsection{Risk-Averse Multistage Formulation}
\label{sec:riskmulti}

In a multistage stochastic dynamic setting, the uncertain demand is revealed gradually, where we need to make both facility location $\boldsymbol{x}_t$ and flow decisions $\boldsymbol{y}_t$ in each stage $t$ based on the currently realized demand $\boldsymbol{d}_t$. Correspondingly, the decision-making process can be described as follows:
\begin{eqnarray*}
\underbrace{\text{decision}\ (\boldsymbol{x}_1, \boldsymbol{y}_1)}_{\text{Stage 1}}&\to& \underbrace{\text{observation}\ (\boldsymbol{d}_2)\to \text{decision}\ (\boldsymbol{x}_2, \boldsymbol{y}_2)}_{\text{Stage 2}}\to \text{observation}\ (\boldsymbol{d}_3)\\
&\to& \cdots\to \text{decision}\ (\boldsymbol{x}_{T-1}, \boldsymbol{y}_{T-1})\to \underbrace{\text{observation}\ (\boldsymbol{d}_{T})\to \text{decision}\ (\boldsymbol{x}_T, \boldsymbol{y}_T)}_{\text{Stage}\ T}.
\end{eqnarray*}

We consider the probability space $(\Xi,\mathcal{F},P)$, and let $\mathcal{F}_1\subset\mathcal{F}_2\subset\ldots\subset\mathcal{F}_T$ be sub-sigma-algebras of $\mathcal{F}$ such that each $\mathcal{F}_t$ corresponds to the information available up to (and including) stage $t$, with $\mathcal{F}_1=\{\emptyset,\Xi\}, \ \mathcal{F}_T=\mathcal{F}$. Let $\mathcal{Z}_t$ denote a space of $\mathcal{F}_t$-measurable functions from $\Xi$ to $\mathbb{R}$, and let $\mathcal{Z}:=\mathcal{Z}_1\times\cdots\times\mathcal{Z}_T$. We define a multiperiod risk function $\mathbb{F}^{MS}$ as a mapping from $\mathcal{Z}$ to $\mathbb{R}$ below:
\begin{equation}
\mathbb{F}^{MS}(Z_1,\ldots,Z_{T})=Z_1+\rho_2(Z_2)+\mathbb{E}_{\boldsymbol d_{[2]}}\left[{\rho_3^{\boldsymbol d_{[2]}}}(Z_3)\right]+\mathbb{E}_{\boldsymbol d_{[3]}}\left[{\rho_4^{\boldsymbol d_{[3]}}}(Z_4)\right]+\cdots+\mathbb{E}_{\boldsymbol d_{[T-1]}}\left[{\rho_{T}^{\boldsymbol d_{[T-1]}}}(Z_{T})\right],\label{ECRMs}
\end{equation}
where ${\rho_t^{\boldsymbol d_{[t-1]}}}$ is a {conditional} risk measure mapping from $\mathcal{Z}_t$ to $\mathcal{Z}_{t-1}$ to represent risk given the information available up to (including) stage $t-1$, {i.e., $\boldsymbol d_{[t-1]}=(\boldsymbol d_1,\ldots,\boldsymbol d_{t-1})$}. Because several risk measures are functions of expectation, their conditional counterparts correspond to replacing the expectation with conditional expectation. This class of multiperiod risk measures is called expected conditional risk measures (ECRMs), first introduced by \cite{homem2016risk}. 
We choose this class of multiperiod risk measures due to the following reasons:
\begin{itemize}
    \item \textbf{[Comparability]} The risk function \eqref{ECRMs} can be seen as a natural extension of the risk function $\mathbb{F}^{TS}$ that we defined in Section \ref{sec:risktwo}, where the demand is revealed gradually rather than all the demand is revealed at once. Because the facility location decisions are stagewise connected (i.e., $\boldsymbol{x}_t$ may affect the cost in future stages), we need to consider the dynamic process of the demand realization when measuring risk. Moreover, as we will show in Lemma \ref{lemma:equivalence}, the risk-averse two-stage and multistage models can be recast in an ``equivalent'' way such that the only differences are the ``two-stage'' constraints.
    \item \textbf{[Tractability]} The risk function \eqref{ECRMs} can be written in a nested form, and the corresponding risk-averse multistage stochastic programs can be recast as a risk-neutral counterpart with additional variables and constraints. Because of these two properties, we can explicitly write down the extensive form of risk-averse multistage models, which lays a foundation of comparing multistage models with two-stage counterparts directly. Note that not all multiperiod risk measures can be written as an extensive form, which we will introduce in Theorem \ref{thm:extensive}.
    \item \textbf{[Time Consistency]} An important property of dynamic risk measures is \textit{time consistency}, which ensures consistent risk preferences over stages. As defined in \cite{ruszczynski2010risk}, if a certain outcome is considered less risky in all states at stage $k$, then it should also be considered less risky at stage $l<k$. According to Section 5.1 of \citet{homem2016risk}, the risk function \eqref{ECRMs} can be recast as a composition of one-step conditional risk mappings in a nested way, and thus we can prove the time consistency of it. The detailed definitions and proof are presented in Appendix \ref{e-companion:ECRMs} and Theorem \ref{thm:ECRMs}.
\end{itemize}

For notation simplicity, denote $g_t(\boldsymbol{x}_{1:t},\boldsymbol{y}_t)= \boldsymbol{f}_{t}^{\mathsf T} \sum_{\tau=1}^t \boldsymbol{x}_{\tau}+\boldsymbol{c}_t^{\mathsf T} \boldsymbol{y}_t$ by $g_t$.
Plugging $g_t,\ t=1,\ldots, T$ into the formulation \eqref{ECRMs} and using tower property of expectations \citep{wasserman2004all}, the risk function \eqref{ECRMs} can be recast as
\begin{align}
\mathbb{F}(g_1,\ldots,g_{T})=g_1+\rho_2(g_2)+\mathbb{E}_{\boldsymbol d_{2}}\Big[{\rho_3^{\boldsymbol d_{[2]}}}(g_3)+\mathbb{E}_{\boldsymbol d_{3}|\boldsymbol d_{[2]}}\Big[{\rho_4^{\boldsymbol d_{[3]}}}(g_4)+\cdots
+\mathbb{E}_{\boldsymbol d_{T-1}|\boldsymbol d_{[T-2]}}\Big[{\rho_{T}^{\boldsymbol d_{[T-1]}}}(g _{T})\Big]\cdots\Big]\Big],\label{eq:obj}
\end{align}
where $\mathbb{E}_{\boldsymbol{d}_{t}|\boldsymbol d_{[t-1]}}$ represents the expectation with respect to the conditional probability distribution of $\boldsymbol{d}_{t}$ given realization $\boldsymbol d_{[t-1]}$.

Given parameters $\lambda_t \in [0, 1]$ and $\alpha_t\in(0,1)$ for all $t=2,\ldots,T$, we consider the conditional counterpart of the risk measure that we used in the two-stage model in Section \ref{sec:risktwo}, i.e.,  for $t=2,\ldots,T$,
\begin{equation}
{\rho_{t}^{\boldsymbol d_{[t-1]}}}(g_t)=(1-\lambda_t)\mathbb{E}[g_t{|\boldsymbol d_{[t-1]}}]+\lambda_t {\text{CVaR}_{\alpha_t}^{\boldsymbol d_{[t-1]}}}[g_t],
\label{eq:rho2}
\end{equation}
where $\text{CVaR}_{\alpha_t}^{\boldsymbol d_{[t-1]}}$ is the CVaR measure given the information $\boldsymbol d_{[t-1]}$, defined as:
\begin{equation}
\text{CVaR}_{\alpha_t}^{\boldsymbol d_{[t-1]}}[g_t]:=\inf_{\eta_t\in\mathbb{R}}\left\lbrace \eta_t+\frac{1}{1-\alpha_t}\mathbb{E}[[g_t-\eta_t]_{+}|\boldsymbol d_{[t-1]}]\right\rbrace.\label{eq:ccvar}
\end{equation}

\begin{figure}[ht!]
\centering
    \begin{subfigure}{0.48\textwidth}
        \centering
        \includegraphics[height = 6.5cm]{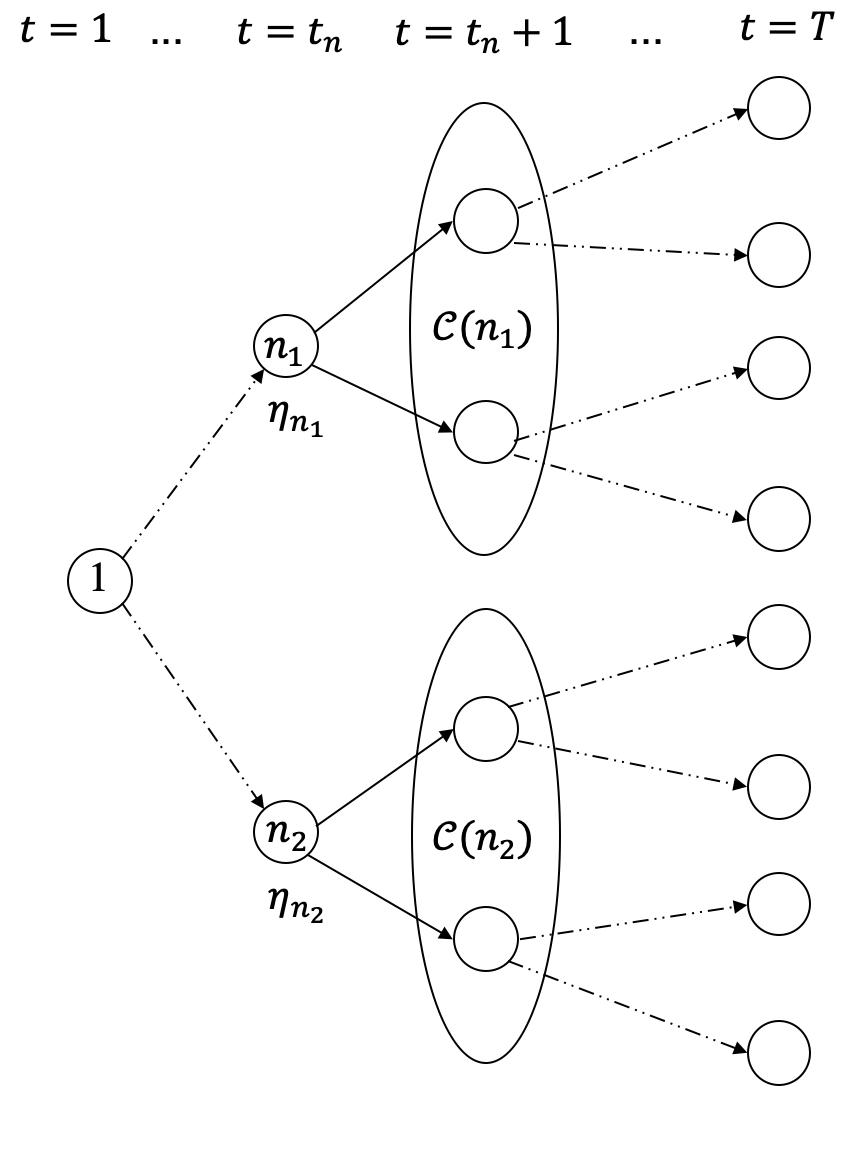}
        \caption{In a multistage setting: $\eta_n\ge \max_{m\in\mathcal{C}(n)}\{\boldsymbol{f}^{\mathsf T} \sum_{l\in \mathcal{P}(m)}\boldsymbol{x}_l+\boldsymbol{c}_m^{\mathsf T} \boldsymbol{y}_m-u_m\}$}
    \end{subfigure}
	\begin{subfigure}{0.48\textwidth}
	\centering
			\includegraphics[height = 6.5cm]{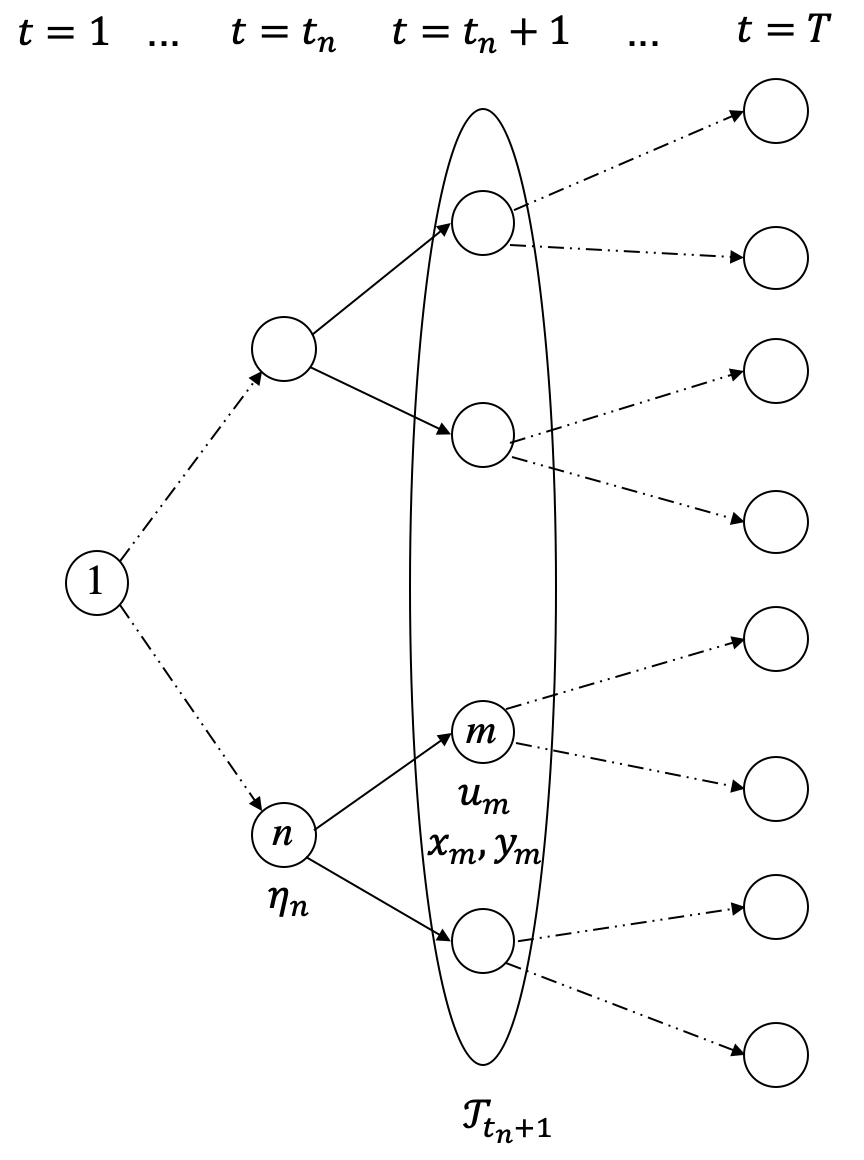}
    	\caption{In a two-stage setting: $\eta_n\ge \max_{m\in\mathcal{T}_{t_n+1}}\{\boldsymbol{f}^{\mathsf T} \sum_{l\in \mathcal{P}(m)}\boldsymbol{x}_l+\boldsymbol{c}_m^{\mathsf T} \boldsymbol{y}_m-u_m\}$}
		\end{subfigure}
    \caption{Illustration of constraints \eqref{eq:constraint_eta} in multistage and two-stage models.}
    \label{fig:scenario-node}
\end{figure}

Combining \eqref{eq:obj}, \eqref{eq:rho2} and \eqref{eq:ccvar}, the objective of the risk-averse multistage model is specified as
\begin{align}
z_R^{MS}=\min_{\substack{\boldsymbol{x}_1,\ldots,\boldsymbol{x}_T,\\\boldsymbol{y}_1,\ldots,\boldsymbol{y}_T,\\u_2,\ldots,u_T}}\quad&g_1+\min_{\eta_2}\lambda_2\eta_2+\mathbb{E}_{\boldsymbol {d}_2}\Big[\frac{\lambda_2}{1-\alpha_2}u_2+(1-\lambda_2)g_2\Big]\nonumber\\
&+\mathbb{E}_{\boldsymbol d_2}\Big[\min_{\eta_3}\lambda_3\eta_3+\mathbb{E}_{\boldsymbol d_3{|\boldsymbol d_{[2]}}}\Big[\frac{\lambda_3}{1-\alpha_3}u_3+(1-\lambda_3)g_3\Big]\nonumber\\
&+\mathbb{E}_{\boldsymbol d_3{|\boldsymbol d_{[2]}}}\Big[\min_{\eta_4}\lambda_4\eta_4+\mathbb{E}_{\boldsymbol d_4{|\boldsymbol d_{[3]}}}\Big[\frac{\lambda_4}{1-\alpha_4}u_4+(1-\lambda_4)g_4\Big]+\cdots\nonumber\\
&+\mathbb{E}_{\boldsymbol d_{T-1}{|\boldsymbol d_{[T-2]}}}\Big[\min_{\eta_T}\lambda_T\eta_T+\mathbb{E}_{\boldsymbol d_T{|\boldsymbol d_{[T-1]}}}\Big[\frac{\lambda_T}{1-\alpha_T}u_T+(1-\lambda_T)g_T\Big]\Big]\cdots\Big]\Big]\label{model:multirisk},
\end{align}
where the auxiliary variable $\eta_t\in\mathbb{R}$ is a function of $\boldsymbol{d}_1,\ldots,\boldsymbol{d}_{t-1}$, i.e., $\eta_t$ is a ``$(t-1)$-stage'' variable similar to $\boldsymbol{x}_{t-1}$ for all $t=2,\ldots,T$, and the auxiliary variable $u_t\in\mathbb{R}$ is a $t$-stage variable to represent the excess of $t$-stage cost of above $\eta_t$, with two additional constraints: $u_t\ge 0,\ u_t+\eta_t\ge g_t, \ \forall t=2,\ldots,T .$ 
\begin{theorem}\label{thm:extensive}
Using scenario-node-based notation {and recalling that $a(n)$ denotes the parent node of node $n$}, the risk-averse multistage model \eqref{model:multirisk} can be written in the following extensive form:
\begin{subequations}\label{model:multiriskaverse}
\begin{align}
z^{MS}_R=\min_{\substack{\boldsymbol{x}_n,\boldsymbol{y}_n, n\in \mathcal{T}\\\eta_n, n\not\in\mathcal{L}, u_n,n\not=1}} \quad& \sum_{n\in \mathcal{T}} p_n\left(\boldsymbol{\tilde{f}}_n^{\mathsf T}\sum_{m\in \mathcal{P}(n)}\boldsymbol{x}_m+\boldsymbol{\tilde{c}}_n^{\mathsf T}\boldsymbol{y}_n + \tilde{\lambda}_n \eta_n + \tilde{\alpha}_n u_n\right)\label{eq:multiriskaverse_obj}\\
\text{s.t.} \quad& {(\boldsymbol{x}_n,\boldsymbol{y}_n)\in X_n,\ \forall n\in\mathcal{T}}\label{eq:multi-constraint_xy}\\
& u_n + \eta_{a(n)} \ge \boldsymbol{f}_{t_n}^{\mathsf T} \sum_{m\in \mathcal{P}(n)}\boldsymbol{x}_m+\boldsymbol{c}_{t_n}^{\mathsf T} \boldsymbol{y}_n,\ \forall n\not= 1 \label{eq:constraint_eta}\\
& \boldsymbol{x}_n\in \mathbb{Z}^{M}_{+},\ \boldsymbol{y}_n\in \mathbb{R}^{M\times N}_{+},\ \forall n\in \mathcal{T},\ u_n \ge 0,\ \forall n\not=1,\nonumber
\end{align}
\end{subequations}
where $\boldsymbol{\tilde{f}}_n = \boldsymbol{f}_{t_n}$ if $n = 1$ and $\boldsymbol{\tilde{f}}_n = (1-\lambda_{t_n})\boldsymbol{f}_{t_n}$ otherwise; $\boldsymbol{\tilde{c}}_n = \boldsymbol{c}_{t_n}$ if $n = 1$ and $\boldsymbol{\tilde{c}}_n = (1-\lambda_{t_n})\boldsymbol{c}_{t_n}$ otherwise; $\tilde{\lambda}_n = 0$ if $n\in \mathcal{L}$ and $\tilde{\lambda}_n = \lambda_{t_n+1}$ otherwise; $\tilde{\alpha}_n = 0$ if $n = 1$ and $\tilde{\alpha}_n = \frac{\lambda_{t_n}}{1-\alpha_{t_n}}$ otherwise.
\end{theorem}
{We refer the interested readers to Appendix \ref{e-companion:proofs} for a detailed proof of Theorem \ref{thm:extensive}. }
\begin{remark}
Note that $\eta_n$-variables are defined for all nodes $n$ except for the leaf nodes $\mathcal{L}$ and $u_n$-variables are defined for non-root nodes $n\not=1$. For notational simplicity, we set the objective coefficients $\tilde{\lambda}_n=0,\ \forall n\in\mathcal{L}$ and $\tilde{\alpha}_n=0,\ n=1$, respectively.
We also depict the related variables of constraints \eqref{eq:constraint_eta} in Figure \ref{fig:scenario-node} to illustrate the difference between multistage and two-stage settings.
\end{remark}

\begin{remark}
When $\lambda_t=0,\ \forall t=2,\ldots, T$, the risk measure $\rho_t$ becomes the expectation and \eqref{model:multiriskaverse} is equivalent to a risk-neutral multistage stochastic dynamic program.
\end{remark}

Note that in the risk-averse multistage model \eqref{model:multiriskaverse}, both facility location decisions $\boldsymbol{x}$ and flow decisions $\boldsymbol{y}$ are recourse variables and thus they are both included in the risk measure (see constraints \eqref{eq:constraint_eta}). However, in the risk-averse two-stage model \eqref{model:tworiskaverse-scenario-node}, only the flow decisions are recourse variables and thus the facility location decisions $\boldsymbol{x}$ are not included in the risk measure (see constraints \eqref{eq:two-stage-constraint_eta}). This results in the differences in the objective functions and constraints between these two models. Nevertheless, after applying a variable transformation, the essential difference between these two models becomes whether the facility location decisions $\boldsymbol{x}$ and risk decisions $\eta$ are static or not, and the two models share the same objective function and constraints except for the ``two-stage'' constraints. Indeed, we can show that \eqref{model:tworiskaverse-scenario-node} can be reformulated as
\begin{subequations}
\label{model:tworiskaverse}
\begin{align}
z^{TS}_R=\min_{\substack{\boldsymbol{x}_n,\boldsymbol{y}_n, n\in \mathcal{T}\\\eta_n, n\not\in\mathcal{L}, u_n,n\not=1}}\quad & \sum_{n\in \mathcal{T}} p_n\left(\boldsymbol{\tilde{f}}_n^{\mathsf T}\sum_{m\in \mathcal{P}(n)}\boldsymbol{x}_m+\boldsymbol{\tilde{c}}_n^{\mathsf T}\boldsymbol{y}_n + \tilde{\lambda}_n \eta_n + \tilde{\alpha}_n u_n\right)\\
\text{s.t.}\quad&\text{\eqref{eq:multi-constraint_xy}--\eqref{eq:constraint_eta}}\nonumber\\
& {\boldsymbol{x}_m}=\boldsymbol{x}_n,\ \forall m,\ n\in \mathcal{T}_t, \ t=1,...,T\label{eq:twostage_x}\\
& {{\eta}_m}={\eta}_n,\ \forall m,\ n\in \mathcal{T}_t, \ t=1,...,T-1\label{eq:twostage_eta}\\
& \boldsymbol{x}_n\in \mathbb{Z}^{M}_{+},\ \boldsymbol{y}_n\in \mathbb{R}^{M\times N}_{+},\ \forall n\in \mathcal{T},\ u_n \ge 0,\ \forall n\not=1,\nonumber
\end{align}
\end{subequations}
in the following lemma, where the detailed proof is presented in Appendix \ref{e-companion:proofs}. 

\begin{lemma}\label{lemma:equivalence}
The risk-averse two-stage formulation \eqref{model:tworiskaverse-scenario-node} is equivalent to \eqref{model:tworiskaverse}. 
\end{lemma}

From Lemma \ref{lemma:equivalence}, we observe that the risk-averse two-stage model \eqref{model:tworiskaverse} is  the multistage model \eqref{model:multiriskaverse} with two additional constraints \eqref{eq:twostage_x} and \eqref{eq:twostage_eta}. Thus, ${\rm VMS_{R}}=z^{TS}_R-z^{MS}_{R}\ge 0$. We illustrate the difference of constraints \eqref{eq:constraint_eta} in the multistage and two-stage models due to the additional two-stage constraints in Figure \ref{fig:scenario-node}, which also leads to different analytical solutions of the substructure problem we will investigate in the following section.

\subsection{Analytical Solutions of the Substructure Problem}
\label{sec:substructure}
We first examine an important substructure of Models \eqref{model:multiriskaverse} and \eqref{model:tworiskaverse} once we fix $(\boldsymbol{y}, \boldsymbol{u})$-variables. We denote the resultant problems with known $(\boldsymbol{y}_n^*, {u}_n^*)$ values as $\mbox{{\bf SP-RMS}}(\boldsymbol{y}_n^{*}, u_n^{*})$ and $\mbox{{\bf SP-RTS}}(\boldsymbol{y}_n^{*}, u_n^{*})$, which are defined as follows:
\begin{subequations}\label{eq:ms-substructure}
\begin{align}
\mbox{{\bf SP-RMS}}(\boldsymbol{y}_n^{*}, u_n^{*}): \ \min_{\substack{\boldsymbol{x}_n, n\in \mathcal{T}\\\eta_n,,n\not\in\mathcal{L}}} \quad& \sum_{n\in \mathcal{T}} p_n\left(\boldsymbol{\tilde{f}}_n^{\mathsf T}\sum_{m\in \mathcal{P}(n)}\boldsymbol{x}_m+ \tilde{\lambda}_n \eta_n\right)\\
\text{s.t.}\quad& \sum_{m\in \mathcal{P}(n)}\boldsymbol{x}_m\ge \boldsymbol{B}_{t_n}\boldsymbol{y}_n^*,\ \forall n\in \mathcal{T}\label{eq:constraint_xy_knowny}\\
& \sum_{m\in \mathcal{P}(n)}\boldsymbol{x}_m\le \mathbf{1},\ \forall n\in \mathcal{T}\label{eq:constraint_x_knowny}\\
& \eta_{a(n)} \ge \boldsymbol{f}_{t_n}^{\mathsf T} \sum_{m\in \mathcal{P}(n)}\boldsymbol{x}_m+\boldsymbol{c}_{t_n}^{\mathsf T} \boldsymbol{y}_n^* - u^{*}_n,\ \forall n\not= 1 \label{eq:constraint_eta_knowny}\\
&\boldsymbol{x}_n\in \mathbb{Z}_{+}^{M},\ \forall n\in\mathcal{T},\nonumber
\end{align}
\end{subequations}
and
\begin{align}
\mbox{{\bf SP-RTS}}(\boldsymbol{y}_n^{*}, u_n^{*}): \ \min_{\substack{\boldsymbol{x}_n, n\in \mathcal{T}\\\eta_n,,n\not\in\mathcal{L}}}\quad & \sum_{n\in \mathcal{T}} p_n\left(\boldsymbol{\tilde{f}}_n^{\mathsf T}\sum_{m\in \mathcal{P}(n)}\boldsymbol{x}_m+ \tilde{\lambda}_n \eta_n\right) \label{eq:ts-substructure}\\
\text{s.t.}\quad& \text{\eqref{eq:constraint_xy_knowny}--\eqref{eq:constraint_eta_knowny}}\nonumber\\
&\text{\eqref{eq:twostage_x}, \eqref{eq:twostage_eta}\ (Two-stage constraints for $\boldsymbol{x}$ and $\boldsymbol{\eta}$)}\nonumber\\
&\boldsymbol{x}_n\in \mathbb{Z}_{+}^{M},\ \forall n\in\mathcal{T},\nonumber
\end{align}
{where $a(n)$ denotes the parent node of node $n$.}
Here, we denote the optimal objective values of Models \eqref{eq:ms-substructure} and \eqref{eq:ts-substructure} as $Q^M(\boldsymbol{y}_n^{*}, u_n^{*}),\ Q^{T}(\boldsymbol{y}_n^{*}, u_n^{*})$, respectively. The next proposition demonstrates the analytical forms of the optimal solutions to $\mbox{{\bf SP-RMS}}(\boldsymbol{y}_n^{*}, u_n^{*})$ and $\mbox{{\bf SP-RTS}}(\boldsymbol{y}_n^{*}, u_n^{*})$, {and we present a detailed proof in Appendix \ref{e-companion:proofs}.}

\begin{proposition}
\label{prop:substructure}
 Given $(\boldsymbol{y}^*, \boldsymbol{u}^*)$ values such that $\boldsymbol{B}_{t_n}\boldsymbol{y}_n^* \le 1,\ \forall n\in \mathcal{T}$, the optimal solutions of \eqref{eq:ms-substructure} and \eqref{eq:ts-substructure} have the following analytical forms:
	\begin{align*}
	&\boldsymbol{x}^{MS}_{1}=\lceil \boldsymbol{B}_{t_1} \boldsymbol{y}_1^{*}\rceil,\ \boldsymbol{x}^{MS}_{n}=\max_{m\in\mathcal{P}(n)}\lceil \boldsymbol{B}_{t_m} \boldsymbol{y}_m^{*}\rceil-\max_{m\in\mathcal{P}(a(n))}\lceil \boldsymbol{B}_{t_m} \boldsymbol{y}_m^{*}\rceil,\ \forall n\not=1,  \\
	&{\eta}^{MS}_n= \max_{m\in \mathcal{C}(n)}\left\lbrace \boldsymbol{f}_{t_m}^{\mathsf T}\sum_{l\in\mathcal{P}(m)}\boldsymbol{x}^{MS}_l+ \boldsymbol{c}_{t_m}^{\mathsf T} \boldsymbol{y}_m^{*}-u_m^{*}\right\rbrace,\ \forall n\not\in\mathcal{L},\\
	&\boldsymbol{x}^{TS}_1 = \lceil \boldsymbol{B}_{t_1} \boldsymbol{y}_1^{*}\rceil,\ \boldsymbol{x}^{TS}_{n}=\max_{m\in\mathcal{P}(n)}\lceil \max_{l\in \mathcal{T}_{t_m}}\boldsymbol{B}_{t_l} \boldsymbol{y}_l^{*}\rceil-\max_{m\in\mathcal{P}(a(n))}\lceil \max_{l\in \mathcal{T}_{t_m}}\boldsymbol{B}_{t_l} \boldsymbol{y}_l^{*}\rceil,\ \forall n\not=1,\\ 
	& \eta^{TS}_n = \max_{m\in \mathcal{T}_{t_n+1}}\left\lbrace \boldsymbol{f}_{t_m}^{\mathsf T}\sum_{l\in\mathcal{P}(m)}\boldsymbol{x}^{TS}_l+ \boldsymbol{c}_{t_m}^{\mathsf T} \boldsymbol{y}^*_m-u_m^{*}\right\rbrace, \ \forall n\not\in\mathcal{L},
	\end{align*}
	and correspondingly we have $Q^M(\boldsymbol{y}_n^{*}, u_n^{*}) = \sum_{n\in \mathcal{T}} p_n\left(\boldsymbol{\tilde{f}}_n^{\mathsf T}\sum_{m\in \mathcal{P}(n)}\boldsymbol{x}^{MS}_m+ \tilde{\lambda}_n \eta^{MS}_n\right), \	Q^{T}(\boldsymbol{y}_n^{*},u_n^{*}) = \sum_{n\in \mathcal{T}} p_n\left(\boldsymbol{\tilde{f}}_n^{\mathsf T}\sum_{m\in \mathcal{P}(n)}\boldsymbol{x}^{TS}_m + \tilde{\lambda}_n \eta^{TS}_n\right)$.
\end{proposition}

From Proposition \ref{prop:substructure}, one can easily verify that $Q^T(\boldsymbol{y}_n^{*}, u_n^{*}) - Q^M(\boldsymbol{y}_n^{*}, u_n^{*})\ge 0$ as $\sum_{m\in \mathcal{P}(n)}\boldsymbol{x}^{TS}_m = \max_{m\in\mathcal{P}(n)}\lceil \max_{l\in \mathcal{T}_{t_m}}\boldsymbol{B}_{t_l} \boldsymbol{y}_l^{*}\rceil \ge \max_{m\in\mathcal{P}(n)}\lceil \boldsymbol{B}_{t_m} \boldsymbol{y}_m^{*}\rceil = \sum_{m\in \mathcal{P}(n)}\boldsymbol{x}^{MS}_m$. Next, we will use the gap between the substructure problems to construct a lower bound on the ${\rm VMS_R}$.

\subsection{${\rm VMS_R}$ for Risk-Averse Facility Location}
We now describe a lower bound on the ${\rm VMS_R}$ for the risk-averse multistage and two-stage facility location models \eqref{model:multiriskaverse} and \eqref{model:tworiskaverse} based on the analysis in the previous section.
\label{sec:VMS}
\begin{theorem}
\label{thm:risk-averse-bound}
	Let $\{\boldsymbol{y}_n^{*}\}_{n\in \mathcal{T}}, \ \{u_n^{*}\}_{n\in \mathcal{T}\setminus\{1\}}$ be the second-stage decisions in an optimal solution to the two-stage model \eqref{model:tworiskaverse}, and let $\boldsymbol{x}^{MS},\boldsymbol{\eta}^{MS},\boldsymbol{x}^{TS},\boldsymbol{\eta}^{TS}$ follow the definitions in Proposition \ref{prop:substructure}, which are constructed by $\{\boldsymbol{y}_n^{*}\}_{n\in \mathcal{T}}, \ \{u_n^{*}\}_{n\in \mathcal{T}\setminus\{1\}}$.
	Then,
	\begin{align}
	{\rm VMS_R}&\ge  \sum_{n\in \mathcal{T}\setminus\{1\}} p_n{(1-\lambda_{t_n})\boldsymbol{f}_{t_n}^{\mathsf T}}\left(\max_{m\in \mathcal{P}(n)} \lceil \max_{l\in \mathcal{T}_{t_m}}\boldsymbol{B}_{t_l} \boldsymbol{y}_l^{*}\rceil - \max_{m\in \mathcal{P}(n)}\lceil \boldsymbol{B}_{t_m} \boldsymbol{y}_m^{*}\rceil\right)+ {\sum_{n\in \mathcal{T}\setminus\mathcal{L}} p_n {\lambda}_{t_n+1}} \left(\eta^{TS}_n-\eta^{MS}_n\right),\nonumber
	\end{align}
	where we denote the right-hand side as ${\rm VMS_R^{LB}}$. 
\end{theorem}
We refer the interested readers to Appendix \ref{e-companion:proofs} for a detailed proof of Theorem \ref{thm:risk-averse-bound}. 

\begin{remark}
    The lower bound ${\rm VMS_R^{LB}}$ provided in Theorem \ref{thm:risk-averse-bound} is \textit{tight}, {and we demonstrate it using Example \ref{eg1} in Appendix \ref{e-companion:example} where the equality holds (i.e., ${\rm VMS_R} = {\rm VMS_R^{LB}}$).}
\end{remark}

Next we derive two more lower bounds ${\rm VMS_R^{LB1}}, {\rm VMS_R^{LB2}}$ that are not necessarily tight but more computationally tractable, where ${\rm VMS_R^{LB1}}$ utilizes the optimal solutions to the linear programming (LP) relaxation of the two-stage model \eqref{model:tworiskaverse}, and ${\rm VMS_R^{LB2}}$ only uses input parameters. {The detailed proofs are presented in Appendix \ref{e-companion:proofs}.}
\begin{corollary}\label{cor:LP}
Let $\{\boldsymbol{y}_n^{LP}\}_{n\in \mathcal{T}}, \ \{u_n^{LP}\}_{n\in \mathcal{T}\setminus\{1\}}$ be the second-stage decisions in an optimal solution to the LP relaxation of the two-stage model \eqref{model:tworiskaverse}, and let
	\begin{align*}
	&\boldsymbol{x}^{MS}_{1}=\lceil\boldsymbol{B}_{t_1} \boldsymbol{y}_1^{LP}\rceil,\ \boldsymbol{x}^{MS}_{n}=\max_{m\in\mathcal{P}(n)} \lceil\boldsymbol{B}_{t_m} \boldsymbol{y}_m^{LP}\rceil-\max_{m\in\mathcal{P}(a(n))} \lceil\boldsymbol{B}_{t_m} \boldsymbol{y}_m^{LP}\rceil,\ \forall n\not=1,  \\
	&{\eta}^{MS}_n= \max_{m\in \mathcal{C}(n)}\left\lbrace \boldsymbol{f}_{t_m}^{\mathsf T}\sum_{l\in\mathcal{P}(m)}\boldsymbol{x}^{MS}_l+ \boldsymbol{c}_{t_m}^{\mathsf T} \boldsymbol{y}_m^{LP}-u_m^{LP}\right\rbrace,\ \forall n\not\in\mathcal{L},\\
	&\boldsymbol{x}^{TS}_1 =  \boldsymbol{B}_{t_1} \boldsymbol{y}_1^{LP},\ \boldsymbol{x}^{TS}_{n}=\max_{m\in\mathcal{P}(n)} \max_{l\in \mathcal{T}_{t_m}}\boldsymbol{B}_{t_l} \boldsymbol{y}_l^{LP}-\max_{m\in\mathcal{P}(a(n))} \max_{l\in \mathcal{T}_{t_m}}\boldsymbol{B}_{t_l} \boldsymbol{y}_l^{LP},\ \forall n\not=1,\\ 
	& \eta^{TS}_n = \max_{m\in \mathcal{T}_{t_n+1}}\left\lbrace \boldsymbol{f}_{t_m}^{\mathsf T}\sum_{l\in\mathcal{P}(m)}\boldsymbol{x}^{TS}_l+ \boldsymbol{c}_{t_m}^{\mathsf T} \boldsymbol{y}^{LP}_m-u_m^{LP}\right\rbrace, \ \forall n\not\in\mathcal{L}.
	\end{align*}
Then,
\small
	\begin{align*}
	{\rm VMS_R}\ge &\boldsymbol{f}_{t_1}^{\mathsf T}\left(\boldsymbol{B}_{t_1} \boldsymbol{y}_1^{LP}-\lceil\boldsymbol{B}_{t_1} \boldsymbol{y}_1^{LP}\rceil\right)+\sum_{n\in \mathcal{T}\setminus\{1\}} p_n(1-\lambda_{t_n})\boldsymbol{{f}}_{t_n}^{\mathsf T}\left(\max_{m\in \mathcal{P}(n)}  \max_{l\in \mathcal{T}_{t_m}}\boldsymbol{B}_{t_l} \boldsymbol{y}_l^{LP} - \max_{m\in \mathcal{P}(n)} \lceil\boldsymbol{B}_{t_m} \boldsymbol{y}_m^{LP}\rceil\right) \\
	&+\sum_{n\in \mathcal{T}\setminus\mathcal{L}} p_n {\lambda}_{t_n+1} \left(\eta^{TS}_n-\eta^{MS}_n\right),
	\end{align*}
	\normalsize
	where we denote the right-hand side as ${\rm VMS_R^{LB1}}$.
\end{corollary}
\begin{corollary}\label{cor:parameter}
We can also derive a lower bound that only depends on the problem parameters as follows:
\begin{align*}
	{\rm VMS_R^{LB}}\ge  \sum_{n\in \mathcal{T}\setminus\{1\}} p_n(1-\lambda_{t_n} )\sum_{i=1}^M f_{t_n,i}\Delta_{ni}
	\end{align*}
	where we denote the right-hand side as ${\rm VMS_R^{LB2}}$.
 For each $i\in[M]$ and $n\in\mathcal{T}\setminus\{1\}$, $\Delta_{ni}=1$ if the following conditions are satisfied:
	\begin{itemize}
	    \item Condition (i): $d_{m,j}=0,\ \forall j\in [N],\ m\in\mathcal{P}(n)$
	    \item Condition (ii): there exists $\bar{n}\in\cup_{m\in\mathcal{P}(n)}\mathcal{T}_{t_m}$ such that $\sum_{j\in [N]}d_{\bar{n},j}>\sum_{i^{\prime}\in[M],i^{\prime}\not=i}h_{t_{\bar{n}}i^{\prime}}$.
	\end{itemize}
	and $\Delta_{ni}=0$ otherwise.
\end{corollary}

\subsection{Approximation Algorithms}
\label{sec:approx}
In this section, we propose approximation algorithms to solve the risk-averse multistage and two-stage programs \eqref{model:multiriskaverse} and \eqref{model:tworiskaverse}. First, both models are inherently hard to solve indicated in the following theorem, {and we present a detailed proof in Appendix \ref{e-companion:proofs}.}

\begin{theorem}
\label{thm:hardness}
    The deterministic facility location problem \eqref{model:deter} and its risk-averse multistage and two-stage counterparts \eqref{model:multiriskaverse} and \eqref{model:tworiskaverse} are NP-hard. 
\end{theorem}

Motivated by the computational intractability of Models \eqref{model:multiriskaverse} and \eqref{model:tworiskaverse}, we proceed to introduce approximation algorithms that can solve the risk-averse multistage and two-stage models efficiently by utilizing the decomposition structure we investigate in Section \ref{sec:substructure}. Next, we describe the main idea of the algorithm for the risk-averse multistage model \eqref{model:multiriskaverse} as follows: we first solve the LP relaxation of Model \eqref{model:multiriskaverse} to obtain a feasible solution $(\boldsymbol{y}_n^{LP},{u}_n^{LP})$, which is fed into the substructure problem $\mbox{{\bf SP-RMS}}(\boldsymbol{y}_n^{LP}, u_n^{LP})$ to obtain an optimal solution $(\boldsymbol{x}_n^1, \eta_n^1)$. We then solve Model \eqref{model:multiriskaverse} with fixed $(\boldsymbol{x}_n, \eta_n) = (\boldsymbol{x}_n^1, \eta_n^1)$ to derive an optimal solution $(\boldsymbol{y}_n^{1},{u}_n^{1})$, which together with $(\boldsymbol{x}_n^1, \eta_n^1)$ constitutes a feasible solution and thus an upper bound to Model \eqref{model:multiriskaverse}. This upper bound can be strengthened iteratively by repeating the process, {and we denote the feasible solution produced at the end of Algorithm \ref{alg:approx-multistage} by $(\boldsymbol{x}_n^{H}, \eta_n^{H}, \boldsymbol{y}_n^{H}, u_n^{H})$}. The detailed steps are described in Algorithm \ref{alg:approx-multistage}. We show the monotonicity of the upper bounds derived in Algorithm \ref{alg:approx-multistage} in Proposition \ref{prop:alg}, and show that the optimality gap of Algorithm \ref{alg:approx-multistage} can be upper bounded in Proposition \ref{prop:ratio}, which will eventually lead to an approximation ratio stated in Theorem \ref{thm:ratio}.
\begin{algorithm}[ht!]
\caption{Approximation Algorithm for Risk-Averse Multistage Facility Location \eqref{model:multiriskaverse}}
\begin{algorithmic}[1]
\label{alg:approx-multistage}
\STATE Solve the LP relaxation of the risk-averse multistage facility location problem \eqref{model:multiriskaverse} and let $(\boldsymbol{x}_n^{LP}, \eta_n^{LP}, \boldsymbol{y}_n^{LP}, u_n^{LP})_{n\in\mathcal{T}}$ be an optimal solution. If $\boldsymbol{x}_n^{LP}$ is integral for all $n\in\mathcal{T}$, stop and return $(\boldsymbol{x}_n^{LP}, \eta_n^{LP}, \boldsymbol{y}_n^{LP}, u_n^{LP})_{n\in\mathcal{T}}$.
\STATE Initialize $k=0$ and $(\boldsymbol{x}_n^0, \eta_n^0, \boldsymbol{y}_n^0, u_n^0)_{n\in\mathcal{T}} = (\boldsymbol{x}_n^{LP}, \eta_n^{LP}, \boldsymbol{y}_n^{LP}, u_n^{LP})_{n\in\mathcal{T}}$.
\WHILE{$||\boldsymbol{x}^k - \boldsymbol{x}^{k-1}|| \ge \epsilon,\ ||\eta^k - \eta^{k-1}|| \ge \epsilon,\ ||\boldsymbol{y}^k - \boldsymbol{y}^{k-1}|| \ge \epsilon,\ ||u^k - u^{k-1}|| \ge \epsilon$}
\STATE \label{alg:step4} Solve Problem $\mbox{{\bf SP-RMS}}(\boldsymbol{y}_n^{k}, u_n^{k})$
and let $\boldsymbol{x}_n^{k+1},\ {\eta}_n^{k+1}$ denote the corresponding optimal solutions. We have the analytical form of the optimal solutions as $\boldsymbol{x}^{k+1}_{1}=\lceil \boldsymbol{B}_{t_1} \boldsymbol{y}_1^{k}\rceil,\ \boldsymbol{x}_n^{k+1} = \max_{m\in\mathcal{P}(n)}\lceil \boldsymbol{B}_{t_m} \boldsymbol{y}_m^{k}\rceil-\max_{m\in\mathcal{P}(a(n))}\lceil \boldsymbol{B}_{t_m} \boldsymbol{y}_m^{k}\rceil,\ \forall n\not=1,\ {\eta}^{k+1}_n= \max_{m\in \mathcal{C}(n)}\left\lbrace \boldsymbol{f}_{t_m}^{\mathsf T}\sum_{l\in\mathcal{P}(m)}\boldsymbol{x}^{k+1}_l+ \boldsymbol{c}_{t_m}^{\mathsf T} \boldsymbol{y}_m^{k}-u_m^{k}\right\rbrace,\ \forall n\not\in\mathcal{L}$.

\STATE \label{alg:step5} Solve the following problem for each $n\in\mathcal{T},\ n\not=1$ independently
\begin{subequations}\label{eq:independent}
\begin{align}
		\min_{{\boldsymbol{y}_n}, u_n} \quad&   \tilde{\boldsymbol{c}}_{n}^{\mathsf T}\boldsymbol{y}_n + \tilde{\alpha}_nu_n \nonumber\\
\text{s.t.}	\quad	&    \boldsymbol{B}_{t_n}\boldsymbol{y}_n \le \sum_{m\in \mathcal{P}(n)}{\boldsymbol{x}^{k+1}_{m}},\\
		& \boldsymbol{A}\boldsymbol{y}_n =\boldsymbol{d}_n,\\
		& u_n - \boldsymbol{c}_{t_n}^{\mathsf T}\boldsymbol{y}_n \ge \boldsymbol{f}_{t_n}^{\mathsf T} \sum_{m\in \mathcal{P}(n)}\boldsymbol{x}^{k+1}_{m} - \eta^{k+1}_{a(n)},\label{indepent-u}\\
		& \boldsymbol{y}_n\in\mathbb{R}_{+},\ u_n\ge 0,\nonumber
\end{align}
\end{subequations}
and when $n=1$ we solve Problem \eqref{eq:independent} without the variables $u_n$ and constraints \eqref{indepent-u}.
Let $\boldsymbol{y}_n^{k+1},\ {u}_n^{k+1}$ be the optimal solutions.
\STATE Update $k = k+1$.
\ENDWHILE
\STATE Return $(\boldsymbol{x}_n^H, \eta_n^H, \boldsymbol{y}_n^H, u_n^H)_{n\in\mathcal{T}} := (\boldsymbol{x}_n^{k+1}, \eta_n^{k+1}, \boldsymbol{y}_n^{k+1}, u_n^{k+1})_{n\in\mathcal{T}}$.
\end{algorithmic}
\end{algorithm}

\begin{proposition}\label{prop:alg}
The objective value at the end of each iteration {in Algorithm \ref{alg:approx-multistage}} provides an upper bound to the optimal objective value of Model \eqref{model:multiriskaverse} and {it is improved or stays the same after each iteration}, i.e.,
    $z_R^{MS}(\boldsymbol{x}_n^{*}, \eta_n^{*}, \boldsymbol{y}_n^{*}, u_n^{*})\le z_R^{MS}(\boldsymbol{x}_n^{k+1}, \eta_n^{k+1}, \boldsymbol{y}_n^{k+1}, u_n^{k+1}) \le z_R^{MS}(\boldsymbol{x}_n^{k}, \eta_n^{k}, \boldsymbol{y}_n^{k}, u_n^{k}),\ \forall k\ge 1$. (With some abuse of notation, we use $(\boldsymbol{x}_n^{*}, \eta_n^{*}, \boldsymbol{y}_n^{*}, u_n^{*})$ to denote an optimal solution to the multistage model \eqref{model:multiriskaverse} and use $z_R^{MS}(\boldsymbol{x}_n, \eta_n, \boldsymbol{y}_n, u_n)$ to denote the corresponding objective value with input values $(\boldsymbol{x}_n, \eta_n, \boldsymbol{y}_n, u_n)$ to Model \eqref{model:multiriskaverse}.)
\end{proposition}

\begin{proposition}\label{prop:ratio}
The optimality gap can be bounded above by a quantity only dependent on the facility location cost, i.e., $z_R^{MS}(\boldsymbol{x}_n^{H}, \eta_n^{H}, \boldsymbol{y}_n^{H}, u_n^{H}) - z_R^{MS}(\boldsymbol{x}_n^{*}, \eta_n^{*}, \boldsymbol{y}_n^{*}, u_n^{*}) \le {\sum_{t=1}^T\sum_{i=1}^M f_{ti}}$.
\end{proposition}

\begin{theorem}\label{thm:ratio}
Algorithm \ref{alg:approx-multistage} has an approximation ratio of $$ {1+\frac{M\sum_{t=1}^Tf_{t,\rm max} }{M_{\rm min}\sum_{t=1}^Tf_{t,\rm min}+\sum_{t=1}^Tc_{t,\rm min}\min_{n\in\mathcal{T}_t}\{\sum_{j=1}^Nd_{n,j}\}}},$$ where $h_{\rm max} = \max_{i=1}^M\{h_{1i}\},\ f_{t, \rm max} = \max_{i=1}^M\{f_{ti}\},\ f_{t,\rm min} = \min_{i=1}^M\{f_{ti}\},\ {c_{t,\rm min}=\min_{i\in[M],j\in[N]}c_{tij}}$ and $M_{\rm min} = \lceil\frac{\sum_{j=1}^N d_{1j}}{h_{\rm max}}\rceil$ measures at least how many facilities we need to cover the first-stage demand. 
\end{theorem}

\begin{corollary}\label{cor:ratio}
Assume that $f_{t,\rm max}=O(1),\ f_{t,\rm min}=O(1),\ c_{t,\rm min}=O(1),\ \min_{n\in\mathcal{T}_t}\sum_{j=1}^Nd_{n,j}=O(t)$ when $t\to\infty$. Then Algorithm \ref{alg:approx-multistage} is asymptotically optimal, i.e.,
\begin{align*}
    \lim_{T\to\infty}\frac{z_R^{MS}(\boldsymbol{x}_n^{H}, \eta_n^{H}, \boldsymbol{y}_n^{H}, u_n^{H})}{z_R^{MS}(\boldsymbol{x}_n^{*}, \eta_n^{*}, \boldsymbol{y}_n^{*}, u_n^{*})}=1.
\end{align*}
\end{corollary}

Note that the assumptions in Corollary \ref{cor:ratio} are not particularly restrictive. They only require that the facility location and unit flow costs are constant with respect to the time stage (i.e., they can be bounded above by a value that does not depend on the time), and the demand grows at least linearly when time increases (the result still holds when the demand grows faster than linearly, e.g., when $\min_{n\in\mathcal{T}_t}\sum_{j=1}^Nd_{n,j}=O(t^2)$). The last condition can be achieved if we have an expanding market and the minimum of the demand is increased for each subsequent year (e.g., $\min_{n\in\mathcal{T}_t}\sum_{j=1}^Nd_{n,j}=\tilde{d}(1+0.2(t-1))$ with $\tilde{d}$ being the nominal demand in the first stage). We will test the linearly increasing demand pattern in Section \ref{sec:real}. Detailed proofs of Propositions \ref{prop:alg}, \ref{prop:ratio}, Theorem \ref{thm:ratio} {and Corollary \ref{cor:ratio}} are given in Appendix \ref{e-companion:proofs}.

One can also tailor Algorithm \ref{alg:approx-multistage} to solve the risk-averse two-stage model \eqref{model:tworiskaverse} by modifying Step \ref{alg:step4} with the analytical solutions of the two-stage model as stated in Proposition \ref{prop:substructure}. We present the detailed steps of Algorithm \ref{alg:approx-twostage} for approximating solutions to the risk-averse two-stage model in Appendix \ref{e-companion:alg}: Algorithm \ref{alg:approx-twostage}. All the other results still hold and follow similar proofs, which we present in Propositions \ref{prop:alg-twostage}, \ref{prop:ratio-twostage} and Theorem \ref{thm:ratio-twostage} without proof in the interest of brevity.

\begin{proposition}\label{prop:alg-twostage}
The objective value at the end of each iteration in Algorithm \ref{alg:approx-twostage} provides an upper bound to the optimal objective value of Model \eqref{model:tworiskaverse} and it never gets worse, i.e.,
    $z_R^{TS}(\boldsymbol{x}_n^{*}, \eta_n^{*}, \boldsymbol{y}_n^{*}, u_n^{*})\le z_R^{TS}(\boldsymbol{x}_n^{k+1}, \eta_n^{k+1}, \boldsymbol{y}_n^{k+1}, u_n^{k+1}) \le z_R^{TS}(\boldsymbol{x}_n^{k}, \eta_n^{k}, \boldsymbol{y}_n^{k}, u_n^{k}),\ \forall k\ge 1$, where we use $(\boldsymbol{x}_n^{*}, \eta_n^{*}, \boldsymbol{y}_n^{*}, u_n^{*})$ to denote an optimal solution to the two-stage model \eqref{model:tworiskaverse}.
\end{proposition}

\begin{proposition}\label{prop:ratio-twostage}
$z_R^{TS}(\boldsymbol{x}_n^{H}, \eta_n^{H}, \boldsymbol{y}_n^{H}, u_n^{H}) - z_R^{TS}(\boldsymbol{x}_n^{*}, \eta_n^{*}, \boldsymbol{y}_n^{*}, u_n^{*}) \le {\sum_{t=1}^T\sum_{i=1}^M f_{ti}}$.
\end{proposition}

\begin{theorem}\label{thm:ratio-twostage}
Algorithm \ref{alg:approx-twostage} has an approximation ratio of 
\begin{align*}
{1+ {\frac{M\sum_{t=1}^Tf_{t,\rm max} }{M_{\rm min}\sum_{t=1}^Tf_{t,\rm min}+\sum_{t=1}^Tc_{t,\rm min}\min_{n\in\mathcal{T}_t}\{\sum_{j=1}^Nd_{n,j}\}}}.}
\end{align*}
\end{theorem}

{
\section{Value of Risk-Averse Multistage Facility Location with Prioritization}\label{sec:VMS-prior}
In this section, we consider both demand and budget uncertainty, and extend the prior work by \citet{kocc2015prioritization} to prioritize the candidate facility sites in both two-stage and multistage settings. A priority list is a many-to-one assignment of candidate sites to priority levels such that each priority level contains at least one candidate site. In a two-stage stochastic program, the priority list is decided in advance and fixed for all periods, while in a multistage stochastic dynamic program, the priority list can change adaptively over time. After realizing the uncertainty, we need to enforce the site selections to obey the priority list, i.e., a lower-priority candidate site cannot be selected unless all higher-priority candidate sites are selected. We set up the problem formulations in Section \ref{sec:problem-prior}, examine a substructure problem in Section \ref{sec:substructure-prior}, and based on that, we derive lower bounds for the gaps between the objective values of the risk-averse two-stage and multistage models with prioritization in Section \ref{sec:proof-VMS-prior}. We also develop a set of cutting planes called \textit{prioritization cuts} to improve the computation in Section \ref{sec:prioritization-cuts}.

\subsection{Problem Formulations}\label{sec:problem-prior}
We first add a stage 0 in the beginning of the planning horizon to decide an initial priority list and consider stages $0,1,\ldots, T$. We define $s_{t,ii^{\prime}}\in\{0,1\}$ as the priority list at stage $t$ for all $t=0,1,\ldots,T-1$ between each pair of facilities $i,\ i^{\prime}=1,\ldots, M,\ i\not=i^{\prime}$, such that $s_{t,ii^{\prime}}=1$ if $i$ does not have lower priority than $i^{\prime}$ in stage $t+1$ and 0 otherwise. Note that the facility location decisions at stage $t+1$ (i.e., $\boldsymbol{x}_{t+1}$) must obey the priority list decided in the previous stage (i.e., $\boldsymbol{s}_{t}$). Denote $F_t$ as the total budget for opening facilities at stage $t$ for all $t=1,\ldots,T$. In this section, for the ease of presentation, we assume that all the candidate sites have the same location setup cost and the same capacity (i.e., $f_{ti}=f_{t},\ h_{ti}=h_{t},\ \forall i\in[M]$), and thus the opening budget $F_t$ can be simplified as the maximum number of facilities to be open in each stage. Note that the following models and results can be extended to include varying setup cost and capacity as well. We aim to minimize $\sum_{t=0}^{T-1}\sum_{i\not=i^{\prime}}s_{t,ii^{\prime}}$ while variables $\boldsymbol{s}$ and $\boldsymbol{x}$ are subject to the following constraints:
\begin{subequations}\label{model:deterministic-prior}
\begin{align}
  &  s_{0,ii^{\prime}}+s_{0,i^{\prime}i}\ge 1,\ \forall 1\le i<i^{\prime}\le M \label{eq:prior-stage1}\\
& s_{t,ii^{\prime}}+s_{t,i^{\prime}i} + \sum_{\tau=1}^t x_{\tau,i} + \sum_{\tau=1}^tx_{\tau,i^{\prime}}\ge 1,\ \forall 1\le i<i^{\prime}\le M,\ 1\le t\le T-1 \label{eq:prior-staget}\\
& \sum_{\tau=1}^tx_{\tau,i} \ge \sum_{\tau=1}^tx_{\tau,i^{\prime}} + s_{t-1, ii^{\prime}} - 1, \ \forall 1\le i\not=i^{\prime}\le M,\ 1\le t\le T \label{eq:prior-xs}\\
& \sum_{i=1}^M x_{t,i} \le F_t,\ \forall 1\le t\le T \label{eq:prior-budget}\\
& s_{t,ii^{\prime}}\in\{0,1\},\ \forall  1\le i\not=i^{\prime}\le M,\ 0\le t\le T-1.
\end{align}
\end{subequations}
Here constraints \eqref{eq:prior-stage1} ensure that given a pair of facility sites $i,\ i^{\prime}$, either they have the same priority or one has higher priority than the other in the initial setting. If $i$ has higher priority than $i^{\prime}$, we have $s_{0,ii^{\prime}}=1, \ s_{0,i^{\prime}i}=0$ and constraints \eqref{eq:prior-xs} read as $x_{1, i}\ge x_{1,i^{\prime}}$ and $x_{1,i^{\prime}}\ge x_{1,i}-1$ where the latter is redundant.   On the other hand, if $s_{0,ii^{\prime}}=s_{0,i^{\prime}i}=1$, then $i$ and $i^{\prime}$ have the same priority and constraints \eqref{eq:prior-xs} yield $x_{1, i}=x_{1,i^{\prime}}$. One can also add cycle-elimination constraints (e.g., $s_{t,i_1i_2}+s_{t,i_2i_1}\le 1,\ s_{t,i_1i_2}+s_{t,i_2i_3}+s_{t,i_3i_1}\le 2,\ \forall i_1\not=i_2\not=i_3,\ t\in[T]$) to Model \eqref{model:deterministic-prior} to avoid multiple candidate sites in one priority level (e.g., $s_{t,i_1i_2}=s_{t,i_2i_3}=s_{t,i_3i_1}=1$). We refer interested readers to Section 4 in \citet{kocc2015prioritization} for detailed comparison between many-to-one and one-to-one assignment of candidate sites to priority levels.
According to constraints \eqref{eq:prior-staget}, when either $\sum_{\tau=1}^t x_{\tau,i}=1$ or $\sum_{\tau=1}^tx_{\tau,i^{\prime}}=1$ (i.e., site $i$ or $i^{\prime}$ is already open in stage $t$), $s_{t,ii^{\prime}}=s_{t,i^{\prime}i}=0$ because we are minimizing $\sum_{i\not=i^{\prime}}s_{t,ii^{\prime}}$. In this case, $i$ and $i^{\prime}$ are incomparable in stage $t+1$ and constraints \eqref{eq:prior-xs} become redundant (i.e., $\sum_{\tau=1}^{t+1}x_{\tau, i}\ge \sum_{\tau=1}^{t+1}x_{\tau,i^{\prime}}-1$ and $\sum_{\tau=1}^{t+1}x_{\tau,i^{\prime}}\ge \sum_{\tau=1}^{t+1}x_{\tau,i}-1$). If $\sum_{\tau=1}^t x_{\tau,i}=\sum_{\tau=1}^tx_{\tau,i^{\prime}}=0$, constraints \eqref{eq:prior-staget} reduce to $s_{t,ii^{\prime}}+s_{t,i^{\prime}i}\ge 1$. This ensures that in each stage, we only prioritize on facilities that have not been opened yet. 
Moreover, constraints \eqref{eq:prior-budget} set the budget for opening facilities in each stage.

In Model \eqref{model:deterministic-prior}, the data we acquire in stage $t$ is the demand $\boldsymbol{d}_t\in\mathbb{R}_+^N$ and the budget $F_t\in\mathbb{R}_+$ for all $t=1,\ldots, T$. For notational simplicity, we denote $\boldsymbol{\xi}_t=(\boldsymbol{d}_t, F_t),\ \forall t=1,\ldots,T$. Suppose that the data series $\{\boldsymbol{\xi}_1,\boldsymbol{\xi}_2,\ldots,\boldsymbol{\xi}_T\}$ evolve according to a known probability distribution and $\boldsymbol{\xi}_1$ is also uncertain in this setting. We use a scenario tree to represent the decision-making process, which is described in Figure \ref{fig:scenario-tree-prior}. Note that because variable $\boldsymbol{s}_t$ is defined for stages $0,1,\ldots, T-1$ and variables $\boldsymbol{x}_t,\boldsymbol{y}_t$ are defined for stages $1,2,\ldots,T$, in the scenario-node-based notation (see Figure \ref{fig:scenario-tree-prior}(b)), $\boldsymbol{s}_n$ is defined for all nodes $n\in\mathcal{T}\setminus\mathcal{L}$ and $\boldsymbol{x}_n,\boldsymbol{y}_n$ are defined for all nodes $n\in\mathcal{T}\setminus\{1\}$. The definitions of the auxiliary variables $\eta_n,u_n$ in Figure \ref{fig:scenario-tree-prior} will be introduced later. To make the subproblem on each node feasible, we need to ensure that $\sum_{j=1}^Nd_{n,j}\le h_{t_n}\sum_{m\in\mathcal{P}(n)\setminus\{1\}}F_m$ for all $n\in\mathcal{T}\setminus\{1\}$.
\begin{figure}[ht!]
    \centering
    \begin{subfigure}{0.45\textwidth}
    \centering
    		\includegraphics[height = 6.5cm]{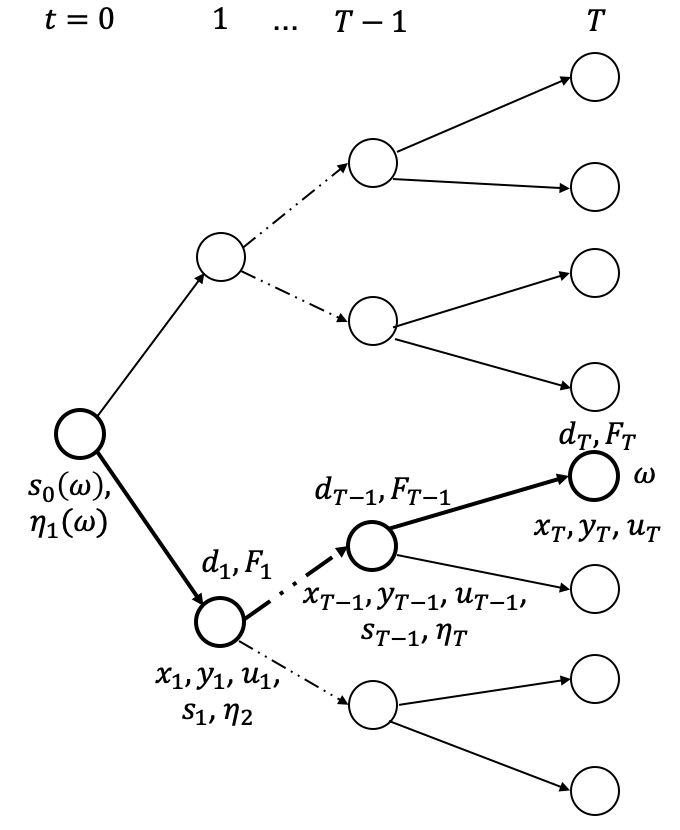}
        \caption{Scenario-path-based notation, where we omit the dependence of these variables on the scenario path $\omega$ for notational simplicity.}
    \end{subfigure}
	\begin{subfigure}{0.45\textwidth}
	\centering
			\includegraphics[height = 6.5cm]{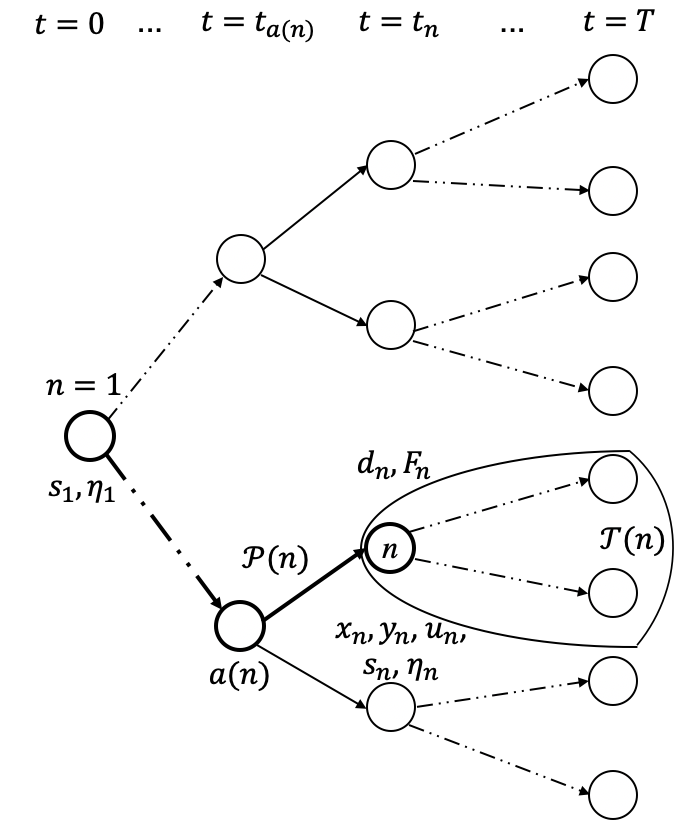}
        \caption{Scenario-node-based notation, where $\mathcal{P}(n)$ is the path from root node to node $n$, and $\mathcal{T}(n)$ is the subtree rooted at node $n$.}
	\end{subfigure}
    \caption{Illustration of a scenario tree and its related notation.}
    \label{fig:scenario-tree-prior}
\end{figure}

In the following Sections \ref{sec:two-stage-prior} and \ref{sec:multistage-prior}, we introduce risk-averse two-stage and multistage stochastic facility location models with prioritization, respectively.

\subsubsection{Risk-Averse Two-Stage Stochastic Facility Location with Prioritization}\label{sec:two-stage-prior}
In a risk-averse two-stage stochastic programming setting, the priority list is decided up front and fixed through all periods, i.e., $\boldsymbol{s}_t=\boldsymbol{s}_0,\ \forall t=0,1,\ldots,T-1$. As a result, constraints \eqref{eq:prior-staget} become redundant. We depict the decision-making process in the following flowchart:
\begin{eqnarray*}
\underbrace{\text{decision}\ (\boldsymbol{s}_0, \ldots, \boldsymbol{s}_{T-1})}_{\textcolor{black}{\text{Stage 1}}}\to\underbrace{\text{observation}\ (\boldsymbol{\xi}_1,\boldsymbol{\xi}_2,\ldots,\boldsymbol{\xi}_T)\to \text{decision}\ ((\boldsymbol{x}_t)_{t=1}^T,(\boldsymbol{y}_t)_{t=1}^T)}_{\textcolor{black}{\text{Stage 2}}}.
\end{eqnarray*}
After fixing the priority lists $\boldsymbol{s}_t,\ t=0,1,\ldots,T-1$, for each realized sample of the demand and budget, we select the facilities from the priority list until we exhaust the budget, and then calculate the operational cost $Q_t(\boldsymbol{s},\omega)$ for each $t=1,\ldots, T$. Similar to Section \ref{sec:risktwo}, we use a multiperiod risk function $\mathbb{F}^{TS}(Q_1,\ldots,Q_T)=\sum_{t=1}^T\rho_t(Q_t)$, which is the summation of the risk values in all periods.  Then, a scenario-path-based formulation of the two-stage risk-averse model with the risk measure $\mathbb{F}^{TS}(Q_1,\ldots,Q_T)$ can be written as follows:
\begin{align}
z_P^{TS}=\min_{\boldsymbol{s}_0,\cdots,\boldsymbol{s}_{T-1}}\quad &\sum_{t=0}^{T-1}\sum_{i\not=i^{\prime}}s_{t,ii^{\prime}}+\sum_{t=1}^T \rho_{t}\left(Q_t(\boldsymbol{s},\omega)\right)\label{model:twostagerisk-prior}\\
\text{s.t.}\quad& \text{\eqref{eq:prior-stage1}\ (constraints for $\boldsymbol{s}_0$)}\nonumber\\
& s_{t,ii^{\prime}}=s_{0,ii^{\prime}},\ \forall t=0,1,\ldots, T-1,\ 1\le i\not=i^{\prime}\le M\nonumber\\
&s_{t,ii^{\prime}}\in\{0,1\},\ \forall  1\le i\not=i^{\prime}\le M,\ 0\le t\le T-1,\nonumber
\end{align}
where for each $t=1,\ldots,T,\ \omega\in\Omega$, 
\begin{align*}
Q_t(\boldsymbol{s},\omega):=\min_{\boldsymbol{x}_t(\omega),\boldsymbol{y}_t(\omega)} \quad&\boldsymbol{c}_t^{\mathsf T} \boldsymbol{y}_t(\omega)\nonumber\\
\text{s.t.}\quad& (\boldsymbol{x}_t(\omega),\boldsymbol{y}_t(\omega))\in X_t(\omega) \nonumber\\
& s_{t,ii^{\prime}}+s_{t,i^{\prime}i} + \sum_{\tau=1}^t x_{\tau,i}(\omega) + \sum_{\tau=1}^tx_{\tau,i^{\prime}}(\omega)\ge 1,\ \forall 1\le i<i^{\prime}\le M,\ 1\le t\le T-1 \\
& \sum_{\tau=1}^tx_{\tau,i}(\omega) \ge \sum_{\tau=1}^tx_{\tau,i^{\prime}}(\omega) + s_{t-1, ii^{\prime}} - 1, \ \forall 1\le i\not=i^{\prime}\le M,\ 1\le t\le T\\
& \sum_{i=1}^M x_{t,i}(\omega) \le F_t(\omega),\ \forall 1\le t\le T \\
& \boldsymbol{x}_t(\omega)\in \mathbb{Z}^{M}_{+},\ \boldsymbol{y}_t(\omega)\in \mathbb{R}^{M\times N}_{+}. \nonumber
\end{align*}
Using the specific risk measure $\rho_t$ defined in \eqref{eq:rho} and auxiliary variables $({\eta}_n)_{n\in\mathcal{T}\setminus\{1\}}, (u_n)_{n\in\mathcal{T}\setminus\{1\}}$, a scenario-node-based formulation of the two-stage risk-averse model is given by
\begin{subequations}\label{model:two-prior-original}
\begin{align}
    z_P^{TS}=\min_{\substack{\boldsymbol{x}_n,\boldsymbol{y}_n, u_n, \eta_n, n\not=1\\\boldsymbol{s}_{n},n\not\in\mathcal{L}}} \quad& \sum_{n\in \mathcal{T}\setminus\{1\}} p_n\left((1-\lambda_{t_n}){\boldsymbol{c}}_{t_n}^{\mathsf T}\boldsymbol{y}_n + \frac{\lambda_{t_n}}{1-\alpha_{t_n}} u_n+{\lambda}_{t_n}\eta_n\right)+\sum_{n\in \mathcal{T}\setminus\mathcal{L}} p_n\sum_{i\not=i^{\prime}}s_{n,ii^{\prime}}\label{eq:two-prior-original-obj}\\
    \text{s.t.} \quad& s_{1,ii^{\prime}}+s_{1,i^{\prime}i}\ge 1,\ \forall 1\le i<i^{\prime}\le M \label{eq:two-prior-equal}\\
    & s_{n,ii^{\prime}}+s_{n,i^{\prime}i} + \sum_{m\in \mathcal{P}(n)\setminus\{1\}}x_{m,i} + \sum_{m\in \mathcal{P}(n)\setminus\{1\}}x_{m,i^{\prime}}\ge 1,\ \forall 1\le i<i^{\prime}\le M,\ n\not=1,\ n\not\in\mathcal{L} \label{eq:two-prior-sx}\\
& \sum_{m\in \mathcal{P}(n)\setminus\{1\}}x_{m,i} \ge \sum_{m\in \mathcal{P}(n)\setminus\{1\}}x_{m,i^{\prime}} + s_{a(n), ii^{\prime}} - 1, \ \forall 1\le i\not=i^{\prime}\le M,\ n\not=1\label{eq:two-prior-xs}\\
& s_{n,ii^{\prime}}\in\{0,1\},\ \forall  1\le i\not=i^{\prime}\le M,\ n\not\in\mathcal{L}\label{eq:two-prior-s}\\
& \sum_{i=1}^M x_{n,i} \le F_n,\ \forall n\in\mathcal{T},\ n\not=1\label{eq:two-prior-budget}\\
& (\boldsymbol{x}_n, \boldsymbol{y}_n)\in X_n,\ u_n \ge 0,\ \forall n\in \mathcal{T},\ n\not=1 \label{eq:two-prior-constraint_x}\\
& u_n + \eta_{n} \ge \boldsymbol{c}_n^{\mathsf T} \boldsymbol{y}_n,\ \forall n\in\mathcal{T}\setminus\{1\} \label{eq:two-prior-eta}\\
& s_{n,ii^{\prime}}=s_{1,ii^{\prime}},\ \forall 1\le i\not=i^{\prime}\le M,\ n\not\in\mathcal{L}\label{eq:two-prior-twostage-s-old}\\
& \eta_{n}=\eta_m,\ \forall n,m\in\mathcal{T}_t,\ 1\le t\le T\label{eq:two-prior-twostage-eta-old},
\end{align}
\end{subequations}
where constraints \eqref{eq:two-prior-twostage-s-old} and \eqref{eq:two-prior-twostage-eta-old} are the ``two-stage'' constraints to ensure that the priority lists $\boldsymbol{s}_n$ are identical for all the nodes, and the risk-related variables $\eta_n$ are identical for all the nodes in the same stage. We denote constraints \eqref{eq:two-prior-equal}--\eqref{eq:two-prior-s} as $(\boldsymbol{s}, \boldsymbol{x})\in S$.

\subsubsection{Risk-Averse Multistage Stochastic Facility Location with Prioritization}\label{sec:multistage-prior}
In a risk-averse multistage stochastic dynamic programming framework,  we consider the following $T+1$-stage decision process: 
\begin{align*}
&\underbrace{\text{decision}\ (\boldsymbol{s}_0)}_{\text{Stage 0}}\to \underbrace{\text{observation}\ (\boldsymbol{d}_1, F_1)\to \text{decision}\ (\boldsymbol{x}_1, \boldsymbol{y}_1, \boldsymbol{s}_1)}_{\text{Stage 1}}\to \text{observation}\ (\boldsymbol{d}_2, F_2)\\
\to& \cdots\to \text{decision}\ (\boldsymbol{x}_{T-1}, \boldsymbol{y}_{T-1},\boldsymbol{s}_{T-1})\to \underbrace{\text{observation}\ (\boldsymbol{d}_{T},F_T)\to \text{decision}\ (\boldsymbol{x}_T, \boldsymbol{y}_T)}_{\text{Stage T}}.
\end{align*}
Before realizing the uncertainty $\boldsymbol{\xi}_t=(\boldsymbol{d}_t, F_t)$, we decide a priority list $\boldsymbol{s}_{t-1}$, and after observing a realization of $\boldsymbol{\xi}_t=(\boldsymbol{d}_t, F_t)$, we choose facilities from the top of the priority list $\boldsymbol{s}_{t-1}$ until we exhaust the budget and then update the priority list $\boldsymbol{s}_t$ to include all unopened facility sites for all stages $t=1,2,\ldots,T$. In the last stage $T$, we do not need to update the priority list $\boldsymbol{s}_T$ because there are no subsequent decisions to be made. Similar to Section \ref{sec:riskmulti}, we denote $\boldsymbol{\xi}_{[t]}=(\boldsymbol{\xi}_1,\ldots,\boldsymbol{\xi}_t)$.

We denote the stagewise cost by $g_t,\ t=0,1,\ldots,T$, where $g_0 = \sum_{i\not=i^{\prime}}s_{0,ii^{\prime}}$, $g_t=\sum_{i\not=i^{\prime}}s_{t,ii^{\prime}}+\boldsymbol{c}_t^{\mathsf T}\boldsymbol{y}_t$ for all $t = 1,\ldots, T-1$, and $g_T=\boldsymbol{c}_T^{\mathsf T}\boldsymbol{y}_T$. Note that only $g_0$ is deterministic, while $g_t,\ t=1,\ldots,T$ are stochastic with respect to the dynamic uncertainty realization process. We extend the ECRMs \eqref{ECRMs} to handle $(T+1)$-stage costs as below: 
\begin{align*}
\mathbb{F}^{MS}(g_0,g_1,\ldots,g_{T})=g_0+\rho_1(g_1)+\mathbb{E}_{\boldsymbol \xi_{1}}\left[{\rho_2^{\boldsymbol\xi_{1}}}(g_2)\right]+\mathbb{E}_{\boldsymbol \xi_{[2]}}\left[{\rho_3^{\boldsymbol \xi_{[2]}}}(g_3)\right]+\cdots+\mathbb{E}_{\boldsymbol \xi_{[T-1]}}\left[{\rho_{T}^{\boldsymbol\xi_{[T-1]}}}(g_{T})\right],
\end{align*}
where the expectation starts from $\mathbb{E}_{\boldsymbol{\xi}_1}$ because $\boldsymbol{\xi}_1$ is also uncertain as defined in Section \ref{sec:problem-prior}. 
Then, a nested risk-averse multistage stochastic programming model using $\mathbb{F}^{MS}(g_0,g_1,\ldots,g_{T})$ can be written as:
\footnotesize
\begin{align*}
z_P^{MS}=\min_{\substack{\boldsymbol{x}_t,\boldsymbol{y}_t,u_t,\\1\le t\le T\\\boldsymbol{s}_t,\eta_{t+1},\\0\le t\le T-1}}\quad&\sum_{i\not=i^{\prime}}s_{0,ii^{\prime}}+\min_{\eta_1}\lambda_1\eta_1+\mathbb{E}_{\boldsymbol{\xi}_1 }\left[\frac{\lambda_1}{1-\alpha_1}u_1+(1-\lambda_1)\left(\sum_{i\not=i^{\prime}}s_{1,ii^{\prime}}+\boldsymbol{c}_1^{\mathsf T}\boldsymbol{y}_1\right)\right]\nonumber\\
&+\mathbb{E}_{\boldsymbol{\xi}_1}\left[\min_{\eta_2}\lambda_2\eta_2+\mathbb{E}_{\boldsymbol{\xi}_2{|\boldsymbol\xi_1}}\left[\frac{\lambda_2}{1-\alpha_2}u_2+(1-\lambda_2)\left(\sum_{i\not=i^{\prime}}s_{2,ii^{\prime}}+\boldsymbol{c}_2^{\mathsf T}\boldsymbol{y}_2\right)\right]+\cdots\right.\nonumber\\
&+\mathbb{E}_{\boldsymbol{\xi}_{T-2}{|\boldsymbol\xi_{[T-3]}}}\left[\min_{\eta_{T-1}}\lambda_{T-1}\eta_{T-1}+\mathbb{E}_{\boldsymbol{\xi}_{T-1}{|\boldsymbol\xi_{[{T-2}]}}}\left[\frac{\lambda_{T-1}}{1-\alpha_{T-1}}u_{T-1}+(1-\lambda_{T-1})\left(\sum_{i\not=i^{\prime}}s_{T-1,ii^{\prime}}+\boldsymbol{c}_{T-1}^{\mathsf T}\boldsymbol{y}_{T-1}\right)\right]\right.\nonumber\\
&\left.\left.+\mathbb{E}_{\boldsymbol{\xi}_{T-1}{|\boldsymbol\xi_{[T-2]}}}\left[\min_{\eta_T}\lambda_T\eta_T+\mathbb{E}_{\boldsymbol{\xi}_T{|\boldsymbol\xi_{[T-1]}}}\left[\frac{\lambda_T}{1-\alpha_T}u_T+(1-\lambda_T)\boldsymbol{c}_T^{\mathsf T}\boldsymbol{y}_T\right]\right]\cdots\right]\right].
\end{align*}
\normalsize
Correspondingly, we present a scenario-node-based formulation as follows:
\begin{subequations}\label{model:multi-prior}
\begin{align}
z^{MS}_P=\min_{\substack{\boldsymbol{x}_n,\boldsymbol{y}_n, u_n, n\not=1\\\boldsymbol{s}_{n},\eta_{n},n\not\in\mathcal{L}}} \quad& \sum_{n\in \mathcal{T}} p_n\left(\tilde{\boldsymbol{c}}_n^{\mathsf T}\boldsymbol{y}_n + \tilde{\lambda}_n\eta_n + \tilde{\alpha}_n u_n+\tilde{1}_n\sum_{i\not=i^{\prime}}s_{n,ii^{\prime}}\right)\label{eq:multi-prior-obj}\\
\text{s.t.} \quad& \text{\eqref{eq:two-prior-equal}--\eqref{eq:two-prior-constraint_x}}\nonumber\\
& u_n + \eta_{a(n)} \ge \boldsymbol{c}_{t_n}^{\mathsf T} \boldsymbol{y}_n,\ \forall n\in\mathcal{L} \label{eq:multi-prior-etaL}\\
& u_n + \eta_{a(n)} \ge \boldsymbol{c}_{t_n}^{\mathsf T} \boldsymbol{y}_n + \sum_{i\not=i^{\prime}}s_{n,ii^{\prime}},\ \forall n\in \mathcal{T},\ n\not=1,\ n\not\in\mathcal{L}, \label{eq:multi-prior-eta}
\end{align}
\end{subequations}
where $\tilde{\boldsymbol{c}}_n=0$ if $n=1$, $\tilde{\boldsymbol{c}}_n=(1-\lambda_{t_n})\boldsymbol{c}_{t_n}$ otherwise; $\tilde{\lambda}_n=0$ if $n\in\mathcal{L}$, $\tilde{\lambda}_n=\lambda_{t_n+1}$ otherwise; $\tilde{\alpha}_n=0$ if $n=1$, $\tilde{\alpha}_n=\frac{\lambda_{t_n}}{1-\alpha_{t_n}}$ otherwise; $\tilde{1}_n=1$ if $n=1$, $\tilde{1}_n=0$ if $n\in\mathcal{L}$, and $\tilde{1}_n=(1-\lambda_{t_n})$ otherwise. Note that variables $\boldsymbol{s}_n, \eta_n$ are defined for all non-leaf nodes, and variables $\boldsymbol{x}_n,\boldsymbol{y}_n,u_n$ are defined for all non-root nodes.
Comparing Models \eqref{model:two-prior-original} with \eqref{model:multi-prior}, the differences are in the objective function and constraints \eqref{eq:two-prior-eta} and \eqref{eq:multi-prior-etaL}--\eqref{eq:multi-prior-eta}. In the next lemma, we prove that the risk-averse two-stage model \eqref{model:two-prior-original} can be recast in the following way:
\begin{subequations}\label{model:two-prior}
\begin{align}
    z^{TS}_P=\min_{\substack{\boldsymbol{x}_n,\boldsymbol{y}_n, u_n, n\not=1\\\boldsymbol{s}_{n},\eta_{n},n\not\in\mathcal{L}}} \quad& \sum_{n\in \mathcal{T}} p_n\left(\tilde{\boldsymbol{c}}_n^{\mathsf T}\boldsymbol{y}_n + \tilde{\lambda}_n\eta_n + \tilde{\alpha}_n u_n+\tilde{1}_n\sum_{i\not=i^{\prime}}s_{n,ii^{\prime}}\right)\nonumber\\
\text{s.t.} \quad& \text{\eqref{eq:two-prior-equal}--\eqref{eq:two-prior-constraint_x}}\nonumber\\
&\text{\eqref{eq:multi-prior-etaL}--\eqref{eq:multi-prior-eta}}\nonumber\\
& s_{n,ii^{\prime}}=s_{1,ii^{\prime}},\ \forall 1\le i\not=i^{\prime}\le M,\ n\not\in\mathcal{L}\label{eq:two-prior-twostage-s}\\
& \eta_{n}=\eta_m,\ \forall n,m\in\mathcal{T}_t,\ 0\le t\le T-1\label{eq:two-prior-twostage-eta}.
\end{align}
\end{subequations}

\begin{lemma}\label{lemma:prior}
The risk-averse two-stage model \eqref{model:two-prior-original} is equivalent to Model \eqref{model:two-prior}.
\end{lemma}

The detailed proof of Lemma \ref{lemma:prior} is presented in Appendix \ref{e-companion:proofs}. From here, we observe that ${\rm VMS_P}=z^{TS}_P - z^{MS}_P\ge 0$. In the next section, we examine a substructure problem, based on which we derive an analytical lower bound for ${\rm VMS_P}$.

\subsection{Analytical Solutions of the Substructure Problem}\label{sec:substructure-prior}
We first examine a substructure of Models \eqref{model:multi-prior} and \eqref{model:two-prior} once we fix the values of the $(\boldsymbol{x},\boldsymbol{y},\boldsymbol{u})$-variables. We denote the resultant problems with known $(\boldsymbol{x}_n^*,\boldsymbol{y}_n^*,{u}_n^*)$-values as $\mbox{{\bf SP-PMS}}(\boldsymbol{x}_n^*,\boldsymbol{y}_n^{*}, u_n^{*})$ and $\mbox{{\bf SP-PTS}}(\boldsymbol{x}_n^*,\boldsymbol{y}_n^{*}, u_n^{*})$, which are respectively defined as follows:
\begin{subequations}\label{eq:ms-substructure-prior}
\begin{align}
\mbox{{\bf SP-PMS}}(\boldsymbol{x}_n^*,\boldsymbol{y}_n^{*}, u_n^{*}):
\min_{\substack{\boldsymbol{s}_n, \eta_n,n\not\in\mathcal{L}}} \quad& \sum_{n\in \mathcal{T}} p_n\left(\tilde{\lambda}_n\eta_n + \tilde{1}_n\sum_{i\not=i^{\prime}}s_{n,ii^{\prime}}\right)\\
\text{s.t.}\quad& (\boldsymbol{s}, \boldsymbol{x}^*)\in S\label{eq:multi-prior-S-sub}\\
& u_n^* + \eta_{a(n)} \ge \boldsymbol{c}_{t_n}^{\mathsf T} \boldsymbol{y}_n^*,\ \forall n\in\mathcal{L} \label{eq:multi-prior-etaL*}\\
& u_n^* + \eta_{a(n)} \ge \boldsymbol{c}_{t_n}^{\mathsf T} \boldsymbol{y}_n^* + \sum_{i\not=i^{\prime}}s_{n,ii^{\prime}},\ \forall n\in \mathcal{T},\ n\not=1,\ n\not\in\mathcal{L} \label{eq:multi-prior-eta*}
\end{align}
\end{subequations}
and
\begin{align}
\mbox{{\bf SP-PTS}}(\boldsymbol{x}_n^*,\boldsymbol{y}_n^{*}, u_n^{*}):\ \min_{\substack{\boldsymbol{s}_n,\eta_n,n\not\in\mathcal{L}}}\quad & \sum_{n\in \mathcal{T}} p_n\left(\tilde{\lambda}_n\eta_n + \tilde{1}_n\sum_{i\not=i^{\prime}}s_{n,ii^{\prime}}\right) \label{eq:ts-substructure-prior}\\
\text{s.t.}\quad& \text{\eqref{eq:multi-prior-S-sub}--\eqref{eq:multi-prior-eta*}}\nonumber\\
&\text{\eqref{eq:two-prior-twostage-s}--\eqref{eq:two-prior-twostage-eta}\ (Two-stage constraints for $\boldsymbol{s}$ and $\boldsymbol{\eta}$)}.\nonumber
\end{align}
We denote the optimal objective values of Models \eqref{eq:ms-substructure-prior} and \eqref{eq:ts-substructure-prior} as $Q_P^M(\boldsymbol{x}_n^{*}, \boldsymbol{y}_n^{*}, u_n^{*}),\ Q_P^{T}(\boldsymbol{x}_n^{*}, \boldsymbol{y}_n^{*}, u_n^{*})$, respectively. The next proposition gives the analytical forms of the optimal solutions to $\mbox{{\bf SP-PMS}}(\boldsymbol{x}_n^{*},\boldsymbol{y}_n^{*}, u_n^{*})$ and $\mbox{{\bf SP-PTS}}(\boldsymbol{x}_n^{*},\boldsymbol{y}_n^{*}, u_n^{*})$, of which a detailed proof is presented in Appendix \ref{e-companion:proofs}.

\begin{proposition}\label{prop:substructure-prior}
Given a feasible solution $(\boldsymbol{x}_n^{*}, \boldsymbol{y}_n^{*}, u_n^{*})_{n\not=1}$ to Model \eqref{model:two-prior}, the optimal solutions $(\boldsymbol{s}_n^{MS},{\eta}_n^{MS})_{n\not\in\mathcal{L}}$ of \eqref{eq:ms-substructure-prior} have the following analytical forms: for all $n\in\mathcal{T}\setminus\mathcal{L}$, we have
\begin{itemize}
    \item if $\sum_{m\in\mathcal{P}(n)\setminus\{1\}}x^{*}_{m,i} = 1$ or $\sum_{m\in\mathcal{P}(n)\setminus\{1\}}x^{*}_{m,i^{\prime}} = 1$, then $s^{MS}_{n,ii^{\prime}}=s^{MS}_{n,i^{\prime}i}=0$;
    \item if $\sum_{m\in\mathcal{P}(n)\setminus\{1\}}x^{*}_{m,i} =\sum_{m\in\mathcal{P}(n)\setminus\{1\}}x^{*}_{m,i^{\prime}} = 0$ and there exists $m\in\mathcal{C}(n)$ such that $\sum_{l\in\mathcal{P}(m)\setminus\{1\}}x^{*}_{l,i} = 0,\  \sum_{l\in\mathcal{P}(m)\setminus\{1\}}x^{*}_{l,i^{\prime}}=1$, then $s^{MS}_{n,ii^{\prime}}=0,\ s^{MS}_{n,i^{\prime}i}=1$;
    \item if $\sum_{m\in\mathcal{P}(n)\setminus\{1\}}x^{*}_{m,i} =\sum_{m\in\mathcal{P}(n)\setminus\{1\}}x^{*}_{m,i^{\prime}} = 0$ and for all nodes $m\in\mathcal{C}(n)$ we have either $\sum_{l\in\mathcal{P}(m)\setminus\{1\}}x^{*}_{l,i} = \sum_{l\in\mathcal{P}(m)\setminus\{1\}}x^{*}_{l,i^{\prime}}=1$ or $\sum_{l\in\mathcal{P}(m)\setminus\{1\}}x^{*}_{l,i} = \sum_{l\in\mathcal{P}(m)\setminus\{1\}}x^{*}_{l,i^{\prime}}=0$, then we randomly pick one of $s^{MS}_{n,ii^{\prime}},\ s^{MS}_{n,i^{\prime}i}$ to be 1 and set the other to 0,
\end{itemize}
and
\begin{subequations}\label{eq:construct_eta_ms}
\begin{align}
    &\eta_n^{MS}=\max_{m\in\mathcal{C}(n)}\{\boldsymbol{c}_{t_m}^{\mathsf T} \boldsymbol{y}^{*}_m + \sum_{i\not=i^{\prime}}s^{MS}_{m,ii^{\prime}}-u_m^{*}\},\ \forall n\in\mathcal{T}_t,\ t=0,\ldots, T-2\\
    &\eta_n^{MS}=\max_{m\in\mathcal{C}(n)}\{\boldsymbol{c}_{t_m}^{\mathsf T} \boldsymbol{y}^{*}_m -u_m^{*}\},\ \forall n\in\mathcal{T}_{T-1}
\end{align}
\end{subequations}
while the optimal solutions $(\boldsymbol{s}_n^{TS},{\eta}_n^{TS})_{n\not\in\mathcal{L}}$ of \eqref{eq:ts-substructure-prior} have the following analytical forms: 
\begin{itemize}
    \item if there exists $m\in\mathcal{T}\setminus\{1\}$ such that $\sum_{l\in\mathcal{P}(m)\setminus\{1\}}x^{*}_{l,i} = 0,\  \sum_{l\in\mathcal{P}(m)\setminus\{1\}}x^{*}_{l,i^{\prime}}=1$, then for all $n\in\mathcal{T}\setminus\mathcal{L}$, we have $s^{TS}_{n,ii^{\prime}}=0,\ s^{TS}_{n,i^{\prime}i}=1$;
    \item if for all nodes $m\in\mathcal{T}\setminus\{1\}$ we have either $\sum_{l\in\mathcal{P}(m)\setminus\{1\}}x^{*}_{l,i} = \sum_{l\in\mathcal{P}(m)\setminus\{1\}}x^{*}_{l,i^{\prime}}=1$ or $\sum_{l\in\mathcal{P}(m)\setminus\{1\}}x^{*}_{l,i} = \sum_{l\in\mathcal{P}(m)\setminus\{1\}}x^{*}_{l,i^{\prime}}=0$, then we randomly pick one of $s^{TS}_{n,ii^{\prime}},\ s^{TS}_{n,i^{\prime}i}$ to be 1 and set the other to 0,
\end{itemize}
and 
\begin{subequations}\label{eq:construct_eta_ts}
\begin{align}
    &\eta_n^{TS}=\max_{m\in\mathcal{T}_{t_n+1}}\{\boldsymbol{c}_{t_m}^{\mathsf T} \boldsymbol{y}^{*}_m + \sum_{i\not=i^{\prime}}s^{TS}_{m,ii^{\prime}}-u_m^{*}\},\ \forall n\in\mathcal{T}_t,\ t=0,\ldots, T-2\\
    &\eta_n^{TS}=\max_{m\in\mathcal{T}_{t_n+1}}\{\boldsymbol{c}_{t_m}^{\mathsf T} \boldsymbol{y}^{*}_m -u_m^{*}\},\ \forall n\in\mathcal{T}_{T-1}.
\end{align}
\end{subequations}
Correspondingly, we have $
    Q_P^M(\boldsymbol{x}_n^{*}, \boldsymbol{y}_n^{*}, u_n^{*}) = \sum_{n\in \mathcal{T}} p_n\left(\tilde{\lambda}_n\eta^{MS}_n + \tilde{1}_n\sum_{i\not=i^{\prime}}s^{MS}_{n,ii^{\prime}}\right),\ 
    Q_P^T(\boldsymbol{x}_n^{*}, \boldsymbol{y}_n^{*}, u_n^{*}) = \sum_{n\in \mathcal{T}} p_n\left(\tilde{\lambda}_n\eta^{TS}_n + \tilde{1}_n\sum_{i\not=i^{\prime}}s^{TS}_{n,ii^{\prime}}\right).
$
\end{proposition}

\subsection{${\rm VMS_P}$ for the Risk-Averse Facility Location Problem with Prioritization}\label{sec:proof-VMS-prior}
We now describe a lower bound on the ${\rm VMS_P}$ for the risk-averse multistage and two-stage facility location models \eqref{model:multi-prior} and \eqref{model:two-prior} based on the analysis in the previous section.
\begin{theorem}
\label{thm:risk-averse-bound-prior}
	Let $(\boldsymbol{x}_n^{*}, \boldsymbol{y}_n^{*}, u_n^{*})_{n\not=1}$ be the second-stage decisions in an optimal solution to the two-stage model \eqref{model:two-prior}, and let $\boldsymbol{s}_n^{TS}, \eta_n^{TS}, \boldsymbol{s}_n^{MS}, \eta_n^{MS}$ be defined in Proposition \ref{prop:substructure-prior}, which are constructed by $(\boldsymbol{x}_n^{*}, \boldsymbol{y}_n^{*}, u_n^{*})_{n\not=1}$.
	Then,
	\begin{align}
	{\rm VMS_P}\ge & \sum_{n\in \mathcal{T}\setminus\mathcal{L}} p_n{\lambda}_{t_n+1}(\eta_n^{TS}-\eta_n^{MS})+ \sum_{n\in \mathcal{T}\setminus(\{1\}\cup\mathcal{L})} p_n (1-\lambda_{t_n})(\sum_{i\not=i^{\prime}}s^{TS}_{n,ii^{\prime}}-\sum_{i\not=i^{\prime}}s^{MS}_{n,ii^{\prime}})\label{eq:VMS_P_LB}\\
	\ge &\sum_{n\in \mathcal{T}\setminus\mathcal{L}} p_n{\lambda}_{t_n+1} \left(\max_{m\in\mathcal{T}_{t_n+1}}\{\boldsymbol{c}_{t_m}^{\mathsf T}\boldsymbol{y}_m^{*}-u_m^{*}\}-\max_{m\in\mathcal{C}(n)}\{\boldsymbol{c}_{t_m}^{\mathsf T}\boldsymbol{y}_m^{*}-u_m^{*}\}\right),\nonumber
	\end{align}
	where we denote the right-hand side of \eqref{eq:VMS_P_LB} as ${\rm VMS_P^{LB}}$.
\end{theorem}

\begin{remark}
    The lower bound ${\rm VMS_P^{LB}}$ provided in Theorem \ref{thm:risk-averse-bound-prior} is \textit{tight} (i.e., it cannot be improved any further).
    In Appendix \ref{e-companion:example}, we will provide Example \ref{eg2} where the equality holds (i.e., ${\rm VMS_P} = {\rm VMS_P^{LB}}$). 
\end{remark}

\begin{remark}\label{remark:parameter}
    To derive a lower bound that is directly related to the uncertain parameters $\boldsymbol{\xi}_t$, we consider a large enough confidence level $\alpha_t$ such that the $\rm{VaR}_{\alpha_t}$ is the maximum value. For example, when we set $\alpha_t=0.95$ and the number of scenarios $|\Omega|<20$, $\rm{VaR}_{\alpha_t}$ is the maximum value among all scenarios $\omega\in\Omega$ and correspondingly, $u_n^*=0,\ \forall n\not=1$. Recall that $c_{t,\rm min}=\min_{i\in[M],j\in[N]}c_{tij},\ c_{t,\rm max}=\max_{i\in[M],j\in[N]}c_{tij},\ \forall t\in[T]$. Then, we have
    \begin{align}
    {\rm VMS_P}\ge & \sum_{n\in \mathcal{T}\setminus\mathcal{L}} p_n{\lambda}_{t_n+1} \left(\max_{m\in\mathcal{T}_{t_n+1}}\{\boldsymbol{c}_{t_m}^{\mathsf T}\boldsymbol{y}_m^{*}-u_m^{*}\}-\max_{m\in\mathcal{C}(n)}\{\boldsymbol{c}_{t_m}^{\mathsf T}\boldsymbol{y}_m^{*}-u_m^{*}\}\right)\nonumber\\
    \ge & \sum_{n\in \mathcal{T}\setminus\mathcal{L}} p_n{\lambda}_{t_n+1} \left(c_{t_n+1,\rm min}\max_{m\in\mathcal{T}_{t_n+1}}\{\sum_{j\in[N]}d_{m,j}\}-c_{t_n+1,\rm max}\max_{m\in\mathcal{C}(n)}\{\sum_{j\in[N]}d_{m,j}\}\right)\nonumber\\
    \overset{(a)}{=} & \sum_{n\in \mathcal{T}\setminus\mathcal{L}} p_n{\lambda}_{t_n+1}{c}_{t_n+1} \left(\max_{m\in\mathcal{T}_{t_n+1}}\{\sum_{j\in[N]}d_{m,j}\}-\max_{m\in\mathcal{C}(n)}\{\sum_{j\in[N]}d_{m,j}\}\right),\label{eq:VMS_P_LB1}
    \end{align}
    where we denote the right-hand side of \eqref{eq:VMS_P_LB1} as ${\rm VMS_P^{LB1}}$, and $(a)$ holds when $c_{t,\rm min}=c_{t,\rm max},\ \forall t\in[T]$.
    Note that ${\rm VMS_P^{LB1}}$ is directly related to the variation of the uncertain demand in the scenario tree.
\end{remark}

\subsection{Prioritization Cuts}\label{sec:prioritization-cuts}
To improve the computational time, we develop a set of cutting planes for Models \eqref{model:multi-prior} and \eqref{model:two-prior} in Theorems \ref{thm:prior-cuts} and \ref{thm:prior-two-stage-cuts}, respectively, named the \textit{prioritization cuts}. The detailed proofs are presented in Appendix \ref{e-companion:proofs}.

\begin{theorem}\label{thm:prior-cuts}
There exists an optimal solution $(\bar{\boldsymbol{s}}_n, \bar{\boldsymbol{x}}_n, \bar{\eta}_n, \bar{\boldsymbol{y}}_n, \bar{u}_n)_{n\in\mathcal{T}}$ to model \eqref{model:multi-prior} such that $\bar{\boldsymbol{x}}_n$ satisfies the set of inequalities
\begin{align}
    \sum_{m\in\mathcal{P}(n)\setminus\{1\}}\bar{x}_{m,i}\ge \sum_{m\in\mathcal{P}(n^{\prime})\setminus\{1\}}\bar{x}_{m,i},\ \forall i=1,\ldots, M, \label{eq:prioritization-cuts-multi}
\end{align}
for all $n\not=1, n^{\prime}\in\mathcal{C}(a(n))$ such that $F_n\ge F_{n^{\prime}}$.
\end{theorem}

\begin{theorem}\label{thm:prior-two-stage-cuts}
There exists an optimal solution $(\bar{\boldsymbol{s}}_n, \bar{\boldsymbol{x}}_n, \bar{\eta}_n, \bar{\boldsymbol{y}}_n, \bar{u}_n)_{n\in\mathcal{T}}$ to model \eqref{model:two-prior} such that $\bar{\boldsymbol{x}}_n$ satisfies the set of inequalities
\begin{align}
    \sum_{m\in\mathcal{P}(n)\setminus\{1\}}\bar{x}_{m,i}\ge \sum_{m\in\mathcal{P}(n^{\prime})\setminus\{1\}}\bar{x}_{m,i},\ \forall i=1,\ldots, M, \label{eq:prioritization-cuts-two}
\end{align}
for all $n\not=1, n^{\prime}\in\mathcal{C}(a(n))$ such that $F_n\ge F_{n^{\prime}}$.
\end{theorem}

We do not use the term \textit{valid inequality} for prioritization cuts \eqref{eq:prioritization-cuts-multi} and \eqref{eq:prioritization-cuts-two} because they may rule out some feasible (even optimal) solutions. However, we ensure that there remains at least one optimal solution that satisfies the proposed cuts. \cite{israeli2002shortest} refer to such inequalities as \textit{super-valid inequalities}. In our experiments, we extract the prioritization cuts \eqref{eq:prioritization-cuts-multi} and \eqref{eq:prioritization-cuts-two} from the problem data before we solve the problem and add all of them to models \eqref{model:multi-prior} and \eqref{model:two-prior}. To avoid increasing the problem size, one can also place the prioritization cuts in a pool and iteratively add those that are violated by the LP relaxation solution in the branch-and-bound tree.
}

\section{Computational Results}
\label{sec:compu}
We test the risk-averse two-stage and multistage models {with or without prioritization} on two types of networks -- a randomly generated grid network where we vary the parameter settings extensively and a real-world network based on the United States map with 49 candidate facilities and 88 customer sites \citep{daskin2011network}.
Specifically, we conduct sensitivity analysis and report results based on the synthetic data to illustrate the tightness of the analytical bound and the efficacy and efficiency of the proposed approximation algorithms {and prioritization cuts} in Section \ref{sec:synthetic}. We also conduct a case study on the United States map-based network to display the solution patterns under different settings of uncertainties in Section \ref{sec:real}.
We use Gurobi 9.0.3 coded in Python 3.6.8 for solving all mixed-integer programming models, where the computational time limit is set to one hour. Our numerical tests are conducted on a Windows 2012 Server with 128 GB RAM and an Intel 2.2 GHz processor. 

\subsection{Result Analysis on Synthetic Data}
\label{sec:synthetic}
 We first introduce the experimental design and setup in Section \ref{sec:setup}, and report sensitivity analysis results in Section \ref{sec:sensi}, which are based on in-sample objective values. In Section \ref{sec:risk-profile}, we conduct out-of-sample test in a rolling horizon way to evaluate the risk profile of the two-stage and multistage solutions.
 Then we examine how tight the analytical bounds derived in Theorems~\ref{thm:risk-averse-bound} {and \ref{thm:risk-averse-bound-prior}} are in Section \ref{sec:tight}, and the performance of the approximation algorithms on the synthetic data set in Section \ref{sec:aa}, respectively. {Finally, in Section \ref{sec:time}, we report the computational time for the two risk-averse models with or without prioritization.} To compare the two-stage and multistage models {without prioritization}, we define the relative value of risk-averse multistage stochastic programming as ${\rm RVMS_R} = \frac{z_R^{TS} - z_R^{MS}}{z_R^{MS}}$, and the relative gap of the analytical bounds as ${\rm RGAP_R} = \frac{{\rm VMS_R} - {\rm VMS_R^{LB}}}{{\rm VMS_R}}$, respectively. {In parallel, we define the relative value of risk-averse multistage stochastic programming with prioritization as ${\rm RVMS_P} = \frac{z_P^{TS} - z_P^{MS}}{z_P^{MS}}$, and the relative gap of the analytical bounds as ${\rm RGAP_P} = \frac{{\rm VMS_P} - {\rm VMS_P^{LB}}}{{\rm VMS_P}}$, respectively.}

\subsubsection{Experimental Design and Setup}
\label{sec:setup}
We randomly sample $M$ potential facilities and $N$ customer sites on a $100\times 100$ grid and in the default setting, we have the number of stages ($T$) being 3, the number of facilities ($M$) being 6, the number of customer sites ($N$) being 10, the number of branches in each non-leaf node ($C$) being 2. The risk attitude parameters are set to $\lambda_t=0.5,\ \alpha_t=0.95,\ \forall t=2,\ldots, T$ at default. The operational costs between facilities and customer sites are calculated by their Manhattan distances times the unit travel cost.
We set the per stage renting costs $f_{ti}=6\times 10^4$ and all the facilities have the same capacity $h_{ti}=h=10^5$ for all $t=1,\ldots,T,\ i=1,\ldots,M$. For each customer site $j=1,\ldots, N$ and stage $t=1,\ldots,T$, we uniformly sample the demand mean from $U(1000(2t-1),5000(2t-1))$ and then multiply each mean by a fixed number ($\sigma=0.8$ at default) to generate its demand standard deviation. Lastly, we sample demand data following a truncated Normal distribution with the generated mean and standard deviation, while negative demand values are deleted. {For the models with prioritization, to keep the problem feasible, we uniformly sample values of budget $F_t$ from $U(\lceil\max_{\omega\in\Omega}\{\sum_{j\in[N]}d_{t,j}(\omega)\}/h\rceil-\lceil\max_{\omega\in\Omega}\{\sum_{j\in[N]}d_{t-1,j}(\omega)\}/h\rceil, \lceil\max_{\omega\in\Omega}\{\sum_{j\in[N]}d_{t,j}(\omega)\}/h\rceil-\lceil\max_{\omega\in\Omega}\{\sum_{j\in[N]}d_{t-1,j}(\omega)\}/h\rceil+2)$, where $d_{0,j}=0,\ \forall j\in[N]$ and $\lceil\max_{\omega\in\Omega}\{\sum_{j\in[N]}d_{t,j}(\omega)\}/h\rceil-\lceil\max_{\omega\in\Omega}\{\sum_{j\in[N]}d_{t-1,j}(\omega)\}/h\rceil$ represents at least how many new facilities are needed to cover the demand in stage $t$. We consider three types of scenario trees listed below:
\begin{itemize}
    \item Stagewise dependent (SD): at every stage $t=1,\ldots,T-1$, every node $n\in\mathcal{T}_t$ is associated with a different set of children nodes $\mathcal{C}(n)$, i.e., $\mathcal{C}(n)\not = \mathcal{C}(m),\ \forall n,m\in\mathcal{T}_t$;
    \item Stagewise independent (SI): at every stage $t=1,\ldots,T-1$, every node $n\in\mathcal{T}_t$ is associated with an identical set of children nodes $\mathcal{C}(n)$, i.e., $\mathcal{C}(n) = \mathcal{C}(m),\ \forall n,m\in\mathcal{T}_t$;
    \item Stagewise dependent with Scenario 0 (SD0): at every stage $t=1,\ldots,T-1$, every node $n\in\mathcal{T}_t$ is associated with a different set of children nodes $\mathcal{C}(n)$, which includes realization $d_{m,j}=0,\ \forall j\in[N]$ in one of the children nodes $m\in\mathcal{C}(n),\ n\in\mathcal{T}_t$.
\end{itemize}
Here, SD represents the most general case where in each stage $t$, we have at most $C^{t}$ different realizations of the uncertainty $\boldsymbol{\xi}_t$; SI assumes that the stochastic process $(\boldsymbol{\xi}_1,\boldsymbol{\xi}_2,\ldots,\boldsymbol{\xi}_T)$ is stagewise independent and thus we have at most $C$ different realizations of the uncertainty in each stage $t$; and SD0 is a special case to illustrate how ${\rm VMS_R}$ depends on the number of nodes $n$ with realization $d_{n,j}=0,\ \forall j\in[N]$ as stated in Corollary \ref{cor:parameter}. To illustrate what types of demand scenarios make the multistage model more valuable, we compare these three types of scenario trees to evaluate ${\rm RVMS_R}$, namely SD-${\rm RVMS_R}$, SI-${\rm RVMS_R}$, and SD0-${\rm RVMS_R}$, respectively. We also compare them to the risk-averse models with prioritization using SD scenario trees, namely SD-${\rm RVMS_P}$.
}

\subsubsection{Sensitivity Analysis on ${\rm RVMS_R}$ {and ${\rm RVMS_P}$}}
\label{sec:sensi}
\begin{figure}[ht!]
    \centering
    \begin{subfigure}{0.45\textwidth}
 \resizebox{0.8\textwidth}{!}{%
\begin{tikzpicture}
  \begin{axis}
  [
    xlabel={Number of branches $C$},
    ylabel={${\rm RVMS}$},
    xtick={2,3,4,5},
    yticklabel=
{\pgfmathparse{\tick*100}\pgfmathprintnumber{\pgfmathresult}\%},
cycle list name=black white,
    legend pos=north west,
    legend cell align={left}
]
    \addplot coordinates {
(2,	0.053879724)
(3, 0.084522619)
(4, 0.101842469)
(5, 0.106378299)
    };
    \addplot coordinates {
(2,	0.023220165)
(3, 0.040960461)
(4, 0.047273132)
(5, 0.055521724)
    };\pgfplotsset{cycle list shift=2}
        \addplot coordinates {
(2,	0.095943103)
(3, 0.1138455)
(4, 0.118128264)
(5, 0.118755166)
    };\pgfplotsset{cycle list shift=5}
    \addplot coordinates {
(2,	0.072144243)
(3, 0.088004343)
(4, 0.083690602)
(5, 0.1014)
    };
  \end{axis}
\end{tikzpicture}%
}
 \caption{different numbers of branches $C$}
\end{subfigure}
\begin{subfigure}{0.45\textwidth}
 \resizebox{0.8\textwidth}{!}{%
 \pgfplotsset{scaled y ticks=false}
\begin{tikzpicture}
  \begin{axis}
  [
    xlabel={Number of stages $T$},
    ylabel={${\rm RVMS}$},
    xtick={3,4,5,6},
    yticklabel= {\pgfmathparse{\tick*100}\pgfmathprintnumber{\pgfmathresult}\%},
    cycle list name=black white,
    legend pos=north west,
    ]
    \addplot coordinates {
(3,	0.0538797)
(4, 0.097359633)
(5, 0.125050732)
(6, 0.157416325)
    };
    \addplot coordinates {
(3,	0.023220165)
(4, 0.031441941)
(5, 0.027240595)
(6, 0.032914946)
    };\pgfplotsset{cycle list shift=2}
        \addplot coordinates {
(3,	0.095943103)
(4, 0.134838054)
(5, 0.161966104)
(6, 0.161726302)
    };\pgfplotsset{cycle list shift=5}
    \addplot coordinates {
(3,	0.072144243)
(4, 0.105720491)
(5, 0.13750083)
(6, 0.157555622)
    };
  \end{axis}
\end{tikzpicture}%
}
\caption{different numbers of stages $T$}
\end{subfigure}
\begin{subfigure}{0.45\textwidth}
 \resizebox{0.8\textwidth}{!}{%
\begin{tikzpicture}
    \begin{axis}
  [
    xlabel={Risk attitude $\lambda$},
    ylabel={${\rm RVMS}$},
    xtick={0,0.2,0.4,0.6,0.8,1},
    yticklabel=
{\pgfmathparse{\tick*100}\pgfmathprintnumber{\pgfmathresult}\%},
    cycle list name=black white,
    ]
    \addplot coordinates {
(0,	0.0904126)
(0.2, 0.074066481)
(0.4, 0.063876339)
(0.6, 0.051600539)
(0.8, 0.052515345)
(1, 0.04020922)
    };
    \addplot coordinates {
(0, 0.055935392)
(0.2, 0.040440322)
(0.4, 0.030788046)
(0.6, 0.021599204)
(0.8, 0.012176273)
(1, 0.000832135)
    };\pgfplotsset{cycle list shift=2}
        \addplot coordinates {
(0,	0.162296678)
(0.2, 0.132706602)
(0.4, 0.104816938)
(0.6, 0.085998544)
(0.8, 0.059408474)
(1, 0.048493188)
    };\pgfplotsset{cycle list shift=5}
    \addplot coordinates {
(0,	0.001733403)
(0.2, 0.031806897)
(0.4, 0.058487404)
(0.6, 0.084012081)
(0.8, 0.102432776)
(1, 0.127836551)
    };
  \end{axis}
\end{tikzpicture}%
}
\caption{different risk attitudes $\lambda$}
\end{subfigure}
\begin{subfigure}{0.45\textwidth}
 \resizebox{0.8\textwidth}{!}{%
 \pgfplotsset{scaled y ticks=false}
\begin{tikzpicture}
    \begin{axis}
  [
    xlabel={Demand standard deviation $\sigma$},
    ylabel={${\rm RVMS}$},
    xtick={0.2,0.4,0.6,0.8},
    yticklabel={\pgfmathparse{\tick*100}\pgfmathprintnumber{\pgfmathresult}\%},
    cycle list name=black white,
    legend pos=north west,
    ]
    \addplot coordinates {
(0.2,	0.013141937)
(0.4, 0.026496662)
(0.6, 0.044739945)
(0.8, 0.053879724)
    };
    \addplot coordinates {
(0.2,	0.009441415)
(0.4, 0.013672261)
(0.6, 0.021568807)
(0.8, 0.023220165)
    };\pgfplotsset{cycle list shift=2}
        \addplot coordinates {
(0.2,	0.03812985)
(0.4, 0.054939824)
(0.6, 0.075508382)
(0.8, 0.095943103)
    };\pgfplotsset{cycle list shift=5}
    \addplot coordinates {
(0.2,	0.030046305)
(0.4, 0.045977127)
(0.6, 0.061281109)
(0.8, 0.072144243)
    };
    \legend{SD-${\rm RVMS_R}$, SI-${\rm RVMS_R}$,SD0-${\rm RVMS_R}$, SD-${\rm RVMS_P}$}
  \end{axis}
\end{tikzpicture}%
}
\caption{different standard deviations $\sigma$}
\end{subfigure}
    \caption{Statistics of ${\rm RVMS}$ over 100 instances with different numbers of branches $C$, stages $T$, risk attitudes $\lambda$ and demand standard deviations $\sigma$.}
    \label{fig:sensitivity_RVMS}
\end{figure}
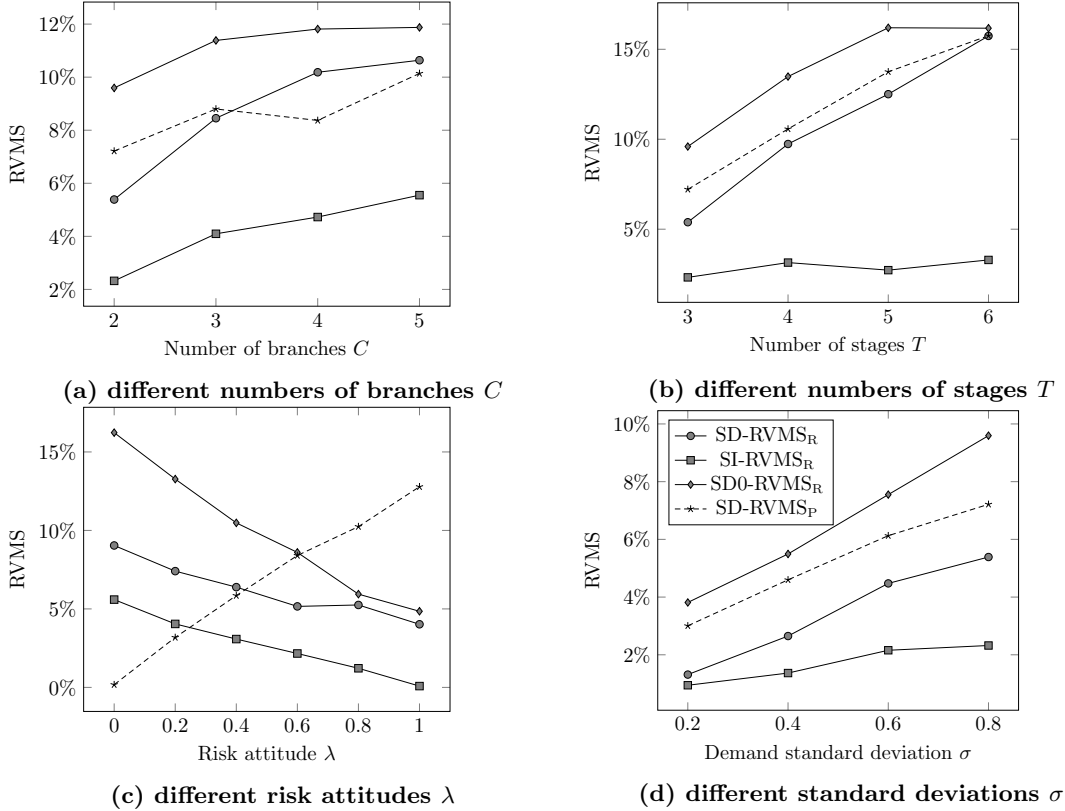

Using {the three types of scenario trees defined in Section \ref{sec:setup}}, we first vary the number of branches $C$ from 2 to 5, the number of stages $T$ from 3 to 6, the risk attitude $\lambda$ from 0 to 1, and the standard deviation $\sigma$ from 0.2 to 0.8 to see how ${\rm RVMS_R}$ {and ${\rm RVMS_P}$} change with respect to different parameter settings. The corresponding results are presented in Figure \ref{fig:sensitivity_RVMS}, where we plot the mean of ${\rm RVMS_R}$ {and ${\rm RVMS_P}$} over 100 independently generated instances.

From the figure, when we increase the demand variability such as the number of branches, {number of stages} and the standard deviation, ${\rm RVMS_R}$ {and ${\rm RVMS_P}$ with stagewise dependent scenario trees} increase dramatically (see Figures \ref{fig:sensitivity_RVMS}(a), (b) and (d)). 
A higher ${\rm RVMS}$ indicates that the multistage model will gain more benefits over the two-stage counterpart. 
{Comparing different types of scenario trees, SD scenario trees always gain much higher ${\rm RVMS}$ than the stagewise independent ones, and adding more stages makes no significant changes in ${\rm RVMS_R}$ when using the SI scenario trees, as can be seen in Figure \ref{fig:sensitivity_RVMS}(b). This is because in SI scenario trees, the number of realizations in each stage does not depend on the number of stages, and having a deeper scenario tree would not necessarily increase the demand variability. }
Notably, from Figure \ref{fig:sensitivity_RVMS}(c), ${\rm RVMS_R}$ decreases approximately linearly with respect to the risk attitude parameter $\lambda$, {while ${\rm RVMS_P}$ increases approximately linearly}, which is due to the linear dependence of ${\rm VMS_R^{LB}}$ and ${\rm VMS_P^{LB}}$ on $\lambda$.

{
\subsubsection{Risk Profile of the Solutions via Out-of-sample Test}\label{sec:risk-profile}
We aim to compare the risk profile of the in-sample solutions generated by two-stage and multistage models via out-of-sample test. For risk-averse two-stage models without prioritization, this raises the first issue -- the facility locations $\boldsymbol{x}_1,\ldots,\boldsymbol{x}_T$ are decided in the first stage based on the in-sample scenarios $\{\boldsymbol{d}_t(\omega)\}_{\omega\in{\Omega}}$; however, under a different set of out-of-sample scenarios $\{\boldsymbol{d}_t(\omega)\}_{\omega\in\tilde{\Omega}}$, the implemented decisions  $\boldsymbol{x}_1,\ldots,\boldsymbol{x}_T$ become infeasible if there exists $\omega\in\tilde{\Omega}$ such that $\sum_{j=1}^Nd_{tj}(\omega)>\sum_{i=1}^Nh_i\sum_{\tau=1}^tx_{\tau i}$. Because of this, we only conduct out-of-sample test on risk-averse two-stage and multistage models with prioritization. In this case, the implemented decisions are the priority lists, which are always feasible under different demand and budget realizations as long as the budget is set large enough to cover all the demand in each scenario. As the multistage models generate a solution on each node of the scenario tree rather than generate a policy, when the scenarios change from in-sample to out-of-sample, we do not have corresponding solutions on hand. To conduct out-of-sample test, we follow a rolling horizon approach and only implement the first-stage decisions, which we shall describe as follows.

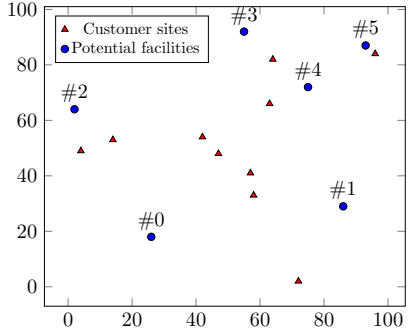
\begin{figure}[ht!]
\centering
\begin{tikzpicture}[scale=0.7]
\begin{axis}[legend style={nodes={scale=0.8, transform shape}}, 
legend pos=north west]
    \addplot[
        scatter/classes={        a={mark=triangle*, fill=red},
        b={mark=*,fill=blue}},
        scatter, mark=*, only marks, 
        scatter src=explicit symbolic,
        nodes near coords*={\Label},
        visualization depends on={value \thisrow{label} \as \Label} 
    ] table [meta=class] {
        x y class label
        42 54 a {}
        96 84 a {}
        58 33 a {}
        14 53 a {}
        72 2 a {}
        4 49 a {}
        47 48 a {}
        64 82 a {}
        57 41 a {}
        63 66 a {}
        26 18 b \#0
        86 29 b \#1
        2 64 b \#2
        55 92 b \#3
        75 72 b \#4
        93 87 b \#5
    };
    \legend{Customer sites, Potential facilities}
\end{axis}
\end{tikzpicture}
\caption{Locations of customer sites and potential facilities on a 100$\times$100 grid} \label{fig:locations}
\end{figure}

Using $T=3$ and the default setting, we plot the locations of $M=6$ potential facility sites and $N=10$ customer sites of a toy example in Figure \ref{fig:locations}. We first construct a scenario tree with a three-stage forecast, based on which we generate solutions $\boldsymbol{s}_0^{TS}$ and $\boldsymbol{s}_0^{MS}$. Now we move one stage forward and observe the uncertainty $\boldsymbol{d}_1(\omega_1),{F}_1(\omega_1)$, which may be different from all the in-sample scenarios that we estimate. For the two-stage model, we always implement $\boldsymbol{s}_0^{TS}$; for the multistage model, we implement $\boldsymbol{s}_0^{MS}$ and update the scenario tree with a two-stage forecast, based on which we optimize a new priority list $\boldsymbol{s}_1^{MS}$. We continue this process until we reach the last stage, and this will give us an operational cost of $\sum_{t=1}^T\boldsymbol{c}_t^{\mathsf T}\boldsymbol{y}_t(\omega_1)$ under the specific out-of-sample scenario path $\omega_1$. We repeat this process for 100 independently generated out-of-sample scenario paths $\{\omega_i\}_{i=1}^{100}$, and record the 95\% percentile, 75\% percentile, and the mean of the operational cost produced by two-stage and multistage solutions ($\boldsymbol{s}_0$) across out-of-sample scenarios in Table \ref{tab:out-of-sample}, where we vary the risk parameter $\lambda$ from 0 to 1 and Columns ``Time'' present the average time for computing in-sample solutions in seconds.

From Table \ref{tab:out-of-sample}, although multistage models require slightly more time for computing solutions, they obtain lower 95\% percentiles and means of the operational cost compared to two-stage models, where the gaps of the 95\% percentile of multistage and two-stage costs are amplified as the risk parameter $\lambda$ increases to 1. Moreover, looking at the optimal priority list $\boldsymbol{s}_0$ and Figure \ref{fig:locations}, the multistage solutions are more consistent across different risk parameters, with facilities \#4, \#2, \#0 always in the top 3. On the other hand, two-stage models change the priority list more often with varying risk parameters, where facilities \#4, \#2, \#1 are always in the top 3 with varying orders.
\begin{table}[ht!]
  \centering
  \caption{Risk profile of two-stage and multistage solutions with prioritization under different risk parameters via out-of-sample test}
  \resizebox{\textwidth}{!}{
    \begin{tabular}{l|rrrcr|rrrcr}
    \hline
          & \multicolumn{5}{c|}{Two-stage}         & \multicolumn{5}{c}{Multistage} \\
    $\lambda$ & 95\%  & 75\%  & Mean  & $\boldsymbol{s}_0$  & Time  & 95\%  & 75\%  & Mean  & $\boldsymbol{s}_0$  & Time \\
    \hline
    0     & \$115,759 & \$101,891 & \$67,064 & [4,2,5,1,3,0] & 0.43  & \$115,713 & \$101,891 & \$67,047 & [4,2,0,3,1,5] & 0.56 \\
    0.2   & \$119,906 & \$47,350 & \$49,456 & [2,4,1,5,3,0] & 0.44  & \$119,762 & \$44,280 & \$47,872 & [4,2,0,5,1,3] & 0.63 \\
    0.4   & \$132,957 & \$110,277 & \$99,602 & [2,1,4,5,3,0] & 0.49  & \$132,957 & \$109,880 & \$98,861 & [2,4,0,1,3,5] & 0.60 \\
    0.6   & \$125,202 & \$110,531 & \$87,660 & [4,1,2,5,3,0] & 0.40  & \$125,202 & \$109,548 & \$87,104 & [4,2,1,0,3,5] & 0.63 \\
    0.8   & \$135,689 & \$111,232 & \$87,567 & [0,4,2,5,1,3] & 0.45  & \$129,309 & \$108,455 & \$84,331 & [4,2,0,1,3,5] & 0.64 \\
    1     & \$133,256 & \$83,849 & \$62,966 & [4,2,1,0,3,5] & 0.39  & \$132,889 & \$85,136 & \$62,956 & [4,0,1,3,2,5] & 0.69 \\
    \hline
    \end{tabular}%
    }
  \label{tab:out-of-sample}%
\end{table}%
}

\subsubsection{Tightness of the Analytical Bound on the Synthetic Data Set}
\label{sec:tight}
Using the default setting and the three scenario trees defined in Section \ref{sec:setup}, we present the mean of ${\rm RGAP_R}$ {and ${\rm RGAP_P}$} over 100 independently generated instances in Table \ref{tab:statis-RGAP}, where in the last three columns, we display the percentage of instances that RGAP does not exceed given thresholds. {From the table, SD0-${\rm RGAP_R}$ obtains the lowest ${\rm RGAP_R}$, where in 98\% of the instances ${\rm RGAP_R}$ is within 10\% and about 81\% of the instances attain nearly tight lower bounds (${\rm RGAP_R}\le 10^{-5}$). A lower ${\rm RGAP_R}$ means that the analytical bound can recover the true $\rm{VMS_R}$ better.} This is because by adding demand 0 to the scenario tree, the utilization rate of each facility will vary significantly across different scenarios and thus the analytical bound $\rm{VMS_R^{LB}}$ will be large enough to recover the true $\rm{VMS_R}$. {Moreover, ${\rm RGAP_P}$ is much lower than ${\rm RGAP_R}$ with a mean of 3.66\%, while SD-${\rm RGAP_R}$ attains a mean of 30.32\%.}

\begin{table}[ht!]
  \centering
  \caption{Statistics of ${\rm RGAP_R}$ and ${\rm RGAP_P}$ over 100 instances with different types of scenario trees}
    \begin{tabular}{lrrrr}
    \hline
    RGAP  & Mean  & $< 10^{-5}$ & $<10\%$ & $ < 50\%$ \\
    \hline
    SD-${\rm RGAP_R}$ & 30.32\% & 58/100 & 61/100 & 70/100 \\
    SI-${\rm RGAP_R}$   & 38.99\% & 51/100 & 61/100 & 61/100 \\
    SD0-${\rm RGAP_R}$ & 1.30\% & 81/100 & 98/100 & 100/100 \\
    SD-${\rm RGAP_P}$ & 3.66\% & 0/100 & 45/100 & 100/100 \\
    \hline
    \end{tabular}%
  \label{tab:statis-RGAP}%
\end{table}%

\subsubsection{Performance of the Approximation Algorithm on the Synthetic Data Set}
\label{sec:aa}
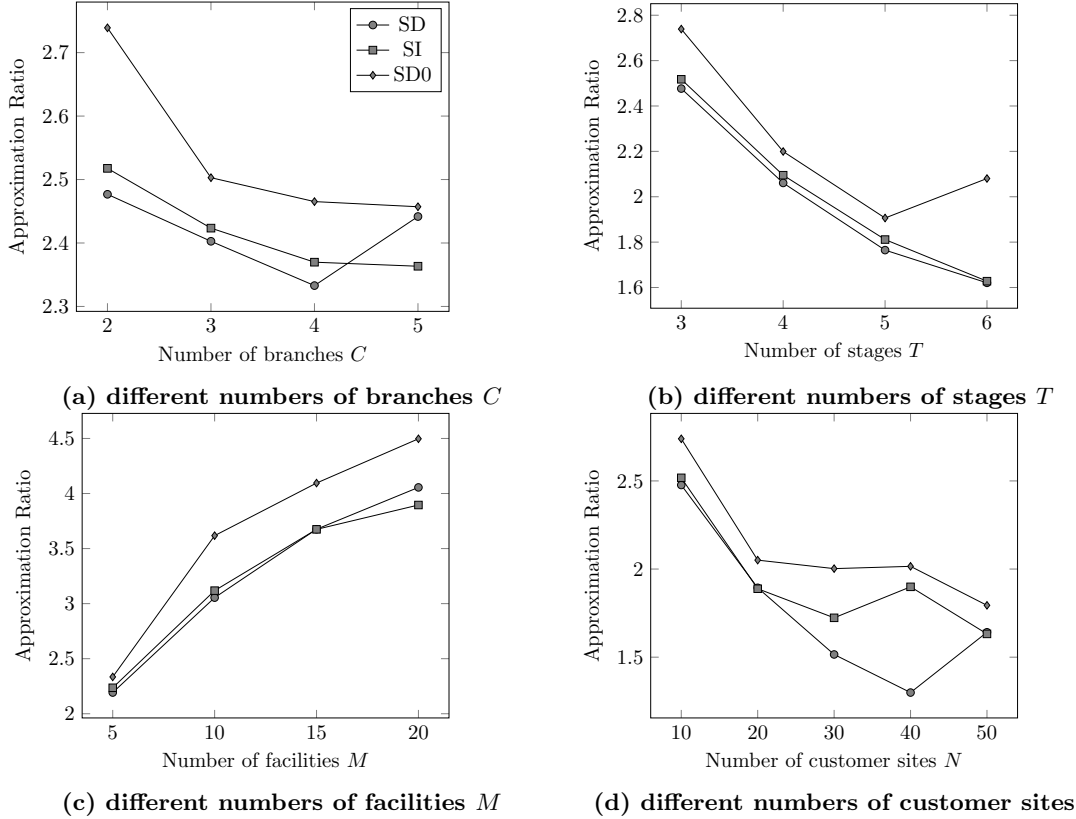
\begin{figure}[ht!]
    \centering
    \begin{subfigure}{0.45\textwidth}
 \resizebox{0.8\textwidth}{!}{%
\begin{tikzpicture}
  \begin{axis}
  [
    xlabel={Number of branches $C$},
    ylabel={Approximation Ratio},
    xtick={2,3,4,5},
cycle list name=black white,
]
    \addplot coordinates {
(2,	2.476730857)
(3, 2.402681307)
(4, 2.332868125)
(5, 2.441697833)
    };
    \addplot coordinates {
(2,	2.517607232)
(3, 2.423413357)
(4, 2.369746646)
(5, 2.3633659)
    };\pgfplotsset{cycle list shift=2}
        \addplot coordinates {
(2,	2.739056031)
(3, 2.50292337)
(4, 2.465232935)
(5, 2.457120184)
    };
    \legend{SD,SI,SD0}
  \end{axis}
\end{tikzpicture}%
}
\caption{different numbers of branches $C$}
\end{subfigure}
\begin{subfigure}{0.45\textwidth}
 \resizebox{0.8\textwidth}{!}{%
 \pgfplotsset{scaled y ticks=false}
\begin{tikzpicture}
  \begin{axis}
  [
    xlabel={Number of stages $T$},
    ylabel={Approximation Ratio},
    xtick={3,4,5,6},
    cycle list name=black white,
    legend pos=north east,
    legend cell align={left}
    ]
    \addplot coordinates {
(3,	2.476730857)
(4, 2.060774757)
(5, 1.764778566)
(6, 1.620526893)
    };
    \addplot coordinates {
(3,	2.517607232)
(4, 2.09495391)
(5, 1.811327221)
(6, 1.627872845)
    };\pgfplotsset{cycle list shift=2}
        \addplot coordinates {
(3,	2.739056031)
(4, 2.199404691)
(5, 1.906228829)
(6, 2.08032125)
    };
  \end{axis}
\end{tikzpicture}%
}
\caption{different numbers of stages $T$}
\end{subfigure}
\begin{subfigure}{0.45\textwidth}
 \resizebox{0.8\textwidth}{!}{%
 \pgfplotsset{scaled y ticks=false}
\begin{tikzpicture}
    \begin{axis}
  [
    xlabel={Number of facilities $M$},
    ylabel={Approximation Ratio},
    xtick={5,10,15,20},
    cycle list name=black white,
    ]
    \addplot coordinates {
(5,	2.192350574)
(10, 3.054235222)
(15, 3.675647517)
(20, 4.055468873)
    };
    \addplot coordinates {
(5,	2.235842619)
(10, 3.117941655)
(15, 3.674414355)
(20, 3.895617572)
    };\pgfplotsset{cycle list shift=2}
        \addplot coordinates {
(5,2.334332363)
(10, 3.617417265)
(15, 4.094921732)
(20, 4.497066408)
    };
  \end{axis}
\end{tikzpicture}%
}
\caption{different numbers of facilities $M$}
\end{subfigure}
\begin{subfigure}{0.45\textwidth}
 \resizebox{0.8\textwidth}{!}{%
 \pgfplotsset{scaled y ticks=false}
\begin{tikzpicture}
    \begin{axis}
  [
    xlabel={Number of customer sites $N$},
    ylabel={Approximation Ratio},
    xtick={10,20,30,40,50},
    cycle list name=black white,
    ]
    \addplot coordinates {
(10,2.476730857)
(20, 1.89368279)
(30, 1.515400313)
(40, 1.299284718)
(50, 1.641334832)
    };
    \addplot coordinates {
(10,	2.517607232)
(20, 1.888373873)
(30, 1.7233382)
(40, 1.899176664)
(50, 1.632115688)
    };\pgfplotsset{cycle list shift=2}
        \addplot coordinates {
(10,	2.739056031)
(20, 2.050585428)
(30, 2.00246752)
(40, 2.015358044)
(50, 1.794068077)
    };
  \end{axis}
\end{tikzpicture}%
}
\caption{different numbers of customer sites $N$}
\end{subfigure}
    \caption{Statistics of the approximation ratios over 100 instances with different numbers of branches $C$, stages $T$, risk attitudes $\lambda$, demand standard deviations $\sigma$, facilities $M$ and customer sites $N$.}
    \label{fig:sensitivity_ratio}
\end{figure}

Next, using {the three types of scenario trees defined in Section \ref{sec:setup}}, we vary the number of branches $C$ from 2 to 5, number of stages $T$ from 3 to 6, number of facilities $M$ from 5 to 20 and number of customer sites $N$ from 10 to 50 to see how the empirical approximation ratio changes with respect to different parameter settings. The results are presented in Figure \ref{fig:sensitivity_ratio}, where we plot the mean of empirical approximation ratios (i.e., $\frac{z_R^{MS}(\boldsymbol{x}_n^{H}, \eta_n^{H}, \boldsymbol{y}_n^{H}, u_n^{H})}{z_R^{MS}(\boldsymbol{x}_n^{*}, \eta_n^{*}, \boldsymbol{y}_n^{*}, u_n^{*})}$) over 100 independently generated instances.
From the figure, the approximation ratios decrease gradually when we increase the number of branches $C$ and the number of stages $T$ (see Figures \ref{fig:sensitivity_ratio}(a) and (b)). 
Moreover, from Figures \ref{fig:sensitivity_ratio}(c) and (d), the approximation ratios are clearly positively related to the number of facilities $M$ and negatively impacted by the number of customer sites $N$, which is because increasing $N$ would increase the {total demand and} first-stage facilities needed $M_{\rm min}$ and our approximation ratio has an upper bound $1+ {\frac{M\sum_{t=1}^Tf_{t,\rm max}}{M_{\rm min}\sum_{t=1}^Tf_{t,\rm min}+\sum_{t=1}^Tc_{t,\rm min}\min_{n\in\mathcal{T}_t}\{\sum_{j=1}^Nd_{n,j}\}}}$.

\subsubsection{Computational Time Comparison}
\label{sec:time}
\begin{figure}[ht!]
    \centering
    \begin{subfigure}{0.45\textwidth}
 \resizebox{0.8\textwidth}{!}{%
\begin{tikzpicture}
  \begin{axis}
  [
    xlabel={Number of branches $C$},
    ylabel={Computational Time (sec.)},
    xtick={2,3,4,5},
cycle list name=black white,
legend pos=north west,
legend cell align={left}
]
    \addplot coordinates {
(2,	0.403188362)
(3, 1.318968632)
(4, 5.881092119)
(5, 27.03570475)
    };\pgfplotsset{cycle list shift=5}
    \addplot coordinates {
    (2,	0.218538175)
(3, 0.408468866)
(4, 0.689959598)
(5, 1.112607012)
    };\pgfplotsset{cycle list shift=11}
        \addplot coordinates {
(2,	0.366190119)
(3, 0.850977123)
(4, 1.748567886)
(5, 3.020830956)
    };
    \legend{$z_R^{MS}$, $z_R^{TS}$, $z_R^{MS}$ via AA}
  \end{axis}
\end{tikzpicture}%
}
\caption{different numbers of branches $C$}
\end{subfigure}
\begin{subfigure}{0.45\textwidth}
 \resizebox{0.8\textwidth}{!}{%
 \pgfplotsset{scaled y ticks=false}
\begin{tikzpicture}
  \begin{axis}
  [
    xlabel={Number of stages $T$},
    ylabel={Computational Time (sec.)},
    xtick={3,4,5,6},
    cycle list name=black white,
    legend pos=north east,
    ]
    \addplot coordinates {
(3,	0.403188362)
(4, 2.158197207)
(5, 16.06150371)
(6, 88.35860449)
    };\pgfplotsset{cycle list shift=5}
    \addplot coordinates {
    (3,	0.218538175)
(4, 0.573158755)
(5, 1.473669429)
(6, 2.879646306)
    };\pgfplotsset{cycle list shift=11}
        \addplot coordinates {
        (3,	0.366190119)
(4, 1.044921248)
(5, 3.27039144)
(6, 9.64644186)
    };
  \end{axis}
\end{tikzpicture}%
}
\caption{different numbers of stages $T$}
\end{subfigure}
\begin{subfigure}{0.45\textwidth}
 \resizebox{0.8\textwidth}{!}{%
\begin{tikzpicture}
    \begin{axis}
  [
    xlabel={Number of facilities $M$},
    ylabel={Computational Time (sec.)},
    xtick={5,10,15,20},
    cycle list name=black white,
    ]
    \addplot coordinates {
(5,	0.28951076)
(10, 1.249616432)
(15, 4.261511939)
(20, 16.52168616)
    };\pgfplotsset{cycle list shift=5}
        \addplot coordinates {
(5,	0.173255529)
(10, 0.551687415)
(15, 0.725851932)
(20, 1.141980596)
    };\pgfplotsset{cycle list shift=11}
    \addplot coordinates {
    (5,	0.31246347)
(10, 0.620138462)
(15, 0.962169313)
(20, 1.30137593)
    };
  \end{axis}
\end{tikzpicture}%
}
\caption{different numbers of facilities $M$}
\end{subfigure}
\begin{subfigure}{0.45\textwidth}
 \resizebox{0.8\textwidth}{!}{%
\begin{tikzpicture}
    \begin{axis}
  [
    xlabel={Number of customer sites $N$},
    ylabel={Computational Time (sec.)},
    xtick={10,20,30,40,50},
    cycle list name=black white,
    ]
    \addplot coordinates {
(10,	0.403188362)
(20, 0.886932755)
(30, 0.795222378)
(40, 0.760848901)
(50, 1.262223878)
    };\pgfplotsset{cycle list shift=5}
        \addplot coordinates {
        (10,	0.218538175)
(20, 0.444063516)
(30, 0.483927853)
(40, 0.592819052)
(50, 0.941403413)
    };\pgfplotsset{cycle list shift=11}
    \addplot coordinates {
    (10,	0.366190119)
(20, 0.553561032)
(30, 0.781796227)
(40, 1.042295749)
(50, 1.353335094)
    };
  \end{axis}
\end{tikzpicture}%
}
\caption{different numbers of customer sites $N$}
\end{subfigure}
    \caption{Computational time comparison of solving the multistage model using Gurobi and the approximation algorithm, solving the two-stage model using Gurobi  with different numbers of branches $C$, stages $T$, facilities $M$ and customer sites $N$.}
    \label{fig:time}
\end{figure}
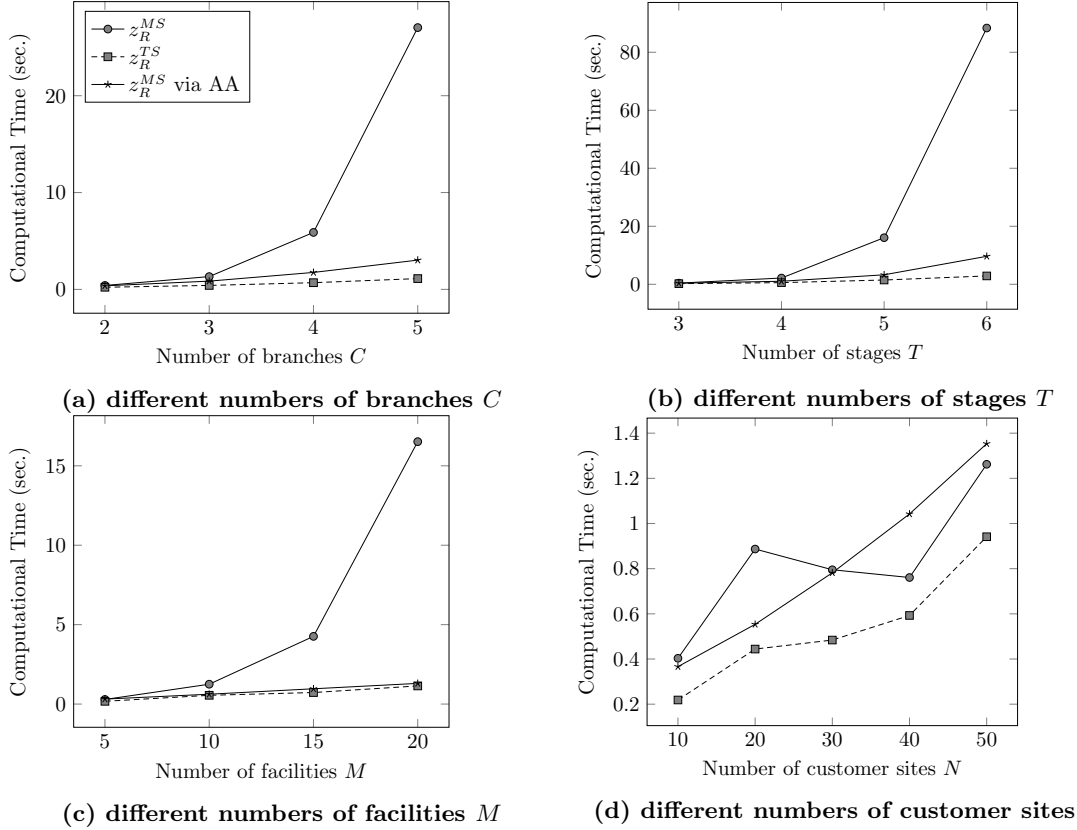

We end this section by showing the computational time for various risk-averse models {with or without prioritization. We first compare the time of two-stage and multistage models without prioritization} solved to optimality, and multistage models solved by the approximation algorithm (AA) in Figure \ref{fig:time}. We consider {SD scenario trees} and fix the risk parameters $\lambda=0.5,\ \alpha=0.95$ and the standard deviation $\sigma=0.8$ while varying the number of branches $C$ from 2 to 5, number of stages $T$ from 3 to 6, number of facilities $M$ from 5 to 20 and number of customer sites $N$ from 10 to 50. 
We record the average time over 100 independently generated instances. Note that we solve Step \ref{alg:step5} for all node $n\in\mathcal{T}$ as a whole in Algorithm \ref{alg:approx-multistage}, and one may further speed up the approximation algorithms by utilizing parallel computing techniques. From the figure, we observe that the computational time of the approximation algorithm grows approximately linearly with respect to the problem size, while the one for solving multistage models to optimality grows exponentially with respect to the number of branches $C$, stages $T$ and facilities $M$. Both models scale well in terms of the number of customer sites $N$, where the approximation algorithm requires more time but can still solve the multistage models within 2 seconds.
Moreover, two-stage models are less computational expensive, and therefore can be preferred over multistage models when their objective gaps are relatively small.

{
Lastly, we record the computational time of the two-stage and multistage models \eqref{model:two-prior} and \eqref{model:multi-prior} with and without employing the prioritization cuts in Table 2, where we vary the number of facilities ($M$) and number of customer sites ($N$), respectively. The CPU time limit is set to 3600 seconds, and the optimality gaps after one hour of computation are marked in the brackets. From the table, as the instance size grows, the benefit of using prioritization cuts is amplified, where we can solve some instances with the prioritization cuts to optimality in one hour but could not solve in this time without using the cuts. Moreover, for larger instances, prioritization cuts can reduce the optimality gaps within the one-hour time frame. 

\begin{table}[ht!]
  \centering
  \caption{Computational time comparison of solving the risk-averse models \eqref{model:two-prior} and \eqref{model:multi-prior} with or without prioritization cuts with different number of facilities $M$ and number of customer sites $N$}
    \begin{tabular}{lrrrr}
    \hline
    $(M, N)$ & $z_P^{TS}$ without cuts     & $z_P^{TS}$ with cuts & $z_P^{MS}$ without cuts     & $z_P^{MS}$ with cuts \\
    \hline
    (50, 10) & 90.59 & 67.10 & 77.73 & 60.59 \\
    (50, 20) & 333.83 & 170.38 & 278.31 & 150.47 \\
    (50, 30) & 413.66 & 214.70 & 367.79 & 208.54 \\
    (50, 40) & 3600 (0.46\%) & 393.88 & 2839.26 & 261.98 \\
    (50, 50) & 3600 (1.47\%) & 1632.88 & 3546.23 & 879.63 \\
    (60, 50) & 3600 (2.81\%) & 3600 (0.65\%) & 1674.38 & 520.57 \\
    (70, 50) & 3600 (1.92\%) & 3600 (0.14\%) & 3600 (0.99\%) & 709.49 \\
    (80, 50) & 3600 (1.90\%) & 3600 (1.23\%) & 3600 (2.02\%) & 3600 (0.54\%) \\
    \hline
    \end{tabular}%
  \label{tab:addlabel}%
\end{table}%

}

\subsection{Case Study on a Real-World Network}
\label{sec:real}
\subsubsection{Experimental Design and Setup}
We consider the 49-node and 88-node data sets described in \citet{daskin2011network} consisting of the capitals
of the continental United States plus Washington, DC, which can be used as candidate facilities ({$M=49$}), and 88 major cities in United States, which can represent the customer sites ({$N=88$}), respectively. 
The capacities of 48 candidate facilities are drawn uniformly between $10^5$ and $10^6$.
Flow costs are set equal to the great-circle distance times the travel cost per mile per unit of demand, i.e., $c_{ij}=\textrm{dist}(i,j) * 0.00001$.
Since these benchmarks are designed for deterministic facility location problems, they do not provide random demand data for each customer site in each year. We generate demand data as follows. We first collect the population data for each city in the United States, and multiply them by 2\% times 12 months, assuming that 2\% of the population will order once per month, which gives nominal demand $\tilde{d}$ in each customer site in the beginning year of the planning horizon. 
We consider four demand patterns (described in Column ``Pattern'' in Table \ref{table:pattern}), all of which follow truncated Normal distributions $\mathcal N(\cdot, *)$ in each stage $t$, where $\cdot$ indicates the mean value and $*$ represents the standard deviation. In Patterns III and IV, the nominal demand $\tilde{d}$ is increased with a rate of $2\%$ for each subsequent year by assuming that the population increase rate is roughly $2\%$; in Patterns II and IV, the standard deviation is increased with the same rate.
\begin{table}[ht!]
	\centering
	\caption{Demand patterns}
	\begin{tabular}{lr}
		\hline
Pattern & Distribution at stage $t$ \\
\hline
I.\ Constant mean, constant standard deviation & $\mathcal N(\tilde{d},\tilde{d}\cdot \sigma)$\\
 II.\ Constant mean, increasing standard deviation & $\mathcal N(\tilde{d},\tilde{d}\cdot (\sigma+0.2(t-1)))$\\
 III.\ Increasing mean, constant standard deviation & $\mathcal N(\tilde{d}\cdot (1+0.2(t-1)),\tilde{d}\cdot \sigma)$\\
 IV.\ Increasing mean, increasing standard deviation & $\mathcal N(\tilde{d}\cdot (1+0.2(t-1)),\tilde{d}\cdot (\sigma+0.2(t-1)))$\\
           \hline
\end{tabular}%
\label{table:pattern}%
\end{table}%

\subsubsection{Results without Prioritization}
With the baseline setting and SD scenario trees, we present the optimal solutions and cost breakdown of two-stage and multistage models under different demand patterns in Table \ref{tab:real-sol}, where Columns ``$|\boldsymbol{x}_1|$'', ``$|\boldsymbol{x}_2|$'' and ``$|\boldsymbol{x}_3|$'' display the number of distinct facilities rented in each stage across all scenarios. Columns ``Renting (\$)'' and  ``Trans. (\$)'' show the renting and operational cost without considering the risk parameters, and Column ``Obj. (\$)'' presents the overall risk-averse optimal objective values, where we mark the lowest ones among the three models in bold.  Note that in multistage models, different stages may rent some facilities in common; this is because multistage models have the flexibility to rent facilities under certain (not all) scenarios in each stage and these common facilities are rented in different scenarios. As a result, although multistage models have more variability in which facilities to rent, they obtain the lowest renting cost as well as the lowest overall objective values. In terms of operational cost, the multistage models solved by approximation algorithms {sometimes} achieve the minimum among the three. The last column displays the ${\rm RVMS_R}$ when comparing the optimal objective values between two-stage and multistage models and approximation ratios (AR) when comparing the ones of multistage models solved to optimality and solved by approximation algorithms. From this column, {demand patterns with increasing standard deviation obtain higher ${\rm RVMS_R}$ compared to the ones with constant standard deviation, which agrees with our findings in Section \ref{sec:sensi} that higher demand variation leads to higher ${\rm RVMS_R}$}.

In Patterns I and II, most facilities are rented in the first stage. In Patterns III and IV where the demand mean is increased with a rate of 2\%, more facilities are built in later stages. Notably, although the two-stage models are less computational expensive than the multistage counterparts, their objective gaps can be as high as {12.37\%}, which can be further increased with more branches and standard deviations. Moreover, the approximation algorithms can solve the multistage models in an extremely quick fashion with an approximation ratio of at most {1.07}.  

\begin{table}[ht!]
  \centering
  \caption{Optimal solutions and cost breakdown of two-stage and multistage models without prioritization under different demand patterns}
  \resizebox{\textwidth}{!}{
    \begin{tabular}{clrrrrrrrr}
    \hline
\multicolumn{1}{l}{Pattern} & Model & $|\boldsymbol{x}_1|$  & $|\boldsymbol{x}_2|$  & $|\boldsymbol{x}_3|$  & Renting (\$) & Trans. (\$) & Obj. (\$)  & Time (sec.) & ${\rm RVMS_R}$/AR \\
\hline
    \multirow{3}[0]{*}{I} & $z_R^{TS}$ via OPT    & 12    & 7     & 0     & \$3,238,500 & \$133,930 & \$3,385,469 & 23.65 & 5.12\% \\
          & $z_R^{MS}$ via OPT    & 12    & 7     & 1     & \textbf{\$2,992,875} & \$146,342 & \textbf{\$3,220,542} & 41.61 & N.A. \\
          & $z_R^{MS}$ via AA & 14    & 5     & 0     & \$3,172,800 & \textbf{\$132,835} & \$3,366,800 & \textbf{17.99} & 1.05 \\
          \hline
    \multirow{3}[0]{*}{II} & $z_R^{TS}$ via OPT    & 13    & 9     & 8     & \$4,404,300 & \textbf{\$162,995} & \$4,589,469 & \textbf{8.76}  & 12.37\% \\
          & $z_R^{MS}$ via OPT    & 13    & 9     & 11    & \textbf{\$3,751,250} & \$185,645 & \textbf{\$4,084,279} & 18.83 & N.A. \\
          & $z_R^{MS}$ via AA & 15    & 8     & 10    & \$4,059,625 & \$170,131 & \$4,366,258 & 18.68 & 1.07 \\
          \hline
    \multirow{3}[0]{*}{III} & $z_R^{TS}$ via OPT    & 13    & 3     & 12    & \$4,001,000 & \$143,189 & \$4,159,788 & \textbf{13.25} & 5.19\% \\
          & $z_R^{MS}$ via OPT    & 13    & 4     & 13    & \textbf{\$3,689,825} & \$151,421 & \textbf{\$3,954,520} & 152.18 & N.A. \\
          & $z_R^{MS}$ via AA & 14    & 4     & 11    & \$3,941,600 & \textbf{\$134,665} & \$4,193,486 & 18.06 & 1.06 \\
          \hline
    \multirow{3}[0]{*}{IV} & $z_R^{TS}$ via OPT    & 13    & 8     & 15    & \$4,920,600 & \textbf{\$136,647} & \$5,066,159 & 23.73 & 5.23\% \\
          & $z_R^{MS}$ via OPT    & 13    & 9     & 21    & \textbf{\$4,425,950} & \$144,496 & \textbf{\$4,814,425} & 1661.48 & N.A. \\
          & $z_R^{MS}$ via AA & 14    & 8     & 16    & \$4,700,425 & \$141,694 & \$5,054,751 & \textbf{18.31} & 1.05 \\
          \hline
    \end{tabular}%
    }
  \label{tab:real-sol}%
\end{table}%

{
\subsubsection{Results with Prioritization}
To make all the candidate sites have the same capacity, we set $h_{ti}=h=10^6,\ \forall i\in[M],\ t\in[T]$. Similar to Section \ref{sec:setup}, we generate budget scenarios according to uniform distribution between $\lceil\max_{\omega\in\Omega}\{\sum_{j\in[N]}d_{t,j}(\omega)\}/h\rceil-\lceil\max_{\omega\in\Omega}\{\sum_{j\in[N]}d_{t-1,j}(\omega)\}/h\rceil$ and $\lceil\max_{\omega\in\Omega}\{\sum_{j\in[N]}d_{t,j}(\omega)\}/h\rceil-\lceil\max_{\omega\in\Omega}\{\sum_{j\in[N]}d_{t-1,j}(\omega)\}/h\rceil+2$. The results are presented in Table \ref{tab:real-sol-prior}, where in the third column, we record the top 6 candidate sites in the priority list $\boldsymbol{s}_0$, and the last column presents the computational time and optimality gaps without or with prioritization cuts. Comparing Tables \ref{tab:real-sol-prior} with \ref{tab:real-sol}, the models with prioritization are more computationally demanding than the ones without prioritization, where the prioritization cuts can reduce the computational time or optimality gaps within the one-hour time frame. Looking at the candidate sites with top priority, although the priority list changes significantly across different demand patterns, CA is always in the top 6 because California has the highest demand among all the customer sites according to \citet{daskin2011network}'s 88-node dataset.
\begin{table}[ht!]
  \centering
  \caption{Optimal solutions and cost breakdown of two-stage and multistage models with prioritization under different demand patterns}
  \resizebox{0.9\textwidth}{!}{
    \begin{tabular}{cllrrr}
    \hline
    Pattern & Model & Top 6 sites in $\boldsymbol{s}_0$    & Tans. (\$) & Obj. (\$)  & Time (sec.) \\
    \hline
    \multirow{4}[0]{*}{I} & $z_P^{TS}$ without cuts    & [CA, TX, OH, NJ, WI, AZ] & \$56,087 & \$65,141 & 3600 (0.93\%) \\
          & $z_P^{TS}$ with cuts & [CA, TX, OH, NJ, WA, CT] & \$55,955 & \$65,078 & 928.64 \\
          & $z_P^{MS}$ without cuts     & [CA, OH, TN, NJ, WI, WA] & \$55,355 & \$59,928 & 2127.05 \\
          & $z_P^{MS}$ with cuts & [NJ, TX, CT, CA, AZ, WI] & \$55,355 & \$59,928 & 516.71 \\
          \hline
    \multirow{4}[0]{*}{II} & $z_P^{TS}$ without cuts     & [NV, CA, TX, CT, DE, WI] & \$56,817 & \$67,269 & 3600 (2.16\%) \\
          & $z_P^{TS}$ with cuts & [CA, OH, NJ, KS, AZ, TX] & \$56,839 & \$67,247 & 1601.24 \\
          & $z_P^{MS}$ without cuts     & [CT, NJ, DE, TX, CA, KS] & \$56,677 & \$61,963 & 3600 (1.99\%) \\
          & $z_P^{MS}$ with cuts & [CA, TX, OH, NJ, WI, WA] & \$56,191 & \$61,797 & 1209.54 \\
          \hline
    \multirow{4}[0]{*}{III} & $z_P^{TS}$ without cuts     & [CA, NM, DC, CT, TX, KS] & \$66,930 & \$76,761 & 3600 (1.11\%) \\
          & $z_P^{TS}$ with cuts & [TX, NJ, CA, AZ, WA, WI] & \$66,930 & \$76,761 & 640.81 \\
          & $z_P^{MS}$ without cuts    & [DC, IN, KS, NM, LA, NJ] & \$66,505 & \$70,068 & 3600 (1.27\%) \\
          & $z_P^{MS}$ with cuts & [AZ, NJ, TX, CA, WA, LA] & \$66,505 & \$69,987 & 183.14 \\
          \hline
    \multirow{4}[0]{*}{IV} & $z_P^{TS}$ without cuts     & [LA, DE, CA, CT, AZ, GA] & \$75,074 & \$89,164 & 3600 (2.48\%) \\
          & $z_P^{TS}$ with cuts & [TX, DC, DE, NM, CT, KS] & \$74,941 & \$88,962 & 3600 (0.37\%) \\
          & $z_P^{MS}$ without cuts    & [CA, IN, NV, OH, CT, SC] & \$74,456 & \$81,271 & 3600 (2.08\%) \\
          & $z_P^{MS}$ with cuts & [AZ, CA, TX, NJ, WI, LA] & \$74,328 & \$81,028 & 660.94 \\
          \hline
    \end{tabular}%
    }
  \label{tab:real-sol-prior}%
\end{table}%
}

\section{Conclusion}
\label{sec:conclu}
We considered a class of multiperiod capacitated facility location problems under uncertain demand and budget in each period. When only demand is uncertain, we compared a multistage stochastic programming model where the locations of facilities can be determined dynamically throughout the uncertainty realization process, with a two-stage model where decision makers have to fix facility locations at the beginning of the horizon. {When both demand and budget are uncertain, we formed a rank-ordered list of all candidate facilities and make sure that the facility-selection decisions obey the priority list. In a two-stage model, the priority list is decided up front and fixed through stages, while the priority list can change adaptively in a multistage model.} Using expected conditional risk measures (ECRMs), we bounded the gaps between risk-averse two-stage and multistage optimal objective values with or without prioritization from below and proved that the lower bounds are tight. Two approximation algorithms are also proposed to solve the risk-averse models without prioritization more efficiently, which are asymptotically optimal under an expanding market. We also proposed prioritization cuts to speed up computation for risk-averse models with prioritization. Our numerical tests indicate that the $\rm{RVMS}$ increases as the uncertainty variation increases and stagewise dependent scenario trees attain higher $\rm{RVMS}$ than the stagewise independent counterparts. On the other hand,
the analytical bounds can recover the true gaps in many cases. Moreover, the approximation ratios are as low as 1.05 in the case study based on a real-world network, and the prioritization cuts can reduce computational time significantly.

There are several interesting directions to investigate for future research. The risk measure we used in this paper is the ECRMs, of which the risk is measured separately for each stage. There are other risk-measure choices that can be applied here, such as the nested risk measures. More research can be done to explore the relationship of these models. Moreover, this paper assumes that the {uncertainty} has a known distribution, while it is more realistic to assume that the probability distribution of {the uncertainty} belongs to an ambiguity set. Therefore, a distributionally robust optimization framework can be considered for two-stage and multistage facility location problems.  Depending on specific constructions of the ambiguity set, we may bound the gaps between the two-stage and multistage distributionally robust formulations.




\begin{appendices}
{
	\section{Literature Review and Our Contributions}
    \label{sec:literature}
We first review the literature on facility location and its variants, mainly focusing on the stochastic facility location with multiperiod uncertainty in Section \ref{sec:review-facility}. In Section \ref{sec:review-risk}, we introduce some popular risk measures used in the literature, including the ones we use in our later analysis. We end this section by reviewing several approximation algorithms proposed for solving facility location problems in Section \ref{sec:review-approx}.
 \subsection{Facility Location and Variants}
 \label{sec:review-facility}
	We refer the interested readers to \citet{owen1998strategic} for a comprehensive review on facility location variants, including the related static, dynamic, deterministic and stochastic models. The methods for solving the NP-hard facility location problems are mainly based on integer programming \citep{gendron2017comparison},  graph theories \citep{teitz1968heuristic}, and heuristics \citep{teitz1968heuristic,kaya2016mixed}). 
	To handle parameter uncertainty, robust optimization and stochastic programming are two approaches for achieving good performance in the worst case and on average, respectively. \citet{snyder2006facility}, \citet{correia2015facility}, and \citet{laporte2016location} reviewed facility location problems under demand uncertainty, summarizing both robust and stochastic facility location models. \citet{aghezzaf2005capacity} discussed stochastic strategic capacity planning and warehouse location problem in supply chains and proposed a decomposition algorithm based on Lagrangian relaxation to solve the robust optimization model. \citet{chan2017robust} proposed a data-driven optimization model for deploying automated external defibrillators in public spaces while accounting for uncertainty in future cardiac arrest locations. They considered both cases when the demand distribution is known and partially known, which leads to a distributionally robust optimization model. When the budget is uncertain, \cite{kocc2015prioritization} were the first to prioritize the activities under resource constraints, where they placed the activities in a rank-ordered list and select those with highest priority. The authors used a two-stage stochastic programming framework and applied the approach to a facility location problem and a knapsack problem. However, their work only considered budget uncertainty and involved a single decision period. To our best knowledge, there is no existing work that models prioritization under both budget and demand uncertainty in a multiperiod/multistage setting.
	
    \cite{nickel2015multi} provided a comprehensive review on the models and algorithms of multiperiod facility location. \citet{borison1984state} decomposed the multiperiod problem into a set of linked static deterministic problems, where the linkages were enforced through Lagrange multipliers, and they designed a primal-dual method to solve it.
	Recently, \citet{marin2018multi} formulated a multiperiod covering location problem as a two-stage mixed-integer linear program and developed a Lagrangian relaxation based heuristic to tackle large-scale problem instances. \cite{castro2017cutting} developed a cutting-plane approach for multiperiod capacitated facility location, in which they used a specialized interior-point method to solve the Benders subproblems.
	
Different from the rich studies on deterministic multiperiod facility location, there are a limited number of papers considering multistage stochastic facility location problems.  \cite{singh2009dantzig} formulated a multistage stochastic mixed-integer programming model for capacity-planning problems and  applied ``variable splitting'' and Dantzig-Wolfe decomposition to tackle the problem. \citet{nickel2012multi} considered a multiperiod facility location problem with random demand and return rate of investment decisions. The problem was formulated as a multistage stochastic mixed-integer linear program, which took extremely long time to solve (e.g., 6 hours by an off-the-shelf solver for 216-scenario instances). \citet{hernandez2012branch} studied a multiperiod facility location problem with stochastic demand. A large-scale model was solved approximately using a heuristic combining branch-and-fix coordination and branch-and-bound algorithm. \cite{albareda2013fix} proposed a fix-and-relax-coordination approximation procedure for tackling a multiperiod facility location problem with both cost and demand uncertainties.
	
	 \subsection{Risk Measures and Multistage Stochastic Programs}
	 \label{sec:review-risk}
 
    \citet{schultz2006conditional,shapiro2009lectures,ahmed2006convexity,miller2011risk} extended two-stage stochastic programs with risk-neutral expectation-based objective functions to risk-averse ones. However, it becomes nontrivial to model multistage risk-averse stochastic programs, as the risk could be measured separately for each stage or in a nested way. 
    \citet{shapiro2012minimax} explored the relations between the minimax, risk-averse and nested formulations of multistage stochastic programs.  
    \citet{philpott2012dynamic} incorporated a time-consistent coherent risk measure to a multistage stochastic programming model in a nested way, when the single-period risk measure was a convex combination of expectation and Conditional Value-at-Risk (CVaR). They derived a variant of Stochastic Dual Dynamic Programming (SDDP) algorithm by adding cuts to approximate recursive functions in the constraints and applied the algorithm to hydrothermal scheduling problems in New Zealand. \citet{philpott2013solving} formalized a general approximation procedure for computing solutions to multistage stochastic programming problems that minimize dynamic coherent risk measure in a nested form. 
    \citet{pflug2005measuring} proposed a class of multiperiod risk measures for a sequence of random incomes adapted to some filtration, which can be calculated by solving a stochastic dynamic linear optimization problem, and they analyzed its convexity and duality structure. \citet{homem2016risk} extended the above risk measures to expected conditional risk measures (ECRMs) and proved some appealing properties. First, ECRMs, originally defined for each stage separately, can be rewritten in a nested form. Second, any risk-averse multistage stochastic programs with ECRMs using CVaR measure can be recast as a simpler risk-neutral multistage stochastic program with additional variables and constraints, which can be efficiently solved using existing algorithms like SDDP. {Third, in this paper, we show that ECRMs are time consistent following the definition of \cite{ruszczynski2010risk}.} Due to these properties, in this paper we base our analysis on the ECRMs and also propose approximation algorithms for solving the resultant problems, which we will review the relevant literature next. 
    
    \subsection{Approximation Algorithms for Facility Location}
    \label{sec:review-approx}
	Motivated by the NP-hardness results, a stream of facility location research focuses on efficient approximation algorithms under different assumptions of network structures or demand patterns. Starting from \cite{shmoys1997approximation}, there has been a series of constant-factor approximation algorithms for either capacitated or uncapacitated deterministic facility location problems utilizing techniques from linear programming (LP) rounding, primal-dual methods, local search, greedy algorithms, etc. \citep{li20131,chudak2003improved,sviridenko2002improved,guha1999greedy}. A class of two-stage stochastic programming formulations have also captured researchers' attention \citep{ravi2006hedging,shmoys2006approximation,gupta2004boosted}, where some facilities can be built in the first stage, and additional facilities can be added in the second stage after observing demand realizations. Among them, \citet{gupta2004boosted} were the first to propose a black-box model and design approximation algorithms for various two-stage stochastic programming models, including the ones for uncapacitated facility location, but they assumed that the second-stage cost must be proportional to the first-stage one. Later, \cite{gupta2005wednesday} extended the prior work by proposing an algorithm for the $k$-stage stochastic program with the same cost restriction based on the boosted sampling framework, but their approximation ratio was \textit{exponential} in $k$. More recently, utilizing the black-box model but without the cost restriction, \cite{swamy2012sampling} successfully obtained an approximation ratio of $1.858(k-1)+1.52$ for the $k$-stage facility location problem.
	{Related to facility location with prioritization, another stream of work focuses on proposing approximation algorithms for incremental $k$-median problems. It aims to produce an ordering of the candidate facilities. For each $k$, consider the ratio of the cost of opening the first $k$ facilities in the ordering to the cost of the optimal $k$-median solution. The goal of the problem is to find an ordering such that the maximum of the ratio over all values of $k$ is minimized. \cite{mettu2003online} introduced the incremental $k$-median problem and gave a 29.86-competitive algorithm. \cite{plaxton2006approximation} introduced the incremental facility location problem and gave a $(4+\epsilon)\alpha$-competitive algorithm, given any $\alpha$-approximation algorithm for the uncapacitated facility location problem.}
	However, all of the above algorithms are designed for  \textit{uncapacitated} facility location problems as the capacity restriction will add another layer of complexity. 
	
	Different from the above studies, our paper focuses on solving risk-averse two-stage and multistage stochastic capacitated facility location problems, to which no one has proposed approximation algorithms to our best knowledge. \citet{huang2009value} studied a similar capacity expansion problem in a risk-neutral setting, and developed an asymptotically optimal approximation algorithm by exploiting a decomposable substructure inherent in the problem. Our approximation algorithms proceed along this line with some crucial differences. First, we derive analytical forms of the optimal solutions of the substructure problem in the risk-averse setting, which can further speed up the algorithm. Second, we repeat the process proposed in \citet{huang2009value} and prove that it can strengthen the upper bound iteratively. Finally, {when the demand is increasing over time, we prove that the approximation schemes are asymptotically optimal}.
	
\section{Time Consistency of ECRMs}\label{e-companion:ECRMs}
We follow the definition in \cite{ruszczynski2010risk} and prove the time consistency of ECRMs \eqref{ECRMs}.
Consider the probability space $(\Xi,\mathcal{F},P)$, and let $\mathcal{F}_1\subset\mathcal{F}_2\subset\ldots\subset\mathcal{F}_T$ be sub-sigma-algebras of $\mathcal{F}$ such that each $\mathcal{F}_t$ corresponds to the information available up to (and including) stage $t$, with $\mathcal{F}_1=\{\emptyset,\Xi\}, \ \mathcal{F}_T=\mathcal{F}$. Let $\mathcal{Z}_t$ denote a space of $\mathcal{F}_t$-measurable functions from $\Xi$ to $\mathbb{R}$, and let $\mathcal{Z}_{1,T}:=\mathcal{Z}_1\times\cdots\times\mathcal{Z}_T$.
\begin{definition}
A mapping $\rho_{t,T}:\mathcal{Z}_{t,T}\to\mathcal{Z}_t$ where $1\le t\le T$ is called a \textit{conditional risk measure}, if it has the following monotonicity property:
$\rho_{t,T}(Z)\le\rho_{t,T}(W)$ for all $Z,W\in\mathcal{Z}_{t,T}$ such that $Z\le W$.
\end{definition}
\begin{definition}
A dynamic risk measure is a sequence of monotone one-step conditional risk measures $\rho_{t,T}: \mathcal{Z}_{t,T}\to\mathcal{Z}_t,\ 1\le t\le T$.
\end{definition}
\begin{definition}\label{def:time-consistent}
A dynamic risk measure $\{\rho_{t,T}\}_{t=1}^{T}$ is called time consistent if, for all $1\le l<k\le T$ and all sequences $Z, \ W\in \mathcal{Z}_{l, T}$, the conditions
\begin{align*}
    Z_i = W_i,\ \forall i=l,\ldots, k-1, \text{and}\ 
    \rho_{k,T}(Z_k,\ldots,Z_T)\le \rho_{k,T}(W_k,\ldots,W_T)
\end{align*}
imply that
\begin{align*}
    \rho_{l,T}(Z_l,\ldots,Z_T)\le 
     \rho_{l,T}(W_l,\ldots,W_T).
\end{align*}
\end{definition}
\begin{theorem}\label{thm:ECRMs}
The ECRMs defined in \eqref{ECRMs} are time consistent, if each $\rho_t^{d_{[1,t-1]}}$ is a coherent one-step conditional risk measure.
\end{theorem}
\proof{Proof of Theorem \ref{thm:ECRMs}}
According to Eq. \eqref{eq:obj} and the translation-invariant property of $\rho_t^{d_{[1,t-1]}}$, the risk function \eqref{ECRMs} can be recast as
\begin{align*}
\mathbb{F}(g_1,\ldots,g_{T})=g_1+\rho_2\Big(g_2+\mathbb{E}_{\boldsymbol d_{2}}\circ{\rho_3^{d_{[1,2]}}}\Big(g_3+\mathbb{E}_{\boldsymbol d_{3}|d_{[1,2]}}\circ{\rho_4^{d_{[1,3]}}}\Big(g_4+\cdots\nonumber\\
+\mathbb{E}_{\boldsymbol d_{T-1}|d_{[1,T-2]}}\circ{\rho_{T}^{d_{[1,T-1]}}}\Big(g _{T}\Big)\Big)\cdots\Big)\Big),
\end{align*}
To simplify the notation, we define $\tilde{\rho}_t^{d_{[1,t-2]}}: = \mathbb{E}_{\boldsymbol{d}_{t-1}|d_{[1,t-2]}}\circ \rho_t^{d_{[1,t-1]}}$, which maps from $\mathcal{Z}_t$ to $\mathcal{Z}_{t-2}$. Then, the multiperiod risk function $\mathbb{F}$ can be recast as
\begin{align*}
\mathbb{F}(g_1,\ldots,g_{T})=g_1+\rho_2\Big(g_2+\tilde{\rho}_3^{d_{[1]}}\Big(g_3+\tilde{\rho}_4^{d_{[1,2]}}\Big(g_4+\cdots
+\tilde{\rho}_T^{d_{[1,T-2]}}\Big(g _{T}\Big)\Big)\cdots\Big)\Big),
\end{align*}
Define a dynamic risk measure $\{\mathbb{F}_{t,T}\}_{t=1}^T$ as
$\mathbb{F}_{1,T}=\mathbb{F}$, and for $2\le t\le T$, we have
\begin{align*}
    \mathbb{F}_{t,T}(g_t,\ldots,g_T) = g_t+\tilde{\rho}_{t+1}^{d_{[1,t-1]}}\Big(g_{t+1}+\tilde{\rho}_{t+2}^{d_{[1,t]}}\Big(g_{t+2}+\cdots
+\tilde{\rho}_T^{d_{[1,T-2]}}\Big(g _{T}\Big)\Big)\cdots\Big)
\end{align*}
Then for $1\le l < k\le T$,
\begin{align*}
    \mathbb{F}_{l,T}(g_l,\ldots,g_T) =g_l + \tilde{\rho}_{l+1}^{d_{[1,l-1]}}\Big(g_{l+1} +\ldots+\tilde{\rho}_k^{d_{[1,k-2]}}\Big(\mathbb{F}_{k,T}(g_k,\ldots,g_T)\Big)\Big).
\end{align*}\normalsize
If $\mathbb{F}_{k,T}(Z_k,\ldots,Z_T)\le \mathbb{F}_{k,T}(W_k,\ldots,W_T)$ and $Z_i=W_i,\ \forall i=l,\ldots, k-1$, we have
$\mathbb{F}_{l,T}(Z_l,\ldots,Z_T)\le \mathbb{F}_{l,T}(W_l,\ldots,W_T)$ because of the monotonicity of $\tilde{\rho}_t^{d_{[1,t-2]}},\ \forall t\ge {l+1}$.
This completes the proof.
$\hfill\Box$
\endproof
}
\section{Details of All Needed Proofs}\label{e-companion:proofs}

\proof{Proof of Theorem \ref{thm:extensive}}
Interchanging the expectation and minimum operator, and then merging the minimum on $\boldsymbol{x}_t,\boldsymbol{y}_t,u_t$ and the minimum on $\eta_{t+1}$, we can rewrite \eqref{model:multirisk} as
\begin{align}
&\min_{\substack{{(\boldsymbol{x}_1,\boldsymbol{y}_1)\in {X}_1(d_1)},\\\eta_2\in\mathbb{R}}}g_1+\lambda_2\eta_2+\mathbb{E}_{\boldsymbol d_2}\Big[\min_{\substack{(\boldsymbol{x}_2,\boldsymbol{y}_2,u_2)\in\tilde{X}_2(\boldsymbol{x}_1,\eta_2,d_2),\\\eta_3\in\mathbb{R}}}\Big\{\frac{\lambda_2}{1-\alpha_2}u_2+(1-\lambda_2)g_2+\lambda_3\eta_3+\ldots\nonumber\\
&+\mathbb{E}_{\boldsymbol d_{T-1}{|\boldsymbol{d}_{[1,T-2]}}}\Big[\min_{\substack{(\boldsymbol{x}_{T-1},\boldsymbol{y}_{T-1},u_{T-1})\in\tilde{X}_{T-1}(\boldsymbol{x_{1:T-2}},\eta_{T-1},d_{T-1}),\\\eta_{T}\in\mathbb{R}}}\Big\{\frac{\lambda_{T-1}}{1-\alpha_{T-1}}u_{T-1}+(1-\lambda_{T-1})g_{T-1}+\lambda_T\eta_T\nonumber\\
&+\mathbb{E}_{\boldsymbol d_T{|{\boldsymbol{d}_{[1,T-1]}}}}\Big[\min_{(\boldsymbol{x}_T,\boldsymbol{y}_T,u_T)\in\tilde{X}_T(\boldsymbol{x_{1:T-1}},\eta_T,d_T)}\Big\{\frac{\lambda_T}{1-\alpha_T}u_T+(1-\lambda_T)g_T\Big\}\Big]\ldots\Big\}\Big],\label{model:multirisk2}
\end{align}
where 
${X}_1(\boldsymbol{d}_1)=\{(\boldsymbol{x}_1, \boldsymbol{y}_1)\in  \mathbb{Z}^{M}_{+}\times \mathbb{R}^{M\times N}_{+}:  \boldsymbol{Ay}_{1}= \boldsymbol{d}_{1},\ \boldsymbol{B}_{1}\boldsymbol{y}_{1}\le \boldsymbol{x}_1, \ \boldsymbol{x}_1\le \mathbf{1} \}$ and $\tilde{X}_t(\boldsymbol{x_{1:t-1}},\eta_t,\boldsymbol{d}_t)={X}_t(\boldsymbol{x_{1:t-1}},\boldsymbol{d}_t)\cap\{u_t\in \mathbb{R}_+: u_t+\eta_t\ge g_t\}=\{(\boldsymbol{x}_t, \boldsymbol{y}_t)\in  \mathbb{Z}^{M}_{+}\times \mathbb{R}^{M\times N}_{+}:  \boldsymbol{A}\boldsymbol{y}_t= {\boldsymbol{d}_{t}},\ \boldsymbol{B}_t\boldsymbol{y}_t-\boldsymbol{x}_t\le \sum_{\tau=1}^{t-1}\boldsymbol{x}_{\tau}, \ \boldsymbol{x}_t\le \mathbf{1}-\sum_{\tau=1}^{t-1}\boldsymbol{x}_{\tau}\}\cap\{u_t\in \mathbb{R}_+: u_t+\eta_t\ge g_t\},\ \forall t=2,\ldots,T$ are the feasible sets in each stage. Plugging $g_t = \boldsymbol{f}^{\mathsf T}_t \sum_{\tau=1}^t \boldsymbol{x}_{\tau}+\boldsymbol{c}_t^{\mathsf T} \boldsymbol{y}_t$ into the formulation \eqref{model:multirisk2} and using scenario-node-based notation, we derive the equivalent extensive form \eqref{model:multiriskaverse}.
$\hfill\Box$
\endproof

\proof{Proof of Lemma \ref{lemma:equivalence}}
Observing that the right-hand sides of constraints \eqref{eq:two-stage-constraint_eta} and \eqref{eq:constraint_eta} are different, we make a variable transformation in \eqref{model:tworiskaverse-scenario-node} --- $\tilde{\eta}_{a(n)}:=\eta_n+\boldsymbol{f}^{\mathsf T}_{t_n}\sum_{m\in \mathcal{P}(n)}\boldsymbol{x}_m, \ \forall n\not=1$. Because both $\eta$ and $\boldsymbol{x}$ are first-stage variables and are identical for all nodes in the same stage, this transformation is without of ambiguity and is equivalent to say $\tilde{\eta}_{t-1}:=\eta_t+\boldsymbol{f}^{\mathsf T}_{t}\sum_{\tau=1}^t\boldsymbol{x}_{\tau}$. With that, the two-stage constraints \eqref{eq:twostage_eta} for the new variable $\tilde{\eta}_n$ are satisfied and constraints \eqref{eq:two-stage-constraint_eta} become \eqref{eq:constraint_eta}.
Moreover, the objective function can be transformed to 
\begin{align*}
&\sum_{n\in \mathcal{T}} p_n\left(\boldsymbol{f}^{\mathsf T}_{t_n}\sum_{m\in \mathcal{P}(n)}\boldsymbol{x}_m\right) + \sum_{n\in \mathcal{T}: n\not=1} p_n\left((1-\lambda_{t_n})\boldsymbol{c}_{t_n}^{\mathsf T}\boldsymbol{y}_n+{\lambda_{t_n}} \eta_n + \frac{\lambda_{t_n}}{1-\alpha_{t_n}} u_n\right) + \boldsymbol{c}_{t_1}^{\mathsf T} \boldsymbol{y}_1\\
=& \sum_{n\in \mathcal{T}} p_n\left(\boldsymbol{f}^{\mathsf T}_{t_n}\sum_{m\in \mathcal{P}(n)}\boldsymbol{x}_m\right) + \sum_{n\in \mathcal{T}: n\not=1} p_n\left((1-\lambda_{t_n})\boldsymbol{c}_{t_n}^{\mathsf T}\boldsymbol{y}_n+{\lambda_{t_n}} \left(\tilde{\eta}_{a(n)} - \boldsymbol{f}^{\mathsf T}_{t_n}\sum_{m\in \mathcal{P}(n)}\boldsymbol{x}_m\right) + \frac{\lambda_{t_n}}{1-\alpha_{t_n}} u_n\right) + \boldsymbol{c}_{t_1}^{\mathsf T} \boldsymbol{y}_1\\
=& \sum_{n\in \mathcal{T}} p_n\left(\boldsymbol{\tilde{f}}_n^{\mathsf T}\sum_{m\in \mathcal{P}(n)}\boldsymbol{x}_m+\boldsymbol{\tilde{c}}_n^{\mathsf T}\boldsymbol{y}_n + \tilde{\alpha}_n u_n\right) + \sum_{n\in \mathcal{T}: n\not=1} p_n{\lambda_{t_n}} \tilde{\eta}_{a(n)}\\
=& \sum_{n\in \mathcal{T}} p_n\left(\boldsymbol{\tilde{f}}_n^{\mathsf T}\sum_{m\in \mathcal{P}(n)}\boldsymbol{x}_m+\boldsymbol{\tilde{c}}_n^{\mathsf T}\boldsymbol{y}_n + \tilde{\alpha}_n u_n\right) + \sum_{n\in \mathcal{T}: n\not\in\mathcal{L}} p_n{\lambda_{t_n+1}} \tilde{\eta}_{n}
\end{align*}
where the last equality holds because $\sum_{m\in\mathcal{C}(n)}p_m\lambda_{t_m} = p_n\lambda_{t_n+1}$. This completes the proof. $\hfill \Box$
\endproof

\proof{Proof of Proposition \ref{prop:substructure}}
We first check the feasibility of solutions $\{\boldsymbol{x}^{MS}_{n}\}_{n\in \mathcal{T}},\ \{\eta_n^{MS}\}_{n\in \mathcal{T}\setminus\mathcal{L}}$ to $\mbox{{\bf SP-RMS}}(\boldsymbol{y}_n^{*}, u_n^{*})$. Note that for all $n\in\mathcal{T}$, we have $\boldsymbol{x}^{MS}_{n} \in \mathbb{Z}_{+}^{M}$ and $\sum_{m\in \mathcal{P}(n)}\boldsymbol{x}^{MS}_m = \max_{m\in\mathcal{P}(n)}\lceil \boldsymbol{B}_{t_m} \boldsymbol{y}_m^{*}\rceil \ge \boldsymbol{B}_{t_n} \boldsymbol{y}_n^{*}$, which also satisfy constraints \eqref{eq:constraint_x_knowny} because of $\boldsymbol{B}_{t_n}\boldsymbol{y}_n^* \le 1$. By the definition of $\eta_n^{MS}$, we have ${\eta}^{MS}_n \ge  \boldsymbol{f}^{\mathsf T}_{t_m}\sum_{l\in\mathcal{P}(m)}\boldsymbol{x}^{MS}_l+ \boldsymbol{c}_{t_m}^{\mathsf T} \boldsymbol{y}_m^{*}-u_m^{*},\ \forall m\in \mathcal{C}(n),\ n\not\in\mathcal{L}$, which satisfy constraints \eqref{eq:constraint_eta_knowny} automatically. As a result, we obtain $Q^M(\boldsymbol{y}_n^{*}, u_n^{*}) \le \sum_{n\in \mathcal{T}} p_n\left(\boldsymbol{\tilde{f}}_n^{\mathsf T}\sum_{m\in \mathcal{P}(n)}\boldsymbol{x}^{MS}_m+ \tilde{\lambda}_n \eta^{MS}_n\right) = \sum_{n\in \mathcal{T}} p_n\left(\boldsymbol{\tilde{f}}_n^{\mathsf T}\max_{m\in\mathcal{P}(n)}\lceil \boldsymbol{B}_{t_m} \boldsymbol{y}_m^{*}\rceil+ \tilde{\lambda}_n \max_{m\in \mathcal{C}(n)}\left\lbrace \boldsymbol{f}^{\mathsf T}_{t_m}\max_{l\in\mathcal{P}(m)}\lceil \boldsymbol{B}_{t_l} \boldsymbol{y}_l^{*}\rceil+ \boldsymbol{c}_{t_m}^{\mathsf T} \boldsymbol{y}_m^{*}-u_m^{*}\right\rbrace\right)$.

Next, we show the optimality of solutions $\{\boldsymbol{x}^{MS}_{n}\}_{n\in \mathcal{T}},\ \{\eta_n^{MS}\}_{n\in \mathcal{T}\setminus\mathcal{L}}$ to $\mbox{{\bf SP-RMS}}(\boldsymbol{y}_n^{*}, u_n^{*})$. Note that for any feasible solution $\boldsymbol{x}_{n}\in \mathbb{Z}_+, \ \forall n\in\mathcal{T}$, from constraints \eqref{eq:constraint_xy_knowny}, we get $\sum_{l\in \mathcal{P}(n)}\boldsymbol{x}_l\ge \sum_{l\in \mathcal{P}(m)}\boldsymbol{x}_l\ge \lceil\boldsymbol{B}_{t_m}\boldsymbol{y}_m^*\rceil,\ \forall m\in \mathcal{P}(n)$ (we can raise $\boldsymbol{B}_{t_m}\boldsymbol{y}_m^*$ to $\lceil\boldsymbol{B}_{t_m}\boldsymbol{y}_m^*\rceil$ because of the integrality of $\boldsymbol{x}_n$-variables) and thus $\sum_{l\in \mathcal{P}(n)}\boldsymbol{x}_l\ge \max_{m\in \mathcal{P}(n)}\lceil\boldsymbol{B}_{t_m}\boldsymbol{y}_m^*\rceil$. From constraints \eqref{eq:constraint_eta_knowny}, we obtain ${\eta}_n \ge  \max_{m\in \mathcal{C}(n)}\{\boldsymbol{f}^{\mathsf T}_{t_m}\sum_{l\in\mathcal{P}(m)}\boldsymbol{x}_l+ \boldsymbol{c}_{t_m}^{\mathsf T} \boldsymbol{y}_m^{*}-u_m^{*}\}\ge \max_{m\in \mathcal{C}(n)}\{\boldsymbol{f}^{\mathsf T}_{t_m}\max_{l\in \mathcal{P}(m)}\lceil\boldsymbol{B}_{t_l}\boldsymbol{y}_l^*\rceil+ \boldsymbol{c}_{t_m}^{\mathsf T} \boldsymbol{y}_m^{*}-u_m^{*}\},\ \forall n\not\in\mathcal{L}$. With that, we conclude $Q^M(\boldsymbol{y}_n^{*}, u_n^{*}) \ge \sum_{n\in \mathcal{T}} p_n\left(\boldsymbol{\tilde{f}}_n^{\mathsf T}\max_{m\in\mathcal{P}(n)}\lceil \boldsymbol{B}_{t_m} \boldsymbol{y}_m^{*}\rceil+ \tilde{\lambda}_n \max_{m\in \mathcal{C}(n)}\left\lbrace \boldsymbol{f}^{\mathsf T}_{t_m}\max_{l\in\mathcal{P}(m)}\lceil \boldsymbol{B}_{t_l} \boldsymbol{y}_l^{*}\rceil+ \boldsymbol{c}_{t_m}^{\mathsf T} \boldsymbol{y}_m^{*}-u_m^{*}\right\rbrace\right)$.

Therefore, $\{\boldsymbol{x}^{MS}_{n}\}_{n\in \mathcal{T}},\ \{\eta_n^{MS}\}_{n\in \mathcal{T}\setminus\mathcal{L}}$ are optimal solutions to the risk-averse multistage problem \eqref{eq:ms-substructure}. In the risk-averse two-stage problem \eqref{eq:ts-substructure}, because the investment decisions $\boldsymbol{x}$ are identical for all nodes in the same stage, one only needs to replace the right-hand side of constraints \eqref{eq:constraint_xy_knowny} with $\max_{l\in \mathcal{T}_{t_n}}\boldsymbol{B}_{t_l} \boldsymbol{y}_l^{*}$. Moreover, to satisfy the two-stage constraints for $\boldsymbol{\eta}$, we need to ensure that ${\eta}_n \ge  \max_{m\in \mathcal{T}_{t_n+1}}\{\boldsymbol{f}^{\mathsf T}_{t_m}\sum_{l\in\mathcal{P}(m)}\boldsymbol{x}_l+ \boldsymbol{c}_{t_m}^{\mathsf T} \boldsymbol{y}_m^{*}-u_m^{*}\},\ \forall n\not\in\mathcal{L}$. Following the same analysis used for risk-averse multistage problem \eqref{eq:ms-substructure} and replacing $\boldsymbol{B}_{t_n} \boldsymbol{y}_n^{*}$ with $\max_{l\in \mathcal{T}_{t_n}}\boldsymbol{B}_{t_l} \boldsymbol{y}_l^{*}$, it can be shown that $\{\boldsymbol{x}^{TS}_{n}\}_{n\in \mathcal{T}},\ \{\eta_n^{TS}\}_{n\in \mathcal{T}\setminus\mathcal{L}}$ are optimal solutions to the risk-averse two-stage problem \eqref{eq:ts-substructure}. This completes the proof. $\hfill \Box$
\endproof

\proof{Proof of Theorem \ref{thm:risk-averse-bound}}
First of all, $\rm VMS_{R}^{LB}\ge 0$ and thus provides a nontrivial lower bound. Let $\left\lbrace {\boldsymbol{y}_n^{*}}\right\rbrace_{{n\in\mathcal{T}}},\ \{u_n^{*}\}_{n\in \mathcal{T}\setminus\{1\}}$ be the values of decisions made in an optimal solution to the two-stage model \eqref{model:tworiskaverse}. We have $
	z^{TS}_R
    =\sum_{n\in \mathcal{T}} p_n(\boldsymbol{\tilde{c}}_n^{\mathsf T}\boldsymbol{y^{*}_n} + \tilde{\alpha}_n u^{*}_n)+Q^{T}(\boldsymbol{y}_n^{*},u_n^{*})$,
	where $Q^{T}(\boldsymbol{y}_n^{*},u_n^{*})$ is the optimal objective value to the substructure problem defined in \eqref{eq:ts-substructure}, and this equality is true because of the optimality of decisions $\boldsymbol{y}_n^{*},\ u_n^{*}$.

Since $\left\lbrace {\boldsymbol{y}_n^{*}}\right\rbrace_{{n\in\mathcal{T}}},\ \{u_n^{*}\}_{n\in \mathcal{T}\setminus\{1\}}$ are feasible solutions for the multistage problem \eqref{model:multiriskaverse}, we have
\begin{align}
z^{MS}_R\le \sum_{n\in \mathcal{T}} p_n(\boldsymbol{\tilde{c}}_n^{\mathsf T}\boldsymbol{y^{*}_n} + \tilde{\alpha}_n u^{*}_n)+Q^{M}(\boldsymbol{y}_n^{*},u_n^{*}),\label{eq:ineq}
\end{align}
where $Q^{M}(\boldsymbol{y}_n^{*},u_n^{*})$ is the optimal objective value to the substructure problem defined in \eqref{eq:ms-substructure}.

Combining all steps above, we have
\begin{align*}
{\rm VMS_R}=&z^{TS}_R-z^{MS}_R\\
\ge& Q^{T}(\boldsymbol{y}_n^{*},u_n^{*})-Q^M(\boldsymbol{y}_n^{*},u_n^{*})\\
\overset{(a)}{=}& \sum_{n\in \mathcal{T}} p_n\left(\boldsymbol{\tilde{f}}_n^{\mathsf T}\left(\max_{m\in \mathcal{P}(n)} \lceil \max_{l\in \mathcal{T}_{t_m}}\boldsymbol{B}_{t_l} \boldsymbol{y}_l^{*}\rceil - \max_{m\in \mathcal{P}(n)}\lceil \boldsymbol{B}_{t_m} \boldsymbol{y}_m^{*}\rceil\right)+ \tilde{\lambda}_n \left(\eta^{TS}_n-\eta^{MS}_n\right)\right)\\
=& {\sum_{n\in \mathcal{T}\setminus\{1\}} p_n(1-\lambda_{t_n})\boldsymbol{f}_{t_n}^{\mathsf T}\left(\max_{m\in \mathcal{P}(n)} \lceil \max_{l\in \mathcal{T}_{t_m}}\boldsymbol{B}_{t_l} \boldsymbol{y}_l^{*}\rceil - \max_{m\in \mathcal{P}(n)}\lceil \boldsymbol{B}_{t_m} \boldsymbol{y}_m^{*}\rceil\right)+ {\sum_{n\in \mathcal{T}\setminus\mathcal{L}} p_n {\lambda}_{t_n+1}} \left(\eta^{TS}_n-\eta^{MS}_n\right)}
\end{align*}
where the equality $(a)$ follows Proposition \ref{prop:substructure}.
This completes the proof. $\hfill \Box$
\endproof

{
\proof{Proof of Corollary \ref{cor:LP}}
Because $\{\boldsymbol{y}_n^{LP}\}_{n\in \mathcal{T}}, \ \{u_n^{LP}\}_{n\in \mathcal{T}\setminus\{1\}}$ is an optimal solution to the LP relaxation of the two-stage model \eqref{model:tworiskaverse}, we have $
	z^{TS}_R
    \ge\sum_{n\in \mathcal{T}} p_n(\boldsymbol{\tilde{c}}_n^{\mathsf T}\boldsymbol{y^{LP}_n} + \tilde{\alpha}_n u^{LP}_n)+Q^{TLP}(\boldsymbol{y}_n^{LP},u_n^{LP})$,
where $Q^{TLP}(\boldsymbol{y}_n^{LP},u_n^{LP})$ stands for the optimal objective value of the LP relaxation of the subproblem $\mbox{{\bf SP-RTS}}(\boldsymbol{y}_n^{LP}, u_n^{LP})$. It can be easily verified that the constructed $\boldsymbol{x}^{TS},\boldsymbol{\eta}^{TS}$ are optimal by replacing all $\lceil\boldsymbol{B}_{t_n} \boldsymbol{y}_n^{LP}\rceil$ with $\boldsymbol{B}_{t_n}\boldsymbol{y}_n^{LP}$, and thus $Q^{TLP}(\boldsymbol{y}_n^{LP},u_n^{LP})=\sum_{n\in \mathcal{T}} p_n\left(\boldsymbol{\tilde{f}}_n^{\mathsf T}\sum_{m\in \mathcal{P}(n)}\boldsymbol{x}^{TS}_m+ \tilde{\lambda}_n \eta^{TS}_n\right)$.

Since $\left\lbrace {\boldsymbol{y}_n^{LP}}\right\rbrace_{{n\in\mathcal{T}}},\ \{u_n^{LP}\}_{n\in \mathcal{T}\setminus\{1\}}$ are feasible solutions for the multistage problem \eqref{model:multiriskaverse}, we have $
z^{MS}_R\le \sum_{n\in \mathcal{T}} p_n(\boldsymbol{\tilde{c}}_n^{\mathsf T}\boldsymbol{y^{LP}_n} + \tilde{\alpha}_n u^{LP}_n)+Q^{M}(\boldsymbol{y}_n^{LP},u_n^{LP})$,
where $Q^{M}(\boldsymbol{y}_n^{LP},u_n^{LP})$ is the optimal objective value to the substructure problem \eqref{eq:ms-substructure}. The rest of the proof follows the proof of Theorem \ref{thm:risk-averse-bound}.
$\hfill\Box$
\endproof

\proof{Proof of Corollary \ref{cor:parameter}} Notice that
	\begin{align*}
	{\rm VMS_R^{LB}}&\ge  \sum_{n\in \mathcal{T}\setminus\{1\}} p_n(1-\lambda_{t_n})\boldsymbol{{f}}_{t_n}^{\mathsf T}\left(\max_{m\in \mathcal{P}(n)} \lceil \max_{l\in \mathcal{T}_{t_m}}\boldsymbol{B}_{t_l} \boldsymbol{y}_l^{*}\rceil - \max_{m\in \mathcal{P}(n)}\lceil \boldsymbol{B}_{t_m} \boldsymbol{y}_m^{*}\rceil\right)\\
	&=\sum_{n\in \mathcal{T}\setminus\{1\}} p_n(1-\lambda_{t_n})\sum_{i=1}^Mf_{t_n,i}\left(\max_{m\in \mathcal{P}(n)} \lceil \max_{l\in \mathcal{T}_{t_m}}\frac{1}{h_{t_l,i}}\sum_{j=1}^Ny_{l,ij}^*\rceil - \max_{m\in \mathcal{P}(n)}\lceil \frac{1}{h_{t_m,i}}\sum_{j=1}^Ny_{m,ij}^*\rceil\right)\\
	&\overset{(a)}{\ge}\sum_{n\in \mathcal{T}\setminus\{1\}} p_n(1-\lambda_{t_n} )\sum_{i=1}^M f_{t_n,i}\Delta_{ni}
	\end{align*}
Here, $(a)$ is true because when Condition (i) holds, we have $\max_{m\in \mathcal{P}(n)}\lceil \frac{1}{h_{t_m,i}}\sum_{j=1}^Ny_{m,ij}^*\rceil=0,\ \forall i\in[M]$; when Condition (ii) holds, we have $\frac{1}{h_{t_{\bar{n}},i}}\sum_{j=1}^Ny_{\bar{n},ij}^*>0$ (otherwise, the rest of the facilities cannot cover the total demand) and thus $\max_{m\in \mathcal{P}(n)} \lceil \max_{l\in \mathcal{T}_{t_m}}\frac{1}{h_{t_l,i}}\sum_{j=1}^Ny_{l,ij}^*\rceil=1$. As a result, when $\Delta_{ni}=1$, we have $\max_{m\in \mathcal{P}(n)} \lceil \max_{l\in \mathcal{T}_{t_m}}\frac{1}{h_{t_l,i}}\sum_{j=1}^Ny_{l,ij}^*\rceil - \max_{m\in \mathcal{P}(n)}\lceil \frac{1}{h_{t_m,i}}\sum_{j=1}^Ny_{m,ij}^*\rceil=1\ge \Delta_{ni}$ and $\max_{m\in \mathcal{P}(n)} \lceil \max_{l\in \mathcal{T}_{t_m}}\frac{1}{h_{t_l,i}}\sum_{j=1}^Ny_{l,ij}^*\rceil - \max_{m\in \mathcal{P}(n)}\lceil \frac{1}{h_{t_m,i}}\sum_{j=1}^Ny_{m,ij}^*\rceil$ is always non-negative.
$\hfill\Box$
\endproof
}

\proof{Proof of Theorem \ref{thm:hardness}}
We first note that any instance of the NP-hard knapsack problem with $M$ items can be polynomially transformed to the deterministic facility location problem \eqref{model:deter} with $T=1$. Moreover, Model \eqref{model:deter} is a single-scenario version of the risk-averse two-stage \eqref{model:tworiskaverse} and multistage \eqref{model:multiriskaverse} counterparts with $\lambda_t=0,\ \forall t=2,\ldots, T$. With that, we prove the NP-hardness of the aforementioned models.
$\hfill \Box$
\endproof

\proof{Proof of Proposition \ref{prop:alg}}
First of all, at the end of each iteration $k\ge 1$, $(\boldsymbol{x}_n^{k}, \eta_n^{k}, \boldsymbol{y}_n^{k}, u_n^{k})$ constitutes a feasible solution to the risk-averse multistage problem \eqref{model:multiriskaverse} and thus provides an upper bound. We have
\begin{align*}
    & z_R^{MS}(\boldsymbol{x}_n^{k+1}, \eta_n^{k+1}, \boldsymbol{y}_n^{k+1}, u_n^{k+1}) - z_R^{MS}(\boldsymbol{x}_n^{k}, \eta_n^{k}, \boldsymbol{y}_n^{k}, u_n^{k})\\
    = & \left(z_R^{MS}(\boldsymbol{x}_n^{k+1}, \eta_n^{k+1}, \boldsymbol{y}_n^{k+1}, u_n^{k+1}) - 
    z_R^{MS}(\boldsymbol{x}_n^{k+1}, \eta_n^{k+1}, \boldsymbol{y}_n^{k}, u_n^{k})\right) +
    \left(z_R^{MS}(\boldsymbol{x}_n^{k+1}, \eta_n^{k+1}, \boldsymbol{y}_n^{k}, u_n^{k}) -
    z_R^{MS}(\boldsymbol{x}_n^{k}, \eta_n^{k}, \boldsymbol{y}_n^{k}, u_n^{k})\right)\\
    \le & 0.
\end{align*}
The first term $z_R^{MS}(\boldsymbol{x}_n^{k+1}, \eta_n^{k+1}, \boldsymbol{y}_n^{k+1}, u_n^{k+1}) - 
    z_R^{MS}(\boldsymbol{x}_n^{k+1}, \eta_n^{k+1}, \boldsymbol{y}_n^{k}, u_n^{k})$ is non-positive because $(\boldsymbol{y}_n^{k+1}, u_n^{k+1})$ is an optimal solution to Step \ref{alg:step5}, while $(\boldsymbol{y}_n^{k}, u_n^{k})$ is a feasible solution to Step \ref{alg:step5} due to Step \ref{alg:step4}; the second term $z_R^{MS}(\boldsymbol{x}_n^{k+1}, \eta_n^{k+1}, \boldsymbol{y}_n^{k}, u_n^{k}) -
    z_R^{MS}(\boldsymbol{x}_n^{k}, \eta_n^{k}, \boldsymbol{y}_n^{k}, u_n^{k})$ is non-positive because $(\boldsymbol{x}_n^{k+1}, \eta_n^{k+1})$ is an optimal solution to Step \ref{alg:step4}, while $(\boldsymbol{x}_n^{k}, \eta_n^{k})$ is feasible to Step \ref{alg:step4} due to Step \ref{alg:step5} at the previous iteration. This completes the proof. $\hfill \Box$
\endproof

\proof{Proof of Proposition \ref{prop:ratio}}
First of all, we have
\begin{align*}
&z_R^{MS}(\boldsymbol{x}_n^{H}, \eta_n^{H}, \boldsymbol{y}_n^{H}, u_n^{H}) - z_R^{MS}(\boldsymbol{x}_n^{*}, \eta_n^{*}, \boldsymbol{y}_n^{*}, u_n^{*})\\
\overset{(a)}{\le} & z_R^{MS}(\boldsymbol{x}_n^{H}, \eta_n^{H}, \boldsymbol{y}_n^{H}, u_n^{H}) - z_R^{MS}(\boldsymbol{x}_n^{LP}, \eta_n^{LP}, \boldsymbol{y}_n^{LP}, u_n^{LP})\\
\overset{(b)}{\le} & z_R^{MS}(\boldsymbol{x}_n^{1}, \eta_n^{1}, \boldsymbol{y}_n^{1}, u_n^{1}) - z_R^{MS}(\boldsymbol{x}_n^{LP}, \eta_n^{LP}, \boldsymbol{y}_n^{LP}, u_n^{LP})\\
= & z_R^{MS}(\boldsymbol{x}_n^{1}, \eta_n^{1}, \boldsymbol{y}_n^{1}, u_n^{1}) - 
z_R^{MS}(\boldsymbol{x}_n^{1}, \eta_n^{1}, \boldsymbol{y}_n^{LP}, u_n^{LP}) +
z_R^{MS}(\boldsymbol{x}_n^{1}, \eta_n^{1}, \boldsymbol{y}_n^{LP}, u_n^{LP}) -
z_R^{MS}(\boldsymbol{x}_n^{LP}, \eta_n^{LP}, \boldsymbol{y}_n^{LP}, u_n^{LP})\\
\overset{(c)}{\le} & z_R^{MS}(\boldsymbol{x}_n^{1}, \eta_n^{1}, \boldsymbol{y}_n^{LP}, u_n^{LP}) -
z_R^{MS}(\boldsymbol{x}_n^{LP}, \eta_n^{LP}, \boldsymbol{y}_n^{LP}, u_n^{LP})\\
= & Q^M(\boldsymbol{y}_n^{LP}, u_n^{LP}) - Q^{MLP}(\boldsymbol{y}_n^{LP}, u_n^{LP})
\end{align*}
where $(a)$ is true because $(\boldsymbol{x}_n^{LP}, \eta_n^{LP}, \boldsymbol{y}_n^{LP}, u_n^{LP})$ provides an optimal solution to the LP relaxation problem, $(b)$ is true because of Proposition \ref{prop:alg}, and $(c)$ is true because $(\boldsymbol{y}_n^{1}, u_n^{1})$ is an optimal solution to Step \ref{alg:step5} in iteration 0 while $(\boldsymbol{y}_n^{LP}, u_n^{LP})$ is feasible. Here, $Q^M(\boldsymbol{y}_n^{LP}, u_n^{LP})$ is defined as the optimal objective value to the subproblem $\mbox{{\bf SP-RMS}}(\boldsymbol{y}_n^{LP}, u_n^{LP})$ \eqref{eq:ms-substructure}, and $Q^{MLP}(\boldsymbol{y}_n^{LP}, u_n^{LP})$ stands for the optimal objective value of its LP relaxation.

Based on Proposition \ref{prop:substructure}, we have a closed form for $Q^M(\boldsymbol{y}_n^{LP}, u_n^{LP})$ and we can apply the same analysis to obtain a closed form for $Q^{MLP}(\boldsymbol{y}_n^{LP}, u_n^{LP})$, where we replace all $\lceil \boldsymbol{B}_{t_n} \boldsymbol{y}_n^{LP}\rceil$ with $\boldsymbol{B}_{t_n} \boldsymbol{y}_n^{LP}$. As a result, 
\begin{align*}
    & Q^M(\boldsymbol{y}_n^{LP}, u_n^{LP}) - Q^{MLP}(\boldsymbol{y}_n^{LP}, u_n^{LP})\\
= & \sum_{n\in \mathcal{T}} p_n\left(\boldsymbol{\tilde{f}}_n^{\mathsf T}\max_{m\in\mathcal{P}(n)}\lceil \boldsymbol{B}_{t_m} \boldsymbol{y}_m^{LP}\rceil+ \tilde{\lambda}_n \max_{m\in \mathcal{C}(n)}\left\lbrace \boldsymbol{f}^{\mathsf T}_{t_m}\max_{l\in\mathcal{P}(m)}\lceil \boldsymbol{B}_{t_l} \boldsymbol{y}_l^{LP}\rceil+ \boldsymbol{c}_{t_m}^{\mathsf T} \boldsymbol{y}_m^{LP}-u_m^{LP}\right\rbrace\right) - \\
&\sum_{n\in \mathcal{T}} p_n\left(\boldsymbol{\tilde{f}}_n^{\mathsf T}\max_{m\in\mathcal{P}(n)} \boldsymbol{B}_{t_m} \boldsymbol{y}_m^{LP}+ \tilde{\lambda}_n \max_{m\in \mathcal{C}(n)}\left\lbrace \boldsymbol{f}^{\mathsf T}_{t_m}\max_{l\in\mathcal{P}(m)} \boldsymbol{B}_{t_l} \boldsymbol{y}_l^{LP}+ \boldsymbol{c}_{t_m}^{\mathsf T} \boldsymbol{y}_m^{LP}-u_m^{LP}\right\rbrace\right)\\
= & \sum_{n\in \mathcal{T}} p_n\boldsymbol{\tilde{f}}_n^{\mathsf T}\left(\max_{m\in\mathcal{P}(n)} \lceil \boldsymbol{B}_{t_m} \boldsymbol{y}_m^{LP}\rceil - \max_{m\in\mathcal{P}(n)} \boldsymbol{B}_{t_m} \boldsymbol{y}_m^{LP}\right) + \\
& \sum_{n\in \mathcal{T}} p_n \tilde{\lambda}_n\left(\max_{m\in \mathcal{C}(n)}\left\lbrace \boldsymbol{f}^{\mathsf T}_{t_m}\max_{l\in\mathcal{P}(m)}\lceil \boldsymbol{B}_{t_l} \boldsymbol{y}_l^{LP}\rceil+ \boldsymbol{c}_{t_m}^{\mathsf T} \boldsymbol{y}_m^{LP}-u_m^{LP}\right\rbrace - \max_{m\in \mathcal{C}(n)}\left\lbrace \boldsymbol{f}^{\mathsf T}_{t_m}\max_{l\in\mathcal{P}(m)} \boldsymbol{B}_{t_l} \boldsymbol{y}_l^{LP}+ \boldsymbol{c}_{t_m}^{\mathsf T} \boldsymbol{y}_m^{LP}-u_m^{LP}\right\rbrace \right)\\
\overset{(a)}{\le} & \sum_{n\in\mathcal{T}}p_n\boldsymbol{\tilde{f}}_n^{\mathsf T}\boldsymbol{1} + \sum_{n\in\mathcal{T}}p_n\tilde{\lambda}_n\boldsymbol{f}^{\mathsf T}_{t_n+1} \boldsymbol{1}\\
= & {\sum_{i=1}^Mf_{1i} + \sum_{t=2}^T(1-\lambda_t)\sum_{i=1}^Mf_{ti} + \sum_{t=2}^T\lambda_t(\sum_{i=1}^Mf_{ti})}\\
= & {\sum_{t=1}^T\sum_{i=1}^M f_{ti}}
\end{align*}
where $(a)$ is true because $\max_{m\in\mathcal{P}(n)}\lceil \boldsymbol{B}_{t_m} \boldsymbol{y}_m^{LP}\rceil - \max_{m\in\mathcal{P}(n)} \boldsymbol{B}_{t_m} \boldsymbol{y}_m^{LP}\le \boldsymbol{1}$ and by letting $m^*(n) = \arg\max_{m\in \mathcal{C}(n)}\left\lbrace \boldsymbol{f}^{\mathsf T}_{t_m}\max_{l\in\mathcal{P}(m)}\lceil \boldsymbol{B}_{t_m} \boldsymbol{y}_l^{LP}\rceil+ \boldsymbol{c}_{t_m}^{\mathsf T} \boldsymbol{y}_m^{LP}-u_m^{LP}\right\rbrace$, we have 
\begin{align*}
&\max_{m\in \mathcal{C}(n)}\left\lbrace \boldsymbol{f}^{\mathsf T}_{t_m}\max_{l\in\mathcal{P}(m)}\lceil \boldsymbol{B}_{t_l} \boldsymbol{y}_l^{LP}\rceil+ \boldsymbol{c}_{t_m}^{\mathsf T} \boldsymbol{y}_m^{LP}-u_m^{LP}\right\rbrace - \max_{m\in \mathcal{C}(n)}\left\lbrace \boldsymbol{f}^{\mathsf T}_{t_m}\max_{l\in\mathcal{P}(m)} \boldsymbol{B}_{t_l} \boldsymbol{y}_l^{LP}+ \boldsymbol{c}_{t_m}^{\mathsf T} \boldsymbol{y}_m^{LP}-u_m^{LP}\right\rbrace \\
\le 
&\left( \boldsymbol{f}^{\mathsf T}_{t_{n}+1}\max_{l\in\mathcal{P}(m^*(n))}\lceil \boldsymbol{B}_{t_l} \boldsymbol{y}_l^{LP}\rceil+ \boldsymbol{c}_{m^*(n)}^{\mathsf T} \boldsymbol{y}_{m^*(n)}^{LP}-u_{m^*(n)}^{LP}\right) - \left( \boldsymbol{f}^{\mathsf T}_{t_{n}+1}\max_{l\in\mathcal{P}(m^*(n))} \boldsymbol{B}_{t_l} \boldsymbol{y}_l^{LP}+ \boldsymbol{c}_{m^*(n)}^{\mathsf T} \boldsymbol{y}_{m^*(n)}^{LP}-u_{m^*(n)}^{LP}\right)\\
= & \boldsymbol{f}^{\mathsf T}_{t_n+1} \left(\max_{l\in\mathcal{P}(m^*(n))}\lceil \boldsymbol{B}_{t_l} \boldsymbol{y}_l^{LP}\rceil - \max_{l\in\mathcal{P}(m^*(n))} \boldsymbol{B}_{t_l} \boldsymbol{y}_l^{LP}\right)\\
\le & \boldsymbol{f}^{\mathsf T}_{t_n+1} \boldsymbol{1}.
\end{align*} 
This completes the proof. $\hfill \Box$


\proof{Proof of Theorem \ref{thm:ratio}}
Based on Proposition \ref{prop:ratio}, we have
\begin{align*}
    \frac{z_R^{MS}(\boldsymbol{x}_n^{H}, \eta_n^{H}, \boldsymbol{y}_n^{H}, u_n^{H})}{z_R^{MS}(\boldsymbol{x}_n^{*}, \eta_n^{*}, \boldsymbol{y}_n^{*}, u_n^{*})} 
    \le& \frac{z_R^{MS}(\boldsymbol{x}_n^{*}, \eta_n^{*}, \boldsymbol{y}_n^{*}, u_n^{*}) + {\sum_{t=1}^T\sum_{i=1}^M f_{ti}}}{z_R^{MS}(\boldsymbol{x}_n^{*}, \eta_n^{*}, \boldsymbol{y}_n^{*}, u_n^{*})} \\
    = &1+ \frac{{\sum_{t=1}^T\sum_{i=1}^M f_{ti}}}{z_R^{MS}(\boldsymbol{x}_n^{*}, \eta_n^{*}, \boldsymbol{y}_n^{*}, u_n^{*})}\\
    \overset{(a)}{\le} &1+ {\frac{M\sum_{t=1}^Tf_{t,\rm max}}{M_{\rm min}\sum_{t=1}^Tf_{t,\rm min}+\sum_{t=1}^Tc_{t,\rm min}\min_{n\in\mathcal{T}_t}\{\sum_{j=1}^Nd_{n,j}\}}}.
\end{align*} 
Next, we show that $(a)$ is true (i.e., $z_R^{MS}(\boldsymbol{x}_n^{*}, \eta_n^{*}, \boldsymbol{y}_n^{*}, u_n^{*})\ge M_{\rm min}\sum_{t=1}^Tf_{t,\rm min}{+\sum_{t=1}^Tc_{t,\rm min}\min_{n\in\mathcal{T}_t}\{\sum_{j=1}^Nd_{n,j}\}}$). 
As mentioned before, the optimal solution of $\eta_n$ is attained at $\eta^*_n = {\rm VaR}_{\alpha_{t_n+1}}[\boldsymbol{{f}}^{\mathsf T}_{t_m}\sum_{l\in \mathcal{P}(m)}\boldsymbol{x}_l^*+\boldsymbol{{c}}_m^{\mathsf T}\boldsymbol{y}_m^*: m\in\mathcal{C}(n)]$ by the definition of $\rm CVaR$ \eqref{eq:cvar}.
Our analysis then follows 
\begin{align*}
    z^{MS}_R(\boldsymbol{x}_n^{*}, \eta_n^{*}, \boldsymbol{y}_n^{*}, u_n^{*})=&\sum_{n\in \mathcal{T}} p_n(\boldsymbol{\tilde{f}}_n^{\mathsf T}\sum_{m\in \mathcal{P}(n)}\boldsymbol{x}_m^*+\boldsymbol{\tilde{c}}_n^{\mathsf T}\boldsymbol{y}_n^* + \tilde{\lambda}_n \eta_n^* + \tilde{\alpha}_n u_n^*)\\
    \ge& \sum_{n\in \mathcal{T}} p_n(\boldsymbol{\tilde{f}}_n^{\mathsf T}\sum_{m\in \mathcal{P}(n)}\boldsymbol{x}_m^*+\boldsymbol{\tilde{c}}_{n}^{\mathsf T}\boldsymbol{y}_n^* + \tilde{\lambda}_n \eta_n^* )\\
    =& \boldsymbol{{f}}^{\mathsf T}_{t_1}\boldsymbol{x}_1^* + \boldsymbol{{c}}_{t_1}^{\mathsf T}\boldsymbol{y}_1^* + \sum_{n\not=1}(1-\lambda_{t_n}) p_n(\boldsymbol{{f}}^{\mathsf T}_{t_n}\sum_{m\in \mathcal{P}(n)}\boldsymbol{x}_m^* + \boldsymbol{{c}}_{t_n}^{\mathsf T}\boldsymbol{y}_n^*) \\
    &+ \sum_{n\not\in\mathcal{L}} \lambda_{t_n+1}p_n {\rm VaR}_{\alpha_{t_n+1}}[\boldsymbol{{f}}^{\mathsf T}_{t_m}\sum_{l\in \mathcal{P}(m)}\boldsymbol{x}_l^*+\boldsymbol{{c}}_{t_m}^{\mathsf T}\boldsymbol{y}_m^*: m\in\mathcal{C}(n)]\\
    \overset{(a)}{\ge} & \boldsymbol{{f}}^{\mathsf T}_{t_1}\boldsymbol{x}_1^* + \boldsymbol{{c}}_{t_1}^{\mathsf T}\boldsymbol{y}_1^* + \sum_{t=2}^T\min_{n\in\mathcal{T}_t}\{\boldsymbol{{f}}^{\mathsf T}_{t_n}\sum_{m\in \mathcal{P}(n)}\boldsymbol{x}_m^* + \boldsymbol{{c}}_{t_n}^{\mathsf T}\boldsymbol{y}_n^*\}\sum_{n\in\mathcal{T}_t}(1-\lambda_{t_n})p_n \\
    & + \sum_{t=1}^{T-1}\min_{n\in\mathcal{T}_{t+1}}\{\boldsymbol{{f}}^{\mathsf T}_{t_n}\sum_{m\in \mathcal{P}(n)}\boldsymbol{x}_m^* + \boldsymbol{{c}}_{t_n}^{\mathsf T}\boldsymbol{y}_n^*\}\sum_{n\in\mathcal{T}_t}\lambda_{t_n+1}p_n \\
    \overset{(b)}{=} & \boldsymbol{{f}}^{\mathsf T}_{t_1}\boldsymbol{x}_1^* + \boldsymbol{{c}}_{t_1}^{\mathsf T}\boldsymbol{y}_1^* + \sum_{t=2}^T \min_{n\in\mathcal{T}_t}\{\boldsymbol{{f}}^{\mathsf T}_{t_n}\sum_{m\in \mathcal{P}(n)}\boldsymbol{x}_m^* + \boldsymbol{{c}}_{t_n}^{\mathsf T}\boldsymbol{y}_n^*\}\\
    \ge & \sum_{t=1}^T \min_{n\in\mathcal{T}_t}\{\boldsymbol{{f}}^{\mathsf T}_{t_n}\sum_{m\in \mathcal{P}(n)}\boldsymbol{x}_m^*\} {+\sum_{t=1}^T\min_{n\in\mathcal{T}_t}\{\boldsymbol{{c}}_{t_n}^{\mathsf T}\boldsymbol{y}_n^*}\}\\
    \overset{(c)}{\ge} & \sum_{t=1}^Tf_{t, \rm min}\sum_{i=1}^M x_{1i}^*{+\sum_{t=1}^Tc_{t,\rm min}\min_{n\in\mathcal{T}_t}\{\sum_{j=1}^Nd_{n,j}\}}\\
    \overset{(d)}{\ge} & \sum_{t=1}^Tf_{t,\rm min}M_{\rm min}{+\sum_{t=1}^Tc_{t,\rm min}\min_{n\in\mathcal{T}_t}\{\sum_{j=1}^Nd_{n,j}\}},
\end{align*}
where $(a)$ is true because $\boldsymbol{{f}}^{\mathsf T}_{t_n}\sum_{m\in \mathcal{P}(n)}\boldsymbol{x}_m^* + \boldsymbol{{c}}_{t_n}^{\mathsf T}\boldsymbol{y}_n^*\ge \min_{m\in\mathcal{T}_{t_n}}\{\boldsymbol{{f}}^{\mathsf T}_{t_m}\sum_{l\in \mathcal{P}(m)}\boldsymbol{x}_l^* + \boldsymbol{{c}}_{t_m}^{\mathsf T}\boldsymbol{y}_m^*\},\ \forall n\in\mathcal{T}$ and ${\rm VaR}_{\alpha_{t_n+1}}[\boldsymbol{{f}}^{\mathsf T}_{t_m}\sum_{l\in \mathcal{P}(m)}\boldsymbol{x}_l^*+\boldsymbol{{c}}_{t_m}^{\mathsf T}\boldsymbol{y}_m^*: m\in\mathcal{C}(n)]\ge \min_{n\in\mathcal{T}_{t_n+1}}\{\boldsymbol{{f}}^{\mathsf T}_{t_n}\sum_{m\in \mathcal{P}(n)}\boldsymbol{x}_m^* + \boldsymbol{{c}}_{t_n}^{\mathsf T}\boldsymbol{y}_n^*\},\ \forall n\not\in\mathcal{L}$; $(b)$ is true because $\sum_{n\in\mathcal{T}_t}(1-\lambda_{t_n})p_n = (1-\lambda_t)\sum_{n\in\mathcal{T}_t}p_n = 1-\lambda_t$ and $\sum_{n\in\mathcal{T}_t}\lambda_{t_n+1}p_n = \lambda_{t+1}\sum_{n\in\mathcal{T}_t}p_n = \lambda_{t+1}$; $(c)$ is true because $\sum_{m\in \mathcal{P}(n)}\boldsymbol{x}_m^*\ge \boldsymbol{x}_1^*$ and {$\boldsymbol{{c}}_{t_n}^{\mathsf T}\boldsymbol{y}_n^*\ge c_{t,\rm min}\sum_{i=1}^M\sum_{j=1}^N{y}_{n,ij}^*=c_{t,\rm min}\sum_{j=1}^Nd_{n,j}$}; $(d)$ is true because $h_{\rm max}\sum_{i=1}^Mx_{1i}^* \ge \sum_{i=1}^Mh_ix_{1i}^* \ge \sum_{i=1}^M\sum_{j=1}^Ny_{1ij}^* = \sum_{j=1}^Nd_{1j}$ from constraints \eqref{eq:2} and \eqref{eq:3}. This completes the proof.  $\hfill \Box$
\endproof

{
\proof{Proof of Corollary \ref{cor:ratio}}
First of all, $\frac{z_R^{MS}(\boldsymbol{x}_n^{H}, \eta_n^{H}, \boldsymbol{y}_n^{H}, u_n^{H})}{z_R^{MS}(\boldsymbol{x}_n^{*}, \eta_n^{*}, \boldsymbol{y}_n^{*}, u_n^{*})}\ge 1$. According to Theorem \ref{thm:ratio}, we have 
\begin{align*}
    \frac{z_R^{MS}(\boldsymbol{x}_n^{H}, \eta_n^{H}, \boldsymbol{y}_n^{H}, u_n^{H})}{z_R^{MS}(\boldsymbol{x}_n^{*}, \eta_n^{*}, \boldsymbol{y}_n^{*}, u_n^{*})} 
{\le} &1+ {\frac{M\sum_{t=1}^Tf_{t,\rm max} }{M_{\rm min}\sum_{t=1}^Tf_{t,\rm min}+\sum_{t=1}^Tc_{t,\rm min}\min_{n\in\mathcal{T}_t}\{\sum_{j=1}^Nd_{n,j}\}}}\\
 =&1+\frac{\sum_{t=1}^TO(1)}{\sum_{t=1}^TO(1)+\sum_{t=1}^TO(t)}\\
 =&1+ \frac{O(T)}{O(T^2)}\to 1\ (\text{when }T\to\infty)
\end{align*} 
$\hfill\Box$
\endproof

\proof{Proof of Lemma \ref{lemma:prior}}
Observing that the right-hand sides of constraints \eqref{eq:two-prior-eta} and \eqref{eq:multi-prior-eta} are different, we make a variable transformation in \eqref{model:two-prior-original} --- $\tilde{\eta}_{a(n)}:=\eta_{n}+\sum_{i\not=i^{\prime}}s_{n,ii^{\prime}}, \ \forall n\not=1,\ n\not\in\mathcal{L},\ \tilde{\eta}_{a(n)}:=\eta_{n},\ \forall n\in\mathcal{L}$. Because both $\eta$ and $\boldsymbol{s}$ are first-stage variables and are identical for all nodes in the same stage, this transformation is without of ambiguity. With that, the two-stage constraints \eqref{eq:two-prior-twostage-eta} for the new variable $\tilde{\eta}_n$ are satisfied and constraints \eqref{eq:two-prior-eta} become \eqref{eq:multi-prior-etaL}--\eqref{eq:multi-prior-eta}.
Moreover, the terms in the objective function \eqref{eq:two-prior-original-obj} that are related to $\eta_n$ and $\boldsymbol{s}_n$ can be transformed to 
\begin{align*}
    &\sum_{n\in \mathcal{T}\setminus\{1\}} p_n {\lambda}_{t_n}\eta_n +\sum_{n\in \mathcal{T}\setminus\mathcal{L}} p_n \sum_{i\not=i^{\prime}}s_{n,ii^{\prime}}\\
    {=}&\sum_{n\in \mathcal{T}\setminus(\{1\}\cup\mathcal{L})} p_n{\lambda}_{t_n}(\tilde{\eta}_{a(n)}-\sum_{i\not=i^{\prime}}s_{n,ii^{\prime}}) +\sum_{n\in\mathcal{L}}p_n\lambda_{t_n}\tilde{\eta}_{a(n)}+\sum_{n\in \mathcal{T}\setminus\mathcal{L}} p_n\sum_{i\not=i^{\prime}}s_{n,ii^{\prime}}\\
    \overset{(a)}{=}&\sum_{n\in\mathcal{T}\setminus\mathcal{L}}p_n\lambda_{t_n+1}\tilde{\eta}_n+\sum_{i\not=i^{\prime}}s_{1,ii^{\prime}}+\sum_{n\in\mathcal{T}\setminus(\{1\}\cup\mathcal{L})}p_n(1-\lambda_{t_n})\sum_{i\not=i^{\prime}}s_{n,ii^{\prime}}
\end{align*}
where the equality $(a)$ is true because $\sum_{m\in\mathcal{C}(n)}p_m\lambda_{t_m} = p_n\lambda_{t_n+1}$. Thus, the objective function \eqref{eq:two-prior-original-obj} becomes the objective function \eqref{eq:multi-prior-obj}. This completes the proof. $\hfill\Box$
\endproof

\proof{Proof of Proposition \ref{prop:substructure-prior}}
First of all, $\boldsymbol{s}^{MS}$ is feasible for Model \eqref{eq:ms-substructure-prior} by construction. According to constraints \eqref{eq:two-prior-equal} and \eqref{eq:two-prior-sx}, for any feasible solution we have $s_{1,ii^{\prime}}+s_{1,i^{\prime}i}\ge 1, \ s_{n,ii^{\prime}}+s_{n,i^{\prime}i} \ge \max\{0,1 - (\sum_{m\in \mathcal{P}(n)\setminus\{1\}}x^{*}_{m,i} + \sum_{m\in \mathcal{P}(n)\setminus\{1\}}x^{*}_{m,i^{\prime}})\},\ \forall n\not=1,\ n\not\in\mathcal{L}$. Because we are minimizing $\sum_{n\in\mathcal{T}}p_n\tilde{1}_n\sum_{i\not=i^{\prime}}s_{n,ii^{\prime}}$ and $s^{MS}_{1,ii^{\prime}}+s^{MS}_{1,i^{\prime}i}= 1, \ s^{MS}_{n,ii^{\prime}}+s^{MS}_{n,i^{\prime}i}= \max\{0,1 - (\sum_{m\in \mathcal{P}(n)\setminus\{1\}}x^{*}_{m,i} + \sum_{m\in \mathcal{P}(n)\setminus\{1\}}x^{*}_{m,i^{\prime}})\},\ \forall n\not=1,\ n\not\in\mathcal{L}$, we conclude that $\boldsymbol{s}^{MS}$ is optimal. Due to constraints \eqref{eq:multi-prior-etaL} and \eqref{eq:multi-prior-eta}, $\boldsymbol{\eta}^{MS}$ is optimal.

As for the risk-averse two-stage model \eqref{eq:ts-substructure-prior}, because $\boldsymbol{s}_{n}$ needs to be identical for all the nodes $n\in\mathcal{T}\setminus\mathcal{L}$, we have $s_{n,ii^{\prime}}=s_{1,ii^{\prime}}\le \min_{m\in\mathcal{T}\setminus\{1\}}\left\{\sum_{l\in \mathcal{P}(m)\setminus\{1\}}x^*_{l,i} - \sum_{l\in \mathcal{P}(m)\setminus\{1\}}x^*_{l,i^{\prime}} + 1\right\},\ \forall n\in\mathcal{T}\setminus\mathcal{L}$ and $s_{n,ii^{\prime}}+s_{n,i^{\prime}i} = s_{1,ii^{\prime}}+s_{1,i^{\prime}i}\ge 1 $. By construction, $\boldsymbol{s}^{TS}$ is optimal. According to constraints \eqref{eq:multi-prior-etaL}, \eqref{eq:multi-prior-eta} and \eqref{eq:two-prior-twostage-eta}, we conclude that $\boldsymbol{\eta}^{TS}$ is optimal. 
$\hfill \Box$
\endproof

\proof{Proof of Theorem \ref{thm:risk-averse-bound-prior}}
First of all, ${\rm VMS_P^{LB}}\ge 0$ and thus provides a nontrivial lower bound. We have $
    z^{TS}_P = \sum_{n\in \mathcal{T}} p_n\left(\tilde{\boldsymbol{c}}_n^{\mathsf T}\boldsymbol{y}^{*}_n +  \tilde{\alpha}_n u^{*}_n\right) + Q_P^T(\boldsymbol{x}_n^{*},\boldsymbol{y}_n^{*},\boldsymbol{u}_n^{*})
$. 
Since $(\boldsymbol{x}_n^{*}, \boldsymbol{y}_n^{*}, u_n^{*})_{n\not=1}$ are feasible solutions for the multistage model \eqref{model:multi-prior}, we have
\begin{align}
    z_P^{MS}\le \sum_{n\in \mathcal{T}} p_n\left(\tilde{\boldsymbol{c}}_n^{\mathsf T}\boldsymbol{y}^{*}_n +  \tilde{\alpha}_n u^{*}_n\right) + Q_P^M(\boldsymbol{x}_n^{*},\boldsymbol{y}_n^{*},\boldsymbol{u}_n^{*})\label{eq:multi-prior-feasible}
\end{align}
Combining the above steps, we obtain
\begin{align*}
    {\rm VMS_P} =& z_P^{TS} - z_P^{MS}\\
    \ge & Q_P^T(\boldsymbol{x}_n^{*},\boldsymbol{y}_n^{*},\boldsymbol{u}_n^{*}) - Q_P^M(\boldsymbol{x}_n^{*},\boldsymbol{y}_n^{*},\boldsymbol{u}_n^{*})\\
    = & \sum_{n\in \mathcal{T}\setminus\mathcal{L}} p_n{\lambda}_{t_n+1}(\eta_n^{TS}-\eta_n^{MS}) + \sum_{n\in \mathcal{T}\setminus(\{1\}\cup\mathcal{L})} p_n (1-\lambda_{t_n})(\sum_{i\not=i^{\prime}}s^{TS}_{n,ii^{\prime}}-\sum_{i\not=i^{\prime}}s^{MS}_{n,ii^{\prime}})\\
    \overset{(a)}{\ge}&\sum_{n\in \mathcal{T}\setminus\mathcal{L}} p_n{\lambda}_{t_n+1}(\eta_n^{TS}-\eta_n^{MS}) \\
    \overset{(b)}{\ge} & \sum_{n\in \mathcal{T}\setminus\mathcal{L}} p_n{\lambda}_{t_n+1} \left(\max_{m\in\mathcal{T}_{t_n+1}}\{\boldsymbol{c}_{t_m}^{\mathsf T}\boldsymbol{y}_m^{*}-u_m^{*}\}-\max_{m\in\mathcal{C}(n)}\{\boldsymbol{c}_{t_m}^{\mathsf T}\boldsymbol{y}_m^{*}-u_m^{*}\}\right).
\end{align*}
where $(a)$ is true because $\sum_{i\not=i^{\prime}}s^{TS}_{n,ii^{\prime}}=\frac{M(M-1)}{2}\ge \sum_{i\not=i^{\prime}}s^{MS}_{n,ii^{\prime}}$, $(b)$ is true because when $n\in\mathcal{T}_t,\ t=0,\ldots, T-2$, $\eta_n^{TS}=\max_{m\in\mathcal{T}_{t_n+1}}\{\boldsymbol{c}_{t_m}^{\mathsf T} \boldsymbol{y}^{*}_m + \frac{M(M-1)}{2}-u_m^{*}\},\ \eta_n^{MS}\le\max_{m\in\mathcal{C}(n)}\{\boldsymbol{c}_{t_m}^{\mathsf T} \boldsymbol{y}^{*}_m + \frac{M(M-1)}{2}-u_m^{*}\}$, and when $n\in\mathcal{T}_{T-1}$, the inequality holds by construction.
$\hfill \Box$
\endproof

\proof{Proof of Theorem \ref{thm:prior-cuts}}
Let $({\boldsymbol{s}}^*_n, {\boldsymbol{x}}^*_n, {\eta}^*_n, {\boldsymbol{y}}^*_n, {u}_n^*)_{n\in\mathcal{T}}$ be an optimal solution to model \eqref{model:multi-prior}. Let $S(\boldsymbol{x}_n)$ denote the set of open facilities at solution $\boldsymbol{x}_n$, i.e., $S(\boldsymbol{x}_n)=\{i\in[M]|\sum_{m\in\mathcal{P}(n)\setminus\{1\}}x_{m,i}=1\}$. Given two nodes $n$ and $n^{\prime}$, having $S(\boldsymbol{x}^*_n)\supseteq S(\boldsymbol{x}^*_{n^{\prime}})$ implies that $\sum_{m\in\mathcal{P}(n)\setminus\{1\}}{x}^*_{m,i}\ge \sum_{m\in\mathcal{P}(n^{\prime})\setminus\{1\}}{x}^*_{m,i},\ \forall i=1,\ldots, M$. If for all $n\not=1, n^{\prime}\in\mathcal{C}(a(n))$ such that $F_n\ge F_{n^{\prime}}$ we have $S(\boldsymbol{x}^*_n)\supseteq S(\boldsymbol{x}^*_{n^{\prime}})$, the proof is complete. Now suppose for a pair of nodes $(\bar{n},\bar{n}^{\prime})$ such that $\bar{n}\not=1, \bar{n}^{\prime}\in\mathcal{C}(a(\bar{n}))$ with $F_{\bar{n}}\ge F_{\bar{n}^{\prime}}$, we have $S(\boldsymbol{x}^*_{\bar{n}})\subset S(\boldsymbol{x}^*_{\bar{n}^{\prime}})$. Then, we construct a new solution $(\bar{\boldsymbol{s}}_n, \bar{\boldsymbol{x}}_n, \bar{\eta}_n, \bar{\boldsymbol{y}}_n, \bar{u}_n)_{n\in\mathcal{T}}$ to model \eqref{model:multi-prior} such that 
\begin{enumerate}[label=(\roman*)]
	\item $S(\bar{\boldsymbol{x}}_n) = S(\boldsymbol{x}^*_n),\ \forall n\in\mathcal{T}\setminus\{\mathcal{T}(\bar{n})\}$, and $S(\bar{\boldsymbol{x}}_n)=S(\boldsymbol{x}^*_n)\cup S(\boldsymbol{x}^*_{\bar{n}^{\prime}}),\ \forall n\in \mathcal{T}(\bar{n})$;
	\item $(\bar{\boldsymbol{s}}_n, \bar{\boldsymbol{x}}_n, \bar{\eta}_n, \bar{\boldsymbol{y}}_n, \bar{u}_n)_{n\in\mathcal{T}}$ is feasible for model \eqref{model:multi-prior};
	\item  $(\bar{\boldsymbol{s}}_n, \bar{\boldsymbol{x}}_n, \bar{\eta}_n, \bar{\boldsymbol{y}}_n, \bar{u}_n)_{n\in\mathcal{T}}$ has no worse objective function value than $({\boldsymbol{s}}^*_n, {\boldsymbol{x}}^*_n, {\eta}^*_n, {\boldsymbol{y}}^*_n, {u}_n^*)_{n\in\mathcal{T}}$.
\end{enumerate}
Recall that $\mathcal{T}(n)$ denotes the subtree rooted at node $n$ (see Figure \ref{fig:scenario-tree-prior}(b)). The above three conditions ensure that $(\bar{\boldsymbol{s}}_n, \bar{\boldsymbol{x}}_n, \bar{\eta}_n, \bar{\boldsymbol{y}}_n, \bar{u}_n)_{n\in\mathcal{T}}$ is also optimal to model \eqref{model:multi-prior} and $S(\bar{\boldsymbol{x}}_{\bar{n}})\supseteq S(\bar{\boldsymbol{x}}_{\bar{n}^{\prime}})$. 
Repeating this argument for all $\bar{n}\not=1, \bar{n}^{\prime}\in\mathcal{C}(a(\bar{n}))$ with $F_{\bar{n}}\ge F_{\bar{n}^{\prime}}$ such that $S(\boldsymbol{x}^*_{\bar{n}})\subset S(\boldsymbol{x}^*_{\bar{n}^{\prime}})$ proves the theorem.

Now consider a specific new solution $(\bar{\boldsymbol{s}}_n, \bar{\boldsymbol{x}}_n, \bar{\eta}_n, \bar{\boldsymbol{y}}_n, \bar{u}_n)_{n\in\mathcal{T}}$ defined as follows:
\begin{itemize}
	\item $\bar{\boldsymbol{x}}_n = \boldsymbol{x}^*_n,\ \forall n\in\mathcal{T}\setminus\mathcal{T}(\bar{n}),\ \bar{\boldsymbol{x}}_{\bar{n}}=\boldsymbol{x}^*_{\bar{n}^{\prime}}$, and for all $n\in \mathcal{T}(\bar{n})\setminus\{\bar{n}\}$ we have
	\begin{align*}
	\bar{x}_{n,i}=\begin{cases}
	0, \ \text{if } i\in S(\bar{\boldsymbol{x}}_{\bar{n}})\setminus S({\boldsymbol{x}}^*_{\bar{n}})\\
	x_{n,i}^*,\ \text{otherwise}
	\end{cases}
	\end{align*}
	\item $\bar{\boldsymbol{y}}_n=\boldsymbol{y}_n^*,\ \forall n\in \mathcal{T}\setminus\{\bar{n}\}$, and $\bar{\boldsymbol{y}}_{\bar{n}}=\boldsymbol{Y}_{\bar{n}}(\boldsymbol{x}^*_{\bar{n}^{\prime}})$ where $\boldsymbol{Y}_{n}(\boldsymbol{x})$ denotes a minimizer of  the problem $\min_{\boldsymbol{y}}\{\tilde{\boldsymbol{c}}_n^{\mathsf T}\boldsymbol{y}:(\boldsymbol{x},\boldsymbol{y})\in X_n\}$.
	\item $\bar{\boldsymbol{s}}_n = \boldsymbol{s}^*_n,\ \forall n\in\mathcal{T}\setminus \mathcal{T}(\bar{n})$, and for all $n\in\mathcal{T}(\bar{n})\setminus \mathcal{L}$ we have
	\begin{align*}
	\bar{s}_{n,ii^{\prime}}=\begin{cases}
	0,\ \text{if } i \in S(\bar{\boldsymbol{x}}_{\bar{n}})\setminus S({\boldsymbol{x}}^*_{\bar{n}}) \text{ or } i^{\prime}\in S(\bar{\boldsymbol{x}}_{\bar{n}})\setminus S({\boldsymbol{x}}^*_{\bar{n}})\\
	s^*_{n,ii^{\prime}},\ \text{otherwise}
	\end{cases}
	\end{align*}
	\item $\bar{\eta}_n = \eta^*_n,\ \forall n\in\mathcal{T}\setminus(\{a(\bar{n})\}\cup\mathcal{T}(\bar{n}))$, $\bar{\eta}_{n}=\text{VaR}_{\alpha_{t_{n}+1}}[\boldsymbol{c}_{t_m}^{\mathsf T}\bar{\boldsymbol{y}}_m+\sum_{i\not=i^{\prime}}\bar{s}_{m, ii^{\prime}}: m\in \mathcal{C}(n)],\ n\in \{a(\bar{n})\}\cup\mathcal{T}(\bar{n})\setminus (\mathcal{T}_{T-1}\cup\mathcal{L})$, and for all $n\in\mathcal{T}_{T-1}$, we have $\bar{\eta}_{n}=\text{VaR}_{\alpha_{t_{n}+1}}[\boldsymbol{c}_{t_m}^{\mathsf T}\bar{\boldsymbol{y}}_m: m\in \mathcal{C}(n)]$
	\item $\bar{u}_n=u^*_n,\ \forall n\in\mathcal{T}\setminus(\{\bar{n}^{\prime}\}\cup\mathcal{T}(\bar{n}))$, $\bar{u}_n = [\boldsymbol{c}_{t_n}^{\mathsf T}\bar{\boldsymbol{y}}_n+\sum_{i\not=i^{\prime}}\bar{s}_{n, ii^{\prime}} - \bar{\eta}_{a(n)}]_+,\ \forall n\in\{\bar{n}^{\prime}\}\cup\mathcal{T}(\bar{n})\setminus\mathcal{L}$, and $\bar{u}_n = [\boldsymbol{c}_{t_n}^{\mathsf T}\bar{\boldsymbol{y}}_n - \bar{\eta}_{a(n)}]_+,\ \forall n\in\mathcal{L}$.
\end{itemize}
Next let us check that this new solution satisfies the conditions (i)-(iii):
\begin{enumerate}[label=(\roman*)]
	\item This holds by construction;
	\item Feasibility of $(\bar{\boldsymbol{s}}_n, \bar{\boldsymbol{x}}_n, \bar{\eta}_n, \bar{\boldsymbol{y}}_n, \bar{u}_n)_{n\in\mathcal{T}}$:
	\begin{itemize}
		\item Constraints \eqref{eq:two-prior-equal}--\eqref{eq:two-prior-s} are satisfied at nodes $n \in\mathcal{T}\setminus\mathcal{T}(\bar{n})$ automatically because we only modified $\boldsymbol{s}^*$ and $\boldsymbol{x}^*$ at nodes $n\in\mathcal{T}(\bar{n})$. For all node $n\in\mathcal{T}(\bar{n})$, if $i \in S(\bar{\boldsymbol{x}}_{\bar{n}})\setminus S({\boldsymbol{x}}^*_{\bar{n}}) \text{ or } i^{\prime}\in S(\bar{\boldsymbol{x}}_{\bar{n}})\setminus S({\boldsymbol{x}}^*_{\bar{n}})$, we have $\bar{s}_{n,ii^{\prime}}=0$ and in this case, $\sum_{m\in\mathcal{P}(n)\setminus\{1\}}\bar{x}_{m,i}=1$ or $\sum_{m\in\mathcal{P}(n)\setminus\{1\}}\bar{x}_{m,i^{\prime}}=1$. Otherwise, $\bar{x}_{n,i}=x^*_{n,i},\ \bar{s}_{n,ii^{\prime}}=s^*_{n,ii^{\prime}}$. Thus, constraints \eqref{eq:two-prior-sx} are satisfied at nodes $n\in\mathcal{T}(\bar{n})$. As for constraints \eqref{eq:two-prior-xs}, it holds at node $\bar{n}$ because $\bar{n}$ shares the same priority list $\bar{\boldsymbol{s}}_{a(\bar{n})}$ and facility location decisions $\bar{\boldsymbol{x}}_{\bar{n}}$ with node $\bar{n}^{\prime}$. For all nodes $n\in \mathcal{T}(\bar{n})\setminus\{\bar{n}\}$, if $i \in S(\bar{\boldsymbol{x}}_{\bar{n}})\setminus S({\boldsymbol{x}}^*_{\bar{n}}) \text{ or } i^{\prime}\in S(\bar{\boldsymbol{x}}_{\bar{n}})\setminus S({\boldsymbol{x}}^*_{\bar{n}})$, we have $\bar{s}_{a(n),ii^{\prime}}=0$ and thus constraint \eqref{eq:two-prior-xs} is satisfied. Otherwise, $\bar{x}_{n,i}=x^*_{n,i},\ \bar{s}_{a(n),ii^{\prime}}=s^*_{a(n),ii^{\prime}}$ and constraint \eqref{eq:two-prior-xs} is satisfied.
		\item Constraint \eqref{eq:two-prior-budget} is true at node $\bar{n}$ because $F_{\bar{n}}\ge F_{\bar{n}^{\prime}}$. For all other nodes, $\bar{\boldsymbol{x}}_n\le\boldsymbol{x}^*_n$ and thus constraints \eqref{eq:two-prior-budget} are satisfied.
		\item Constraint \eqref{eq:two-prior-constraint_x} is satisfied, i.e., $(\bar{\boldsymbol{x}}_n,\bar{\boldsymbol{y}}_n)\in X_n$, at node $\bar{n}$ because of the construction of $\bar{\boldsymbol{y}}_{\bar{n}}$. At all nodes $n\in\mathcal{T}(\bar{n})\setminus\{\bar{n}\}$, we have $S(\bar{\boldsymbol{x}}_n)\supseteq S(\boldsymbol{x}^*_n)$ and thus $\sum_{m\in\mathcal{P}(n)\setminus\{1\}}\bar{\boldsymbol{x}}_m\ge \sum_{m\in\mathcal{P}(n)\setminus\{1\}}{\boldsymbol{x}}^*_m$ while $\boldsymbol{y}$ stays the same. Moreover, if $i\in S(\bar{\boldsymbol{x}}_{\bar{n}})\setminus S(\boldsymbol{x}^*_{\bar{n}})$, we have $\bar{x}_{n,i}=0$ and thus $\sum_{m\in\mathcal{P}(n)\setminus\{1\}}\bar{{x}}_{m,i}\le 1$. Otherwise, $\bar{x}_{n,i}=x_{n,i}^*$ and $\sum_{m\in\mathcal{P}(n)\setminus\{1\}}\bar{{x}}_{m,i}=\sum_{m\in\mathcal{P}(n)\setminus\{1\}}{{x}}^*_{m,i}\le 1$.
		\item Finally, constraints \eqref{eq:multi-prior-etaL}--\eqref{eq:multi-prior-eta} are satisfied by construction of $\bar{\boldsymbol{\eta}}$ and $\bar{\boldsymbol{u}}$.
	\end{itemize}
	\item Optimality of $(\bar{\boldsymbol{s}}_n, \bar{\boldsymbol{x}}_n, \bar{\eta}_n, \bar{\boldsymbol{y}}_n, \bar{u}_n)_{n\in\mathcal{T}}$:\\ Since $\bar{\boldsymbol{y}}_n=\boldsymbol{y}_n^*,\ \forall n\in \mathcal{T}\setminus\{\bar{n}\}$, $\bar{\eta}_n = \eta^*_n,\ \forall n\in\mathcal{T}\setminus(\{a(\bar{n})\}\cup\mathcal{T}(\bar{n})),\ \bar{u}_n=u^*_n,\ \forall n\in\mathcal{T}\setminus(\{\bar{n}^{\prime}\}\cup\mathcal{T}(\bar{n})),\ \bar{\boldsymbol{s}}_n = \boldsymbol{s}^*_n,\ \forall n\in\mathcal{T}\setminus \mathcal{T}(\bar{n})$, we have
	\small
	\begin{align*}
	&z_P^{MS}(\bar{\boldsymbol{s}}_n, \bar{\boldsymbol{x}}_n, \bar{\eta}_n, \bar{\boldsymbol{y}}_n, \bar{u}_n) - z_P^{MS}({\boldsymbol{s}}^*_n, {\boldsymbol{x}}^*_n, {\eta}^*_n, {\boldsymbol{y}}^*_n, {u}_n^*) \\
	=& p_{\bar{n}}(\tilde{\boldsymbol{c}}_{\bar{n}}^{\mathsf T}\bar{\boldsymbol{y}}_{\bar{n}}-\tilde{\boldsymbol{c}}_{\bar{n}}^{\mathsf T}\boldsymbol{y}^*_{\bar{n}})  +\sum_{n\in \mathcal{T}(\bar{n})\setminus \mathcal{L}}p_n(1-\lambda_{t_n})(\sum_{i\not=i^{\prime}}\bar{s}_{n,ii^{\prime}}-\sum_{i\not=i^{\prime}}s_{n,ii^{\prime}}^*)\\
	&+\sum_{n\in \{a(\bar{n})\}\cup\mathcal{T}(\bar{n})\setminus \mathcal{L}}p_n\lambda_{t_n+1}(\bar{\eta}_n-\eta_n^*) + \sum_{n\in \{\bar{n}^{\prime}\}\cup\mathcal{T}(\bar{n})}p_n\frac{\lambda_{t_n}}{1-\alpha_{t_n}}(\bar{u}_n-u^*_n)\\
	\overset{(a)}{\le}& \sum_{n\in \{a(\bar{n})\}\cup\mathcal{T}(\bar{n})\setminus \mathcal{L}}p_n\lambda_{t_n+1}\left(\left(\bar{\eta}_n+\frac{1}{1-\alpha_{t_n+1}}\sum_{m\in\mathcal{C}(n)}\frac{p_m}{p_n}\bar{u}_m\right) - \left({\eta^*_n}+\frac{1}{1-\alpha_{t_n+1}}\sum_{m\in\mathcal{C}(n)}\frac{p_m}{p_n}{u}^*_m\right)\right)\\
	\overset{(b)}{=} & \sum_{n\in \{a(\bar{n})\}\cup\mathcal{T}(\bar{n})\setminus \mathcal{L}}p_n\lambda_{t_n+1}\left({\rm CVaR}_{\alpha_{t_n+1}}[\boldsymbol{c}_{t_m}^{\mathsf T}\bar{\boldsymbol{y}}_m+\sum_{i\not=i^{\prime}}\bar{s}_{m,ii^{\prime}}: m\in\mathcal{C}(n)]-{\rm CVaR}_{\alpha_{t_n+1}}[\boldsymbol{c}_{t_m}^{\mathsf T}{\boldsymbol{y}}^*_m+\sum_{i\not=i^{\prime}}{s}^*_{m,ii^{\prime}}: m\in\mathcal{C}(n)]\right)\\
	\overset{(c)}{\le} &0
	\end{align*}
	\normalsize
	where in $(a)$, $(\tilde{\boldsymbol{c}}_{\bar{n}}^{\mathsf T}\bar{\boldsymbol{y}}_{\bar{n}}-\tilde{\boldsymbol{c}}_{\bar{n}}^{\mathsf T}\boldsymbol{y}^*_{\bar{n}})\le 0$ because $\sum_{m\in\mathcal{P}(\bar{n})\setminus\{1\}}x_{m,i}^*\le \sum_{m\in\mathcal{P}(\bar{n})\setminus\{1\}}\bar{x}_{m,i}$ and thus ${\boldsymbol{y}}^*_{\bar{n}}$ is a feasible solution to the problem $\min_{\boldsymbol{y}}\{\tilde{\boldsymbol{c}}_{\bar{n}}^{\mathsf T}\boldsymbol{y}:(\bar{\boldsymbol{x}}_{\bar{n}},\boldsymbol{y})\in X_{\bar{n}}\}$ while $\bar{\boldsymbol{y}}_{\bar{n}}$ is a minimizer; $(\sum_{i\not=i^{\prime}}\bar{s}_{n,ii^{\prime}}-\sum_{i\not=i^{\prime}}s_{n,ii^{\prime}}^*)\le 0$ because $\bar{\boldsymbol{s}}$ changes some elements to 0 and thus the sum becomes smaller. Equation (b) is true because of the definition of CVaR and the constructions of $\bar{\boldsymbol{\eta}}$ and $\bar{\boldsymbol{u}}$, and $(c)$ is due to the monotonicity of CVaR where $\boldsymbol{c}_{t_n}^{\mathsf T}\bar{\boldsymbol{y}}_n+\sum_{i\not=i^{\prime}}\bar{s}_{n,ii^{\prime}}\le \boldsymbol{c}_{t_n}^{\mathsf T}{\boldsymbol{y}}^*_n+\sum_{i\not=i^{\prime}}{s}^*_{n,ii^{\prime}}$ for all $n\in \{\bar{n}^{\prime}\}\cup\mathcal{T}(\bar{n})$.
\end{enumerate}
This completes the proof. $\hfill \Box$
\endproof

\proof{Proof of Theorem \ref{thm:prior-two-stage-cuts}}
We follow the proof of Theorem \ref{thm:prior-cuts} where we define $\bar{\boldsymbol{x}}_n$ and $\bar{\boldsymbol{y}}_n$ the same as in the proof of Theorem \ref{thm:prior-cuts} and only modify $\bar{\boldsymbol{s}}_n, \bar{{\eta}}_n$ and $\bar{{u}}_n$ as follows:
\begin{itemize}
	\item $\bar{\boldsymbol{s}}_n = \boldsymbol{s}_n^*,\ \forall n\in\mathcal{T}\setminus\mathcal{L}$
	\item $\bar{\eta}_n = \eta^*_n,\ \forall n\in\mathcal{T}\setminus(\mathcal{T}_{t_{\bar{n}}-1}\cup\mathcal{L})$, and for all nodes $n\in \mathcal{T}_{t_{\bar{n}}-1}$, we set $\bar{\eta}_{n}=\text{VaR}_{\alpha_{t_{\bar{n}}}}[\boldsymbol{c}_{t_m}^{\mathsf T}\bar{\boldsymbol{y}}_m+\sum_{i\not=i^{\prime}}\bar{s}_{m, ii^{\prime}}: m\in \mathcal{T}_{t_{\bar{n}}}]$
	\item $\bar{u}_n=u^*_n,\ \forall n\in\mathcal{T}\setminus\mathcal{T}_{t_{\bar{n}}}$, and $\bar{u}_n = [\boldsymbol{c}_{t_n}^{\mathsf T}\bar{\boldsymbol{y}}_n+\sum_{i\not=i^{\prime}}\bar{s}_{n, ii^{\prime}} - \bar{\eta}_{a(n)}]_+,\ \forall n\in\mathcal{T}_{t_{\bar{n}}}$.
\end{itemize}
Next let us check that this new solution satisfies the conditions (i)-(iii):
\begin{enumerate}[label=(\roman*)]
	\item This holds by construction;
	\item Feasibility of $(\bar{\boldsymbol{s}}_n, \bar{\boldsymbol{x}}_n, \bar{\eta}_n, \bar{\boldsymbol{y}}_n, \bar{u}_n)_{n\in\mathcal{T}}$:
	\begin{itemize}
		\item Constraints \eqref{eq:two-prior-equal}--\eqref{eq:two-prior-s} are satisfied automatically because in the two-stage model, $\bar{s}_{n,ii^{\prime}}+\bar{s}_{n,i^{\prime}i}={s}^*_{n,ii^{\prime}}+{s}^*_{n,i^{\prime}i}={s}^*_{1,ii^{\prime}}+{s}^*_{1,i^{\prime}i}\ge 1$ for all $n\in\mathcal{T}\setminus\mathcal{L}$.
		As for constraints \eqref{eq:two-prior-xs}, it holds at nodes $n\in\mathcal{T}\setminus\mathcal{T}(\bar{n})$ because we only modify $\boldsymbol{x}$ at nodes $n\in\mathcal{T}(\bar{n})$. It also holds at node $\bar{n}$ because $\bar{n}$ shares the same priority list and facility location decisions with node $\bar{n}^{\prime}$. For all nodes $n\in \mathcal{T}(\bar{n})\setminus\{\bar{n}\}$, if $i \in S(\bar{\boldsymbol{x}}_{\bar{n}})\setminus S({\boldsymbol{x}}^*_{\bar{n}})$ and $\bar{s}_{a(n),ii^{\prime}}=1$, we have $\sum_{m\in\mathcal{P}(n)\setminus\{1\}}x_{m,i}=1$ and thus constraint \eqref{eq:two-prior-xs} is satisfied. Otherwise, $ \bar{s}_{a(n),ii^{\prime}}=0$ and constraint \eqref{eq:two-prior-xs} is satisfied. For all other sites $i\not\in S(\bar{\boldsymbol{x}}_{\bar{n}})\setminus S({\boldsymbol{x}}^*_{\bar{n}})$, $\bar{\boldsymbol{x}}$ and $\bar{\boldsymbol{s}}$ stay the same and thus constraint \eqref{eq:two-prior-xs} is satisfied.
		\item Constraint \eqref{eq:two-prior-budget} is true at node $\bar{n}$ because $F_{\bar{n}}\ge F_{\bar{n}^{\prime}}$. For all other nodes, $\bar{\boldsymbol{x}}_n\le\boldsymbol{x}^*_n$ and thus constraints \eqref{eq:two-prior-budget} are satisfied.
		\item Constraint \eqref{eq:two-prior-constraint_x} is satisfied, i.e., $(\boldsymbol{x}_n,\boldsymbol{y}_n)\in X_n$, because of the same reasoning as in the multistage case.
		\item Constraints \eqref{eq:multi-prior-etaL}--\eqref{eq:multi-prior-eta} are satisfied by construction of $\bar{\boldsymbol{\eta}}$ and $\bar{\boldsymbol{u}}$.
		\item Finally, constraints \eqref{eq:two-prior-twostage-s}--\eqref{eq:two-prior-twostage-eta} are satisfied because we did not modify $\boldsymbol{s}_n^*$ and $\bar{\eta}_n$ is identical for all the nodes in stage $t_{\bar{n}}-1$.
	\end{itemize}
	\item Since $\bar{\boldsymbol{y}}_n=\boldsymbol{y}_n^*,\ \forall n\in \mathcal{T}\setminus\{\bar{n}\}$, $\bar{\eta}_n = \eta^*_n,\ \forall n\in\mathcal{T}\setminus\mathcal{T}_{t_{\bar{n}}-1},\ \bar{u}_n=u^*_n,\ \forall n\in\mathcal{T}\setminus\mathcal{T}_{t_{\bar{n}}}$, we have
	\small
	\begin{align*}
	&z_P^{MS}(\bar{\boldsymbol{s}}_n, \bar{\boldsymbol{x}}_n, \bar{\eta}_n, \bar{\boldsymbol{y}}_n, \bar{u}_n) - z_P^{MS}({\boldsymbol{s}}^*_n, {\boldsymbol{x}}^*_n, {\eta}^*_n, {\boldsymbol{y}}^*_n, {u}_n^*) \\
	=& p_{\bar{n}}(\tilde{\boldsymbol{c}}_{\bar{n}}^{\mathsf T}\bar{\boldsymbol{y}}_{\bar{n}}-\tilde{\boldsymbol{c}}_{\bar{n}}^{\mathsf T}\boldsymbol{y}^*_{\bar{n}})+\sum_{n\in\mathcal{T}_{t_{\bar{n}}-1}}p_{n}\lambda_{t_{n}+1}(\bar{\eta}_{n}-\eta_{n}^*) + \sum_{n\in \mathcal{T}_{t_{\bar{n}}}}p_n\frac{\lambda_{t_n}}{1-\alpha_{t_n}}(\bar{u}_n-u^*_n)\\
	\overset{(a)}{\le}& \lambda_{t_{\bar{n}}}\left(\left(\bar{\eta}_{a(\bar{n})}+\frac{1}{1-\alpha_{t_{\bar{n}}}}\sum_{m\in\mathcal{T}_{t_{\bar{n}}}}p_m\bar{u}_m\right) - \left({\eta^*_{a(\bar{n})}}+\frac{1}{1-\alpha_{t_{\bar{n}}}}\sum_{m\in\mathcal{T}_{t_{\bar{n}}}}p_m{u}^*_m\right)\right)\\
	= & \lambda_{t_{\bar{n}}}\left({\rm CVaR}_{\alpha_{t_n+1}}[\boldsymbol{c}_{t_m}^{\mathsf T}\bar{\boldsymbol{y}}_m+\sum_{i\not=i^{\prime}}\bar{s}_{m,ii^{\prime}}: m\in\mathcal{T}_{t_{\bar{n}}}]-{\rm CVaR}_{\alpha_{t_n+1}}[\boldsymbol{c}_{t_m}^{\mathsf T}{\boldsymbol{y}}^*_m+\sum_{i\not=i^{\prime}}{s}^*_{m,ii^{\prime}}: m\in\mathcal{T}_{t_{\bar{n}}}]\right)\\
	\le &0,
	\end{align*}
	\normalsize
	where $(a)$ is true because $\sum_{n\in\mathcal{T}_{t_{\bar{n}}-1}}p_{n}\lambda_{t_{{n}}+1}\bar{\eta}_{n}=\lambda_{t_{\bar{n}}}\bar{\eta}_{a(\bar{n})}\sum_{n\in\mathcal{T}_{t_{\bar{n}}-1}}p_{n}=\lambda_{t_{\bar{n}}}\bar{\eta}_{a(\bar{n})}$.
\end{enumerate}
This completes the proof. $\hfill \Box$
\endproof
}

\section{Algorithm for Solving Risk-Averse Two-Stage Models}\label{e-companion:alg}
In this appendix, we provide an approximation algorithm (Algorithm \ref{alg:approx-twostage}) for solving risk-averse two-stage models \eqref{model:tworiskaverse}.
\begin{algorithm}[ht!]
\caption{Approximation Algorithm for Risk-Averse Two-Stage Facility Location \eqref{model:tworiskaverse}}
\begin{algorithmic}[1]
\label{alg:approx-twostage}
\STATE Solve the LP relaxation of the risk-averse multistage facility location problem \eqref{model:tworiskaverse} and let $(\boldsymbol{x}_n^{LP}, \eta_n^{LP}, \boldsymbol{y}_n^{LP}, u_n^{LP})_{n\in\mathcal{T}}$ be an optimal solution. If $\boldsymbol{x}_n^{LP}$ is integral for all $n\in\mathcal{T}$, stop and return $(\boldsymbol{x}_n^{LP}, \eta_n^{LP}, \boldsymbol{y}_n^{LP}, u_n^{LP})_{n\in\mathcal{T}}$.
\STATE Initialize $k=0$ and $(\boldsymbol{x}_n^0, \eta_n^0, \boldsymbol{y}_n^0, u_n^0)_{n\in\mathcal{T}} = (\boldsymbol{x}_n^{LP}, \eta_n^{LP}, \boldsymbol{y}_n^{LP}, u_n^{LP})_{n\in\mathcal{T}}$.
\WHILE{$||\boldsymbol{x}^k - \boldsymbol{x}^{k-1}|| \ge \epsilon,\ ||\eta^k - \eta^{k-1}|| \ge \epsilon,\ ||\boldsymbol{y}^k - \boldsymbol{y}^{k-1}|| \ge \epsilon,\ ||u^k - u^{k-1}|| \ge \epsilon$}
\STATE Solve Problem $\mbox{{\bf SP-RTS}}(\boldsymbol{y}_n^{k}, u_n^{k})$
and let $\boldsymbol{x}_n^{k+1},\ {\eta}_n^{k+1}$ denote the corresponding optimal solutions. We have the analytical form of the optimal solutions as $\boldsymbol{x}^{k+1}_{1}=\lceil \boldsymbol{B}_{t_1} \boldsymbol{y}_1^{k}\rceil,\ \boldsymbol{x}_n^{k+1} = \max_{m\in\mathcal{P}(n)}\lceil \max_{l\in \mathcal{T}_{t_m}}\boldsymbol{B}_{t_l} \boldsymbol{y}_l^{k}\rceil-\max_{m\in\mathcal{P}(a(n))}\lceil \max_{l\in \mathcal{T}_{t_m}}\boldsymbol{B}_{t_l} \boldsymbol{y}_l^{k}\rceil,\ \forall n\not=1$ and $\eta^{k+1}_n = \max_{m\in \mathcal{T}_{t_n+1}}\left\lbrace \boldsymbol{f}^{\mathsf T}_{t_m}\sum_{l\in\mathcal{P}(m)}\boldsymbol{x}^{k+1}_l+ \boldsymbol{c}_{t_m}^{\mathsf T} \boldsymbol{y}^k_m-u_m^{k}\right\rbrace, \ \forall n\not\in\mathcal{L}$.

\STATE Solve the following problem for each $n\in\mathcal{T},\ n\not=1$ independently
\begin{subequations}\label{eq:independent-ts}
\begin{align}
		\min_{{\boldsymbol{y}_n}, u_n}\quad &   \tilde{\boldsymbol{c}}_{n}^{\mathsf T}\boldsymbol{y}_n + \tilde{\alpha}_nu_n \nonumber\\
\text{s.t.}	\quad	&    \boldsymbol{B}_{t_n}\boldsymbol{y}_n \le \sum_{m\in \mathcal{P}(n)}{\boldsymbol{x}^{k+1}_{m}},\\
		& \boldsymbol{A}\boldsymbol{y}_n =\boldsymbol{d}_n\\
		& u_n - \boldsymbol{c}_{t_n}^{\mathsf T}\boldsymbol{y}_n \ge \boldsymbol{f}^{\mathsf T}_{t_n} \sum_{m\in \mathcal{P}(n)}\boldsymbol{x}^{k+1}_{m} - \eta^{k+1}_{a(n)},\label{eq:independent-u}\\
		& \boldsymbol{y}_n\in\mathbb{R}_{+},\ u_n\ge 0, \nonumber
\end{align}
\end{subequations}
and when $n=1$ we solve Problem \eqref{eq:independent-ts} without the variables $u_n$ and constraints \eqref{eq:independent-u}.
Let $\boldsymbol{y}_n^{k+1},\ {u}_n^{k+1}$ be the optimal solutions.
\STATE Update $k = k+1$.
\ENDWHILE
\STATE Return $(\boldsymbol{x}_n^H, \eta_n^H, \boldsymbol{y}_n^H, u_n^H)_{n\in\mathcal{T}} := (\boldsymbol{x}_n^{k+1}, \eta_n^{k+1}, \boldsymbol{y}_n^{k+1}, u_n^{k+1})_{n\in\mathcal{T}}$.
\end{algorithmic}
\end{algorithm}

{
\section{Examples to Illustrate the Tightness of the Lower Bounds}\label{e-companion:example}
In this appendix, we provide Example \ref{eg1} to illustrate the tightness of lower bound ${\rm VMS_R^{LB}}$ presented in Section \ref{sec:VMS}, and Example \ref{eg2} to illustrate the tightness of lower bound ${\rm VMS_P^{LB}}$ presented in Section \ref{sec:proof-VMS-prior}, respectively.
\begin{example}\label{eg1}
    Here, we consider two facilities and one customer site, where the customer site is closer to Facility \#1 having smaller capacity. Both facilities have the same investment costs, and the customer site's demand $d_1$ when $t=2$ is unknown, which equals to 50 in Scenario $\omega_1$, and 150 in Scenario $\omega_2$. In a two-stage model, to cover the demand in both scenarios, one must invest both facilities for $t=2$ up front. On the contrary, in a multistage model, the decision maker has the chance to wait and see the revealed demand and then make investment decisions, and as a result, the optimal solution will choose to invest Facility \#1 in Scenario $\omega_1$ and invest both facilities in Scenario $\omega_2$, which obviously costs less than the two-stage counterpart.
\begin{figure}[ht!]
    \centering
    \includegraphics[width=0.6\textwidth]{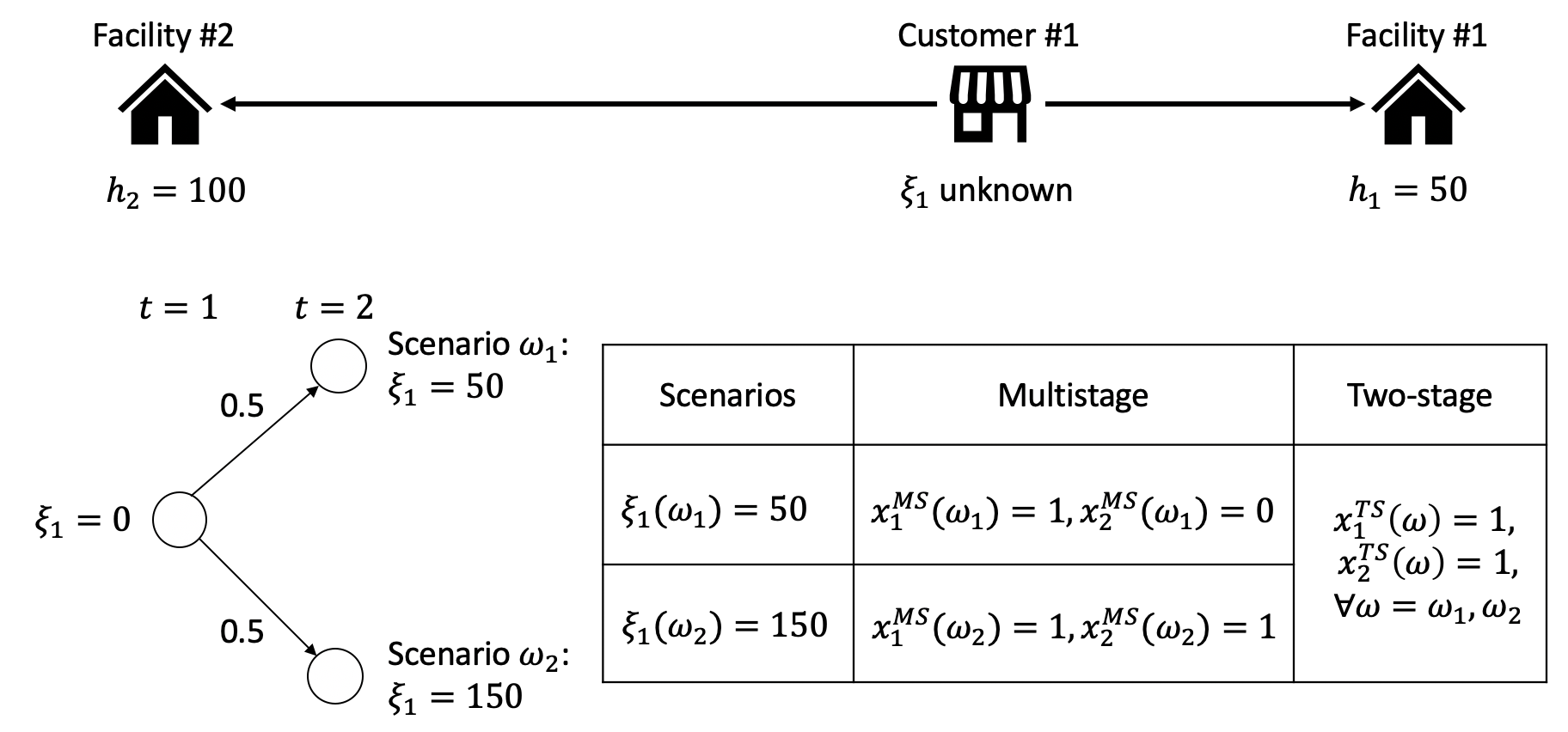}
    \caption{An instance to illustrate the gap between the optimal objective values of the multistage and two-stage facility location models.}
    \label{fig:eg}
\end{figure}

 Because $d_1=0$ when $t=1$, all costs will occur at $t=2$ and thus we only analyze the cost at $t=2$ and omit the index $t$ for all variables here. Denote $f$ as the investment cost of both facilities for one stage and $c_{11},\ c_{21} (c_{11}\le c_{21})$ as the operational cost from Facility \#1 and \#2 to Customer \#1, respectively. From previous discussions, in the two-stage model, the decision maker will choose to invest both facilities when $t=2$ up front. Because the customer site has different distances to the facilities, the optimal flow decisions at the second stage will be $y^{*}_{11} = 50,\ y^{*}_{21}=0$ in Scenario $\omega_1$, and $y^{*}_{11} = 50,\ y^{*}_{21}=100$ in Scenario $\omega_2$. It can be easily verified that $\boldsymbol{y}^{*}$ also constitute the flow decisions in an optimal solution to the multistage model. Next, we calculate the optimal objective values of the risk-averse two-stage and multistage models while assuming $\alpha=0.5$ as follows:
    \begin{itemize}
        \item Two-stage model \eqref{model:tworiskaverse}:
        \begin{itemize}
            \item Scenario $\omega_1$: $g_2(\omega_1)=\boldsymbol{f}^{\mathsf T}\boldsymbol{x}^{*}(\omega_1)+\boldsymbol{c}^{\mathsf T}\boldsymbol{y}^{*}(\omega_1) = 2f + 50c_{11}$ with probability 0.5;
            \item Scenario $\omega_2$: $g_2(\omega_2)=\boldsymbol{f}^{\mathsf T}\boldsymbol{x}^{*}(\omega_2)+\boldsymbol{c}^{\mathsf T}\boldsymbol{y}^{*}(\omega_2) = 2f + 50c_{11} + 100c_{21}$ with probability 0.5.
        \end{itemize}
        Because $\alpha = 0.5$, we have ${\rm VaR_{\alpha}}[g_2] = 2f + 50c_{11},\ {\rm CVaR_{\alpha}}[g_2] = {\rm VaR_{\alpha}}[g_2] + \frac{1}{1-\alpha}\mathbb{E}[g_2 - {\rm VaR_{\alpha}}[g_2]]_+ = 2f + 50c_{11} + \frac{1}{0.5}[0.5 * 100c_{21}] = 2f + 50c_{11} + 100c_{21}$ and correspondingly $u^{*}(\omega_1)=0,\ u^{*}(\omega_2)=100c_{21}$. Therefore, $z_R^{TS} = \rho_2(g_2) = (1-\lambda)\mathbb{E}[g_2] + \lambda {\rm CVaR}(g_2) = (1-\lambda)\frac{4f+100c_{11}+100c_{21}}{2} + \lambda(2f + 50c_{11} + 100c_{21})$.
        \item Multistage model \eqref{model:multiriskaverse}:
            \begin{itemize}
            \item Scenario $\omega_1$: $g_2(\omega_1)=\boldsymbol{f}^{\mathsf T}\boldsymbol{x}^{MS}(\omega_1)+\boldsymbol{c}^{\mathsf T}\boldsymbol{y}^{MS}(\omega_1) = f + 50c_{11}$ with probability 0.5;
            \item Scenario $\omega_2$: $g_2(\omega_2)=\boldsymbol{f}^{\mathsf T}\boldsymbol{x}^{MS}(\omega_2)+\boldsymbol{c}^{\mathsf T}\boldsymbol{y}^{MS}(\omega_2) = 2f + 50c_{11} + 100c_{21}$ with probability 0.5.
        \end{itemize}
        Because $\alpha = 0.5$, we have ${\rm VaR_{\alpha}}[g_2] = f + 50c_{11}$ and ${\rm CVaR_{\alpha}}[g_2] = {\rm VaR_{\alpha}}[g_2] + \frac{1}{1-\alpha}\mathbb{E}[g_2 - {\rm VaR_{\alpha}}[g_2]]_+ = f + 50c_{11} + \frac{1}{0.5}[0.5 * (f + 100c_{21})] = 2f + 50c_{11} + 100c_{21}$. Therefore, $z_R^{MS} = \rho_2(g_2) = (1-\lambda)\mathbb{E}[g_2] + \lambda {\rm CVaR}(g_2) = (1-\lambda)\frac{3f+100c_{11}+100c_{21}}{2} + \lambda(2f + 50c_{11} + 100c_{21})$.

    \end{itemize}
    
    As a result, the gap between the optimal objective values of the two-stage and multistage models only exists in the investment cost, which is ${\rm VMS_R} = z_R^{TS} - z_R^{MS} = 0.5(1-\lambda)f$. Based on equations \eqref{eq:construct_eta_ts} and \eqref{eq:construct_eta_ms}, the constructed $\eta$-solutions for the root node are $\eta^{TS}=\max\{2f + 50c_{11}-0, 2f + 50c_{11} + 100c_{21} - 100c_{21}\} = 2f + 50c_{11}$ and $\eta^{MS}=\max\{f + 50c_{11} - 0, 2f + 50c_{11} + 100c_{21} - 100c_{21}\} = 2f + 50c_{11} = \eta^{TS}$ and thus the lower bound ${\rm VMS_R^{LB}}$ only depends on the variation of utilization rates $\boldsymbol{B} \boldsymbol{y}^{*}$ across different scenarios. We notice that Facility \#1 is always fully utilized under both scenarios, and thus it does not contribute to ${\rm VMS_R^{LB}}$. On the other hand, Facility \#2 has utilization rate $\boldsymbol{B} \boldsymbol{y}^{*} = 0$ under Scenario $\omega_1$, and $\boldsymbol{B} \boldsymbol{y}^{*} = 1$ under Scenario $\omega_2$, which leads to ${\rm VMS_R^{LB}} = 0.5 * (1-\lambda) * f * (1 - 0) =0.5(1-\lambda)f = {\rm VMS_R}$. This example also provides a special case where both ${\rm VMS_R^{LB}}$ and ${\rm VMS_R}$ are positively related to the investment cost and negatively impacted by the risk attitude $\lambda$.
    \end{example}
    
    \begin{example}\label{eg2}
    Now let us consider three facilities and one customer site as shown in Figure \ref{fig:eg2}, where Facility \#1 is the closest to the customer site and Facility \#3 is the farthest. All facilities have the same capacity $h=100$, and the customer site's demand $d_t$ is uncertain. We consider $T=2$ stages and plot the possible realizations of $(d_t,F_t)_{t=1}^2$ on the scenario tree in Figure \ref{fig:eg2}. The confidence level of $\rm{CVaR}$ is set to $\alpha_t=0.95,\ \forall t=1,2$. In a two-stage model, the optimal priority list is $s_{t,12}^{*}=s_{t,13}^{*}=1,\  s_{t,23}^{*}=1,\ \forall t=0,1$ while all other $\boldsymbol{s}^{*}$-values are 0, the optimal opening decision is $x^*_{1,1}(\omega)=1,\ x^*_{2,2}(\omega)=1,\ \forall \omega\in\{\omega_1,\omega_2,\omega_3,\omega_4\},\ x^*_{2,3}(\omega_2)=x^*_{2,3}(\omega_3)=x^*_{2,3}(\omega_4)=1$, and the optimal flow decision is $\boldsymbol{c}_1^{\mathsf T}\boldsymbol{y}_1^{*}(\omega_1)=\boldsymbol{c}_1^{\mathsf T}\boldsymbol{y}_1^{*}(\omega_2)=50c_{11},\ \boldsymbol{c}_1^{\mathsf T}\boldsymbol{y}_1^{*}(\omega_3)=\boldsymbol{c}_1^{\mathsf T}\boldsymbol{y}_1^{*}(\omega_4)=100c_{11},\ \boldsymbol{c}_2^{\mathsf T}\boldsymbol{y}_2^{*}(\omega_1)=\boldsymbol{c}_2^{\mathsf T}\boldsymbol{y}_2^{*}(\omega_3)=100c_{11}+100c_{21},\ \boldsymbol{c}_2^{\mathsf T}\boldsymbol{y}_2^{*}(\omega_2)=\boldsymbol{c}_2^{\mathsf T}\boldsymbol{y}_2^{*}(\omega_4)=100c_{11}+100c_{21}+100c_{31}$. Because $\alpha_t=0.95$ and we only have four scenarios, the optimal $\boldsymbol{\eta}$-values are the maximum cost in each stage, i.e., $\eta_1^*(\omega)={\rm VaR}_{\alpha_1}[\boldsymbol{c}_1^{\mathsf T}\boldsymbol{y}^*_1(\omega)+\sum_{i_1\not=i_2}s_{1,i_1i_2}^*(\omega):\omega\in\{\omega_1,\omega_2,\omega_3,\omega_4\}]=100c_{11}+3,\ \eta_2^*(\omega)={\rm VaR}_{\alpha_2}[\boldsymbol{c}_2^{\mathsf T}\boldsymbol{y}^*_2(\omega):\omega\in\{\omega_1,\omega_2,\omega_3,\omega_4\}]=100c_{11}+100c_{21}+100c_{31},\ \forall \omega\in\{\omega_1,\omega_2,\omega_3,\omega_4\}$. Correspondingly, the optimal $\boldsymbol{u}$-values are $u_t^{*}(\omega)=0,\ \forall t=1,2,\ \omega\in\{\omega_1,\omega_2,\omega_3,\omega_4\}$. It can be easily verified that $(\boldsymbol{x}^*,\boldsymbol{y}^*,\boldsymbol{u}^*)$ also constitutes an optimal solution to the multistage model, and we use $(\boldsymbol{s}^{\prime},\boldsymbol{\eta}^{\prime},\boldsymbol{x}^{\prime},\boldsymbol{y}^{\prime},\boldsymbol{u}^{\prime})$ to denote the optimal solution to the multistage model. Then the optimal priority list in the multistage model is $s_{0,12}^{\prime}=s_{0,13}^{\prime}=1,\ s_{0,23}^{\prime}=1,\ s_{1,23}^{\prime}=1$ while all the other $\boldsymbol{s}^{\prime}$-values are 0, and the optimal $\boldsymbol{\eta}$-values are $\eta_1^{\prime}(\omega)={\rm VaR}_{\alpha_1}[\boldsymbol{c}_1^{\mathsf T}\boldsymbol{y}^{\prime}_1(\omega)+\sum_{i_1\not=i_2}s_{1,i_1i_2}^{\prime}(\omega):\omega\in\{\omega_1,\omega_2,\omega_3,\omega_4\}]=100c_{11}+1,\ \eta_2^{\prime}(\omega_1)=\eta_2^{\prime}(\omega_2)={\rm VaR}_{\alpha_2}[\boldsymbol{c}_2^{\mathsf T}\boldsymbol{y}^{\prime}_2(\omega):\omega\in\{\omega_1,\omega_2\}]=100c_{11}+100c_{21}+100c_{31},\ \eta_2^{\prime}(\omega_3)=\eta_2^{\prime}(\omega_4)={\rm VaR}_{\alpha_2}[\boldsymbol{c}_2^{\mathsf T}\boldsymbol{y}^{\prime}_2(\omega):\omega\in\{\omega_3,\omega_4\}]=100c_{11}+100c_{21}+100c_{31}$.
            \begin{figure}[ht!]
    \centering
    \includegraphics[width=0.4\textwidth]{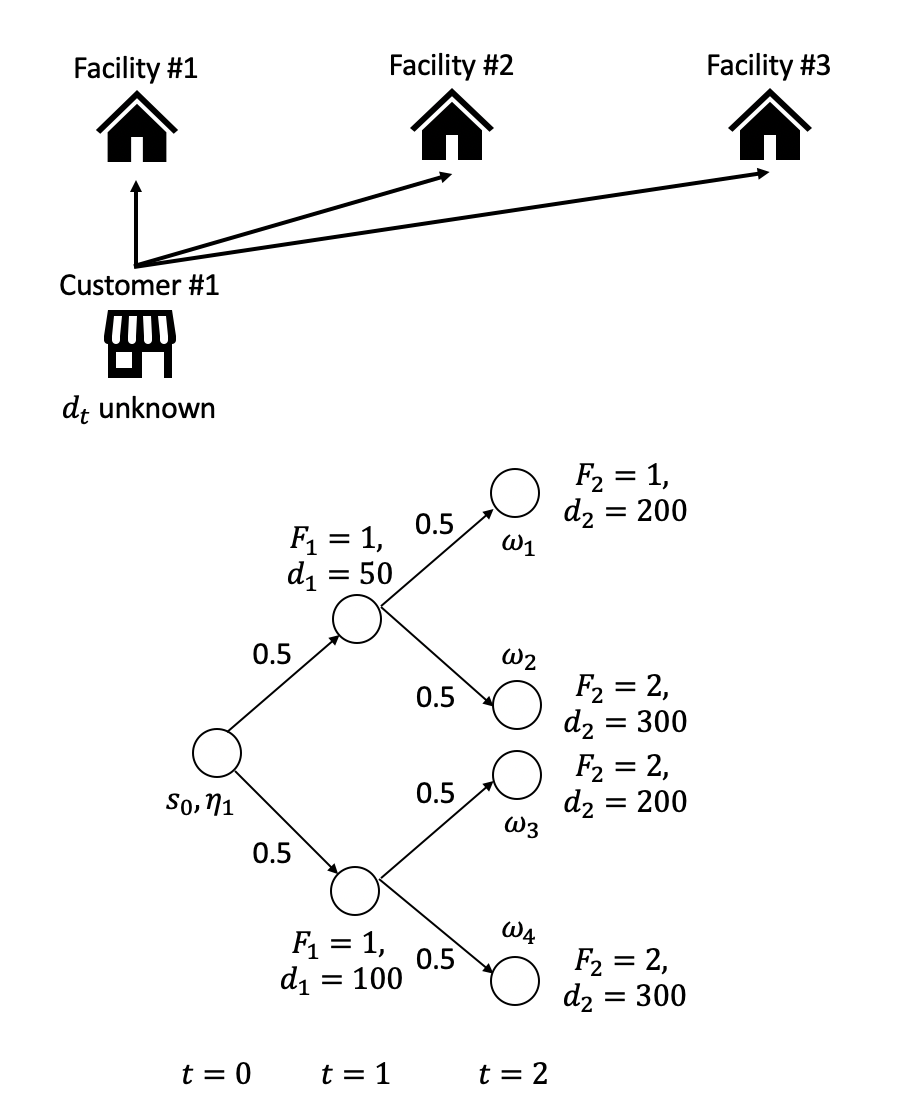}
    \caption{An instance to illustrate the gap between the optimal objective values of the multistage and two-stage facility location models with prioritization.}
    \label{fig:eg2}
\end{figure}

Comparing the optimal solutions to the two models, the only difference is $\boldsymbol{s}^{\prime}_1\not=\boldsymbol{s}^*_1,\ \eta_1^{\prime}\not=\eta_1^*$, and thus ${\rm VMS_P}=z_P^{TS}-z_P^{MS}=\sum_{\omega\in\Omega}0.25\left((1-\lambda)(\sum_{i_1\not=i_2}\boldsymbol{s}_1^*(\omega)-\sum_{i_1\not=i_2}\boldsymbol{s}_1^{\prime}(\omega))+\lambda(\eta_1^*(\omega)-\eta_1^{\prime}(\omega))\right)=2(1-\lambda)+2\lambda=2$.
    Based on Proposition \ref{prop:substructure-prior}, the constructed $\boldsymbol{\eta}^{TS},\ \boldsymbol{s}^{TS}$ coincide with the optimal solutions to the two-stage model (i.e., $\boldsymbol{\eta}^{*},\ \boldsymbol{s}^{*}$), and the constructed $\boldsymbol{\eta}^{MS},\ \boldsymbol{s}^{MS}$ coincide with the optimal solutions to the multistage model (i.e., $\boldsymbol{\eta}^{\prime},\ \boldsymbol{s}^{\prime}$). According to \eqref{eq:VMS_P_LB}, the analytical lower bound can be calculated by ${\rm VMS_P^{LB}}=\sum_{\omega\in\Omega}0.25\left(\lambda(\eta_1^{TS}(\omega)-\eta_1^{MS}(\omega))+\lambda(\eta_2^{TS}(\omega)-\eta_2^{MS}(\omega))+(1-\lambda)(\sum_{i_1\not=i_2}\boldsymbol{s}_1^{TS}(\omega)-\sum_{i_1\not=i_2}\boldsymbol{s}_1^{MS}(\omega))\right)=2\lambda+2(1-\lambda)=2={\rm VMS_P}$.
    \end{example}
    
\section{Two-Stage Distributionally Robust Facility Location}
When the distribution of the uncertain demand is unknown, one can resort to a distributionally robust approach by constructing an ambiguity set and selecting the optimal solution against the worst-case distribution within the set. We follow \citet{xie2019tractable} and use the $\infty-$Wasserstein distance to define the ambiguity set. As a result, the two-stage distributionally robust model becomes
\begin{align*}
    z_{DRO}^{TS} = \min_{\boldsymbol{x},\boldsymbol{y}}\sum_{t=1}^T \boldsymbol{f}^{\mathsf T}_{t} \sum_{\tau=1}^t \boldsymbol{x}_{\tau}+\max_{P\in\mathcal{P}}\mathbb{E}_P[\sum_{t=1}^TQ_t(\boldsymbol{x},\boldsymbol{d})],
\end{align*}
where $\mathcal{P}=\{\mathbb{P}:\mathbb{P}\{\boldsymbol{\xi}\in\Xi\}=1, W^{\infty}(\mathbb{P},\hat{\mathbb{P}}_K)\le \epsilon\}$ and 
\begin{align*}
    W^{\infty}(\mathbb{P}_1,\mathbb{P}_2)=\inf_{\mathbb{Q}}\{\text{ess.sup}||\boldsymbol{\xi}_1-\boldsymbol{\xi}_2||_p\mathbb{Q}(d\boldsymbol{\xi}_1,d\boldsymbol{\xi}_2): \mathbb{Q} \text{ is a joint distribution of } \tilde{\boldsymbol{\xi}}_1 \text{ and } \tilde{\boldsymbol{\xi}}_2 \text{with marginals } \mathbb{P}_1, \mathbb{P}_2\}.
\end{align*}
Here, ess.sup denotes the essential supremum, norm $||\cdot||_p$ denotes the reference distance with $p\in[1,\infty]$ and $\hat{\mathbb{P}}_K$ denotes a discrete empirical distribution generated by $K$ i.i.d. samples, i.e., $\hat{\mathbb{P}}_K\{\boldsymbol{d}=\boldsymbol{d}^j\}=\frac{1}{K},\ \forall j\in[K]$. 

In the following theorem, we show that the inner worst-case expectation cost has a tractable representation when the reference distance $p=1$.
\begin{theorem}\label{thm:DRO}
Suppose that $p=1$. Then the inner worst-case expectation cost $\max_{P\in\mathcal{P}}\mathbb{E}_P[\sum_{t=1}^TQ_t(\boldsymbol{x},\boldsymbol{d})]$ is equivalent to
\begin{align*}
    \frac{1}{K}\sum_{k=1}^K\max_{r\in\{1,-1\}}\max_{\hat{t}\in[T],\hat{j}\in[N]}\min_{\boldsymbol{\beta}}\quad&\sum_{t=1}^T\sum_{i=1}^M\sum_{j=1}^Nc_{tij}\beta_{tij}\\
    \text{s.t.}\quad&\sum_{i=1}^M\beta_{tij}= d_{tj}^k+\epsilon r, \ \forall t, j = \hat{t}, \hat{j}\\
    & \sum_{i=1}^M\beta_{tij} = d_{tj}^k,\ \forall t \not = \hat{t} \text{ or } j \not=\hat{j}\\
    & \sum_{j=1}^N\beta_{tij} \le h_{ti}\sum_{\tau=1}^tx_{\tau i},\ \forall t\in[T],\ i\in[M]\\
    & \beta_{tij} \ge 0,\ \forall t\in[T],\ i\in[M],\ j\in[N].
\end{align*}
\end{theorem}
\proof{Proof of Theorem \ref{thm:DRO}}
Following the proof of Theorem 3 in \citet{xie2019tractable}, we have
\begin{align*}
    \max_{P\in\mathcal{P}}\mathbb{E}_P[\sum_{t=1}^TQ_t(\boldsymbol{x},\boldsymbol{d})] =& \frac{1}{K}\sum_{k=1}^K\sup\{\sum_{t=1}^TQ(\boldsymbol{x},\boldsymbol{d}): ||\boldsymbol{d}-\boldsymbol{d}^k||_1\le \epsilon\}\\
    \overset{(a)}{=}&\frac{1}{K}\sum_{k=1}^K\sup\{\sum_{t=1}^T\sum_{j=1}^N\gamma_{tj}d_{tj}+\sum_{t=1}^T\sum_{i=1}^M\mu_{tj}(h_{ti}\sum_{\tau=1}^tx_{\tau i}):\\ 
    &\hspace{2cm}\gamma_{tj}+\mu_{tj}\le c_{tij},\ \forall t\in[T],\ i\in[M],\ j\in[N]\\
    &\hspace{2cm}\mu_{tj}\le 0,\ \forall t\in[T],\ j\in[J]\\
    &\hspace{2cm}||\boldsymbol{d}-\boldsymbol{d}^k||_1\le \epsilon\}\\
    \overset{(b)}{=}& \frac{1}{K}\sum_{k=1}^K\sup\{\sum_{t=1}^T\sum_{j=1}^N\gamma_{tj}d^k_{tj}+\epsilon||\gamma||_{\infty}+\sum_{t=1}^T\sum_{i=1}^M\mu_{tj}(h_{ti}\sum_{\tau=1}^tx_{\tau i}):\\ 
    &\hspace{2cm}\gamma_{tj}+\mu_{tj}\le c_{tij},\ \forall t\in[T],\ i\in[M],\ j\in[N]\\
    &\hspace{2cm}\mu_{tj}\le 0,\ \forall t\in[T],\ j\in[J]\}\\
    =&\frac{1}{K}\sum_{k=1}^K\sup\{\sum_{t=1}^T\sum_{j=1}^N\gamma_{tj}d^k_{tj}+\epsilon\max_{r\in\{-1,1\}}\max_{\hat{t}\in[T],\hat{i}\in[M]}r\gamma_{\hat{t}\hat{i}}+\sum_{t=1}^T\sum_{i=1}^M\mu_{tj}(h_{ti}\sum_{\tau=1}^tx_{\tau i}):\\ 
    &\hspace{2cm}\gamma_{tj}+\mu_{tj}\le c_{tij},\ \forall t\in[T],\ i\in[M],\ j\in[N]\\
    &\hspace{2cm}\mu_{tj}\le 0,\ \forall t\in[T],\ j\in[J]\}\\
    \overset{(c)}{=}&\frac{1}{K}\sum_{k=1}^K\max_{r\in\{-1,1\}}\max_{\hat{t}\in[T],\hat{i}\in[M]}\min_{\beta}\sum_{t=1}^T\sum_{i=1}^M\sum_{j=1}^Nc_{tij}\beta_{tij}\\
    &\hspace{4.5cm}\sum_{i=1}^M\beta_{tij}=d_{tj}^k+\epsilon r,\ \forall t,j=\hat{t},\hat{j}\\
    &\hspace{4.5cm}\sum_{i=1}^M\beta_{tij}=d_{tj}^k,\ \forall t\not=\hat{t},\text{ or } j \not=\hat{j}\\
    &\hspace{4.5cm} \sum_{j=1}^N\beta_{tij} \le h_{ti}\sum_{\tau=1}^tx_{\tau i},\ \forall t\in[T],\ i\in[M]\\
    &\hspace{4.5cm} \beta_{tij} \ge 0,\ \forall t\in[T],\ i\in[M],\ j\in[N].
\end{align*}
Here, $(a)$ and $(c)$ are true because of strong duality, $(b)$ is true because of the definition of dual norm, i.e., $||r||_{p^*}=\max_{s}\{r^{\mathsf T}s: ||s||_p\le 1\}$.
\endproof

\begin{corollary}\label{cor:DRO}
Suppose $p=1$. Then the two-stage distributionally robust facility location model is equivalent to
\begin{align*}
    z_{DRO}^{TS} = \min_{\boldsymbol{x},\boldsymbol{\eta},\boldsymbol{\beta}} \quad& \sum_{t=1}^T \boldsymbol{f}^{\mathsf T}_{t} \sum_{\tau=1}^t \boldsymbol{x}_{\tau}+\frac{1}{K}\sum_{k=1}^K\eta_k\\
    \text{s.t.}\quad & \eta_k\ge \sum_{t=1}^T\sum_{i=1}^M\sum_{j=1}^Nc_{tij}\beta_{tij}^{\hat{t}\hat{j}rk},\ \forall \hat{t}\in[T],\ \hat{j}\in[N],\ r\in\{-1,1\},\ k\in[K]\\
    & \sum_{i=1}^M\beta_{\hat{t}i\hat{j}}^{\hat{t}\hat{j}rk} = d_{\hat{t}\hat{j}}^k+\epsilon r, \ \forall \hat{t}\in[T],\ \hat{j}\in[N],\ r\in\{-1,1\},\ k\in[K]\\
    & \sum_{i=1}^M\beta_{tij}^{\hat{t}\hat{j}rk} = d_{tj}^k,\ \forall t \not = \hat{t}\ \mbox{or}\ j \not=\hat{j},\ \hat{t}\in[T],\ \hat{j}\in[N],\ r\in\{-1,1\},\ k\in[K]\\
    & \sum_{j=1}^N\beta_{tij}^{\hat{t}\hat{j}rk} \le h_{ti}\sum_{\tau=1}^tx_{\tau i},\ \forall t\in[T],\ i\in[M],\ \hat{t}\in[T],\ \hat{j}\in[N],\ r\in\{-1,1\},\ k\in[K]\\
    & \sum_{\tau=1}^tx_{\tau i}\le 1,\ \forall t \in [T], \ i\in[M]\\
    & \beta_{tij}^{\hat{t}\hat{j}rk} \ge 0,\ \forall t\in[T],\ i\in[M],\ j\in[N],\ \hat{t}\in[T],\ \hat{j}\in[N],\ r\in\{-1,1\},\ k\in[K]\\
    & x_{ti}\in\{0,1\},\ \forall t \in [T], \ i\in[M].
\end{align*}
\end{corollary}
\proof{Proof of Corollary \ref{cor:DRO}}
Follows immediately from Theorem \ref{thm:DRO}.
\endproof
}

\end{appendices}
%
%
%






\end{document}